\documentclass[12pt]{article}
\usepackage{latexsym}
\usepackage[all]{xy}
\usepackage{amsmath,amssymb}
\usepackage{amsfonts}

\newtheorem{thm}{Theorem}[section]
\newtheorem{prop}[thm]{Proposition}

\newtheorem{df}[thm]{Definition}
\newtheorem{cor}[thm]{Corollary}
\newtheorem{conj}[thm]{Conjecture}
\newtheorem{rmk}[thm]{Remark}

\begin{document}

\title{Homotopical and Higher Categorical Structures in Algebraic Geometry}
\author{\textbf{Bertrand To\"en}}\date{}

\maketitle

\vspace{5cm}

\begin{center}
CNRS UMR 6621

Universit\'e de Nice, Math\'ematiques

Parc Valrose

06108 Nice Cedex 2

France

\end{center}

\vspace{2cm}

\begin{center}

\textbf{M\'emoire de Synth\`ese}

en vue de l'obtention de l'

\textbf{Habilitation a Diriger les Recherches}

\end{center}

\newpage

\vspace*{8cm}

\textit{Je d\'edie ce m\'emoire \`a Beussa la mignonne.}

\newpage

\vspace*{4cm}

\begin{center}

\textit{Pour ce qui est des choses humaines, ne pas rire, ne pas pleurer, ne pas
s'indigner, mais comprendre.}

\bigskip

Spinoza

\end{center}

\vspace*{4cm}

\begin{center}

\textit{Ce dont on ne peut parler, il faut le taire.}

\bigskip

Wittgenstein

\end{center}

\vspace*{4cm}

\begin{center}

\textit{Ce que ces gens l\`a sont cryptiques !}

\bigskip

Simpson

\end{center}

\newpage

\vspace*{3cm}

\begin{center} \textit{Avertissements} \end{center}

\bigskip

La pr\'esent m\'emoire est une version
\'etendue du m\'emoire de synth\`ese
d'habilitation originellement r\'edig\'e
pour la soutenance du 16 Mai 2003. Le
chapitre \S 5 n'apparaissait pas
dans la version originale, bien que son
contenu \'etait pr\'esent\'e sous forme
fortement r\'esum\'e.

\newpage

\tableofcontents

\newpage

\begin{center} \textit{Remerciements} \end{center}

\bigskip

Par le hazard des choses, la soutenance de cette habilitation coincide avec mon
d\'epart du laboratoire J.A. Dieudonn\'e de Nice, ce qui rend \`a mes yeux cette page
de remerciements doublement importante. \\

Mes tous premiers remerciements vont \`a Carlos Simpson, pour de tr\`es
nombreuses conversations et correspondences math\'ematiques au cours
des quelles j'ai immens\'ement appris, mais aussi pour ses
conseils et suggestions lors de mes questionements vari\'es. Il ne fait aucun doute
pour moi que son influence a \'et\'e d\'eterminante pour
beaucoup de mes orientations de recherches.

Je remercie tout particuli\`erement Andr\'e Hirschowitz, qui au cours
de ces deux ann\'ees pass\'ees \`a Nice m'a fait partager son enthousiame
l\'egendaire. Les heures pass\'ees dans le bureau W703 \`a d\'ecortiquer
et r\'einventer des notions math\'ematiques vari\'ees (qui pourtant
ne nous avaient rien fait !) furent pour moi un r\'eel
plaisir. Il aura marqu\'e de fa\c{c}on d\'efinitive mon s\'ejour
au JAD. Je voudrais au passage profiter de ces quelques lignes
pour lui d\'edier le fameux
proverbe Bantou suivant:
\textit{Ne cite jamais Alexandre Grothendieck. Si tu ne sais pas pourquoi,
Andr\'e Hirschowitz le sait certainement.}\\

C'est un honneur pour moi de remercier Maxim Kontsevich et Andr\'e Joyal
pour avoir accept\'e d'\^etre rapporteurs de cette habilitation. Leurs
remarques, conseils et encouragements me furent tr\`es utiles.

Je remercie Arnaud Beauville, Yan Soibelman et Vadim Schechtman
pour avoir accept\'e de faire partie du jury. \\

Ma gratitute va aussi \`a mes co-auteurs. Gabriele Vezzosi, qui en plus
d'\^etre un collaborateur d'exception est aussi un ami. Ludmil Katzarkov et
Tony Pantev pour leurs encouragements au tout d\'ebut de mon
int\'eret pour la th\'eorie de Hodge non-ab\'elienne, et pour leur
dynamisme scientifique. Et enfin, Joachim Kock, collaborateur mais aussi
compagnon jazzistique
durant ces ann\'ees \`a Nice. \\

Je remercie enfin les membres du JAD pour l'atmosph\`ere toute
particuli\`ere qui plane dans ce laboratoire. Et plus particuli\`erement,
Olivier Penacchio, Mathieu Bernard, Alessandro Chiodo, Marco Maggesi, Angela Ortega, Marie Virat, Fabien Herbaut \dots
pour repas, caf\'es, pizzas, parties de Billards, bi\`eres et verres de vins \dots (je suis aujourd'hui
convaincu qu'il est illusoire d'essayer de faire de la
recherche math\'ematique sans au moins quelques de ces activit\'es).

\newpage

\setcounter{section}{-1}

\section{Introduction}

During my PhD. thesis I have been working on
cohomology and K-theory of algebraic stacks.
Since then, my research interests moved to
different questions concerning homotopical and
higher categorical structures related to algebraic geometry.
These questions have been originally motivated on one side
by the work of C. Simpson on non-abelian Hodge theory
and on the other side by the reading of
A. Grothendieck's \textit{Pursuing stacks}. From them, I have
learned that homotopical and higher categorical methods
can be used in order to define and study interesting
invariants of algebraic varieties, and it is the main purpose of this short
m\'emoire to present my recent researches
in that direction (several of these works are joint works).
I make my apologizes to the reader for the rather informal
point of view I have decided to follow, first of all because
some
of these works have not been written up and will be presented without proofs
(that hopefully will appear elsewhere), but also because even when written up
only few of them have been published. In any case, I have tried to give
complete and precise definitions, statements and/or references, so that the reader can at least
have a rough idea on the state of my research activities at the present time.
I also apologize as several of these works are still in progress, and
therefore the results I present in this text are probably not the most
definitive: in my opinion the interesting work still has to be done.
And finally, I apologize for my poor English (this
was my last apologizes).

All along this work, I have also tried
to show that the results of this m\'emoire are not at all
independent from each others and that they all belong to the realm of
\textit{homotopical mathematics}\footnote{I have learned this
expression from M. Kontsevich.}. Mathematics are based on sets theory and the notion
of structures (in Bourbaki's sense), while in homotopical mathematics
sets are replaced by homotopy types, and structures are then enriched
over the homotopy theory of spaces (e.g. groups are then replaced by
simplicial groups, categories by $S$-categories, presheaves by simplicial
presheaves, algebras by dga's \dots) . The general philosophy (which is probably
quite old and I guess
goes back to Boardman, Dwyer, Kan, Quillen, Thomason, Waldhausen, Vogt, \dots),
seems to be that a huge part of mathematics does possess interesting and useful
extensions to the context of homotopical mathematics. Of course, several examples
of extensions of fundamental notions have already been studied by many authors,
as $S$-categories, simplicial presheaves, $A_{\infty}$-algebras and $A_{\infty}$-categories \dots.
I like to consider the results of the present work as part of
the possible extensions of algebraic geometry to the realm of homotopical mathematics,
what we call together with Gabriele \textit{homotopical algebraic geometry}.
To be a bit more explicit let me give the following board in which I express
various notions appearing in this work as extensions of
well known notions to the homotopical mathematics context. I apologize for
the rather artificial flavor of this comparison, but I hope
it helps anyway to give a unity to the results of this m\'emoire. \\

$$\hspace{-1cm}
\begin{array}{l|l}
\mathbf{Mathematics} & \mathbf{Homotopical \; Mathematics} \\
 & \\
& \\
Sets & Simplicial\; sets,\; homotopy \;  types \\
& \\
Categories & Segal \; categories,\; or\;  model \; categories \\
& \\
Categories \; of \; functors \; \underline{Hom}(A,B) &
Segal \; categories \; \mathbb{R}\underline{Hom}(A,B) \\
& \\
The\; category\; of \; sets,\; Sets & The\; Segal\; category,\; or\;   model\; category, \\
& of\; simplicial\; sets,\; Top \\
& \\
Abelian \; categories & Stable \; Segal \; categories\; (Def. \ref{d1}), \\
&  or\; stable\; model \; categories \\
& \\
Quillen's\; K-theory\; of \; abelian \; categories &
K-theory\; of\; stable \; Segal \; categories ,\; (\S 2.1) \\
& or\; Waldhausen's \; K-theory \; of\; model\; categories \\
& \\
Presheaves,\; sheaves & Prestacks, \; stacks \; (\S 3.1) \\
& \\
Topos & Segal \; topos,\; model \; topos \; (\S 3.1) \\
& \\
Grothendieck's \; \Pi_{1} \; of \; a \; topos & Pro-homotopy \; type
\; of \; a \; Segal \; topos \; (Def. \ref{d7}) \\
& \\
Algebraic \; completion \; \Gamma^{alg} \; of \; a \; group \;
\Gamma & Schematization \; (X\otimes k)^{sch} \; of
\; a \; space \; X \; (Thm. \ref{t9}) \\
& \\
Affine \; gerbes & Segal \; affine \; gerbes \; (Def. \ref{dtan7})\\
& or\; Schematic\; homotopy \; types \; (Def. \ref{d8})\\
& \\
Tannakian \; category & Tannakian\; Segal\; category \; (Def. \ref{dtan7})\\
& \\
Algebraic \; geometry & Derived \; algebraic \; geometry \; (\S 4.2) \\
& \\
Algebraic\; geometry\; relative & Homotopical \; algebraic \; geometry \; over \\
to \; a \; base \; monoidal \; category \; C &
 \; a \;  base \; monoidal \; Segal\; category, \\
 & or\;  model \; category, \;
M \; (\S 4.1)
\end{array}$$

\bigskip

Before describing a short overview of the content of the present work
I would like to mention that there are several ways of
dealing with homotopy theory and therefore of
\textit{doing homotopical mathematics}, not all of them
being equivalent. I therefore had to make a choice. It is
probably not so easy to explain this choice (though I will try in section $1$), but I wish to mention that
it has been on one hand a pragmatic choice (i.e. I have chosen the theory
that seemed to me the best suited for my particular purpose) but also
a purely personal and psychological choice (i.e. I have chosen the theory I liked the most
and with which I felt the most comfortable). As the reader will see (and might have
already noticed in the previous comparison board), I have decided
to work with essentially two theories: model categories and Segal categories.
No doubts that one could also work with other theories. No doubts also that
one could not work with some of them. In some sense, one of the main slogan
of this work is: \textit{the combination model categories/Segal categories
is wonderful}.  \\

\bigskip

This work is organized in 6 sections, following a rather arbitrary splitting
into different themes, and which does not reflect at all
a chronological order.

In the very first section I present some general remarks about
homotopy theory, and in particular I try to compare various
theories and approaches. This part is supposed to explain
why one should not work with a unique theory but rather with
(at least) two at the same time. This part does not contain any
original mathematics.

In the second section, I present some results concerning Segal categories in their own.
These results (partially joint with
D.-C. Cisinski, A. Hirschowitz, J. Kock, C. Simpson and G. Vezzosi)
have been proved \textit{on the way}, in the sense that
there were not part of my main objective, and are sometimes quite unrelated to algebraic geometry.
However, from my point of view they are interesting examples showing the real flavor of
Segal categories, and how they can be useful. This section also includes a brief overview
of the relations between Segal categories and model categories, which will be used
all along this work.

Section $3$ is devoted to the notion
of \textit{Segal topos}, a natural extension
of the notion of topos to the Segal setting, and its application to an extension of
Artin-Mazur's style homotopy theory (this part is mainly a joint work with G. Vezzosi).
Originally, the notion of Segal topos appeared for the need of our work
on \textit{homotopical algebraic geometry}, but seems to us of independent
interests.

The fourth section is concerned with \textit{homotopical algebraic geometry}, a joint
work with G. Vezzosi, still in progress. I have included a short overview of the main
definitions. The main results are then some applications of the general formalism
to the construction of certain \textit{derived moduli spaces}, and to
what we like to call \textit{brave new algebraic geometry}, an algebraico-geometric point of view
on stable homotopy theory. Because the formalism of ``HAG'' is so widely ramified and
provide very natural settings for several other results of this m\'emoire, I personally
like to consider it as the deep heart of my recent works.

In section 5 I have included the present state of some
thoughts about a Segal version of Tannakian duality. This part
is unfortunately highly conjectural, but I have tried to
present precise definitions and conjectures which are all
mathematically meaningful. In fact the formalism of Tannakian
duality for Segal categories was my starting point of my
interests in homotopical and higher categorical structures, and
it provides a point of view that explains much more naturally
several constructions and statement (that can be proved
by complete different methods). It has been a great source of inspirations
for me during the last few years.

Section 6 is devoted to my work on homotopy types and Hodge theory.
The first part of this work is around the \emph{schematization problem}
of posed by Grothendieck in \cite{gr1}, for which I have proposed
a solution based on the notion of \emph{schematic homotopy types}.
In a second part, I present some joint work with T. Pantev and L. Katzarkov
in which we define a Hodge decomposition of the schematization
of a smooth projective complex manifold extending all the previously defined
Hodge structures on cohomology, rational homotopy and
fundamental group. I also propose a Tannakian interpretation, using
the notion of Tannakian Segal categories of the last
chapter, of theses constructions.

Finally, in an appendix I have reproduced a letter to P. May, containing
some general thoughts about higher categories as well as a short
note in which I construct comparison functors between
certain homotopy theories.

\bigskip

\bigskip

\begin{center} \textbf{Notations and conventions} \end{center}

First of all, I have decided to ignore universes considerations, and I
will assume implicitly that certain objects are \emph{small}
when required. I apologize for this choice but I hope it will help
keeping the redaction as light as possible.

I will use the books \cite{ho,hir} as references for model category theory.
For a model category $M$ I will use the notation $Map_{M}(x,y)$ for
the simplicial mapping space between two objects $x$ and $y$, as defined in \cite[\S 5.2]{ho}.
The homotopy category of $M$, i.e. the localization of $M$ along its
sub-category of equivalences, is denoted by $Ho(M)$.

For the notion and basic definitions of $S$-categories I refer to
\cite{dk1,dk2}. For an $S$-category $T$ I will use the notation
$\underline{Hom}_{T}(x,y)$ for the simplicial set of morphisms
between two objects $x$ and $y$ of $T$, and $Ho(T)$ for its homotopy category.
The simplicial localization
of a model category $M$ along its sub-category of equivalences $W$
is simply denoted by $LM:=L(M,W)$.

The references for the theory of Segal categories are \cite{hi-si,pe,to-ve1}.
For a Segal category $A$ I will denote by $A_{(x,y)}$ the simplicial set
of morphisms between two objets $x$ and $y$, and by $Ho(A)$ the homotopy category
of $A$. For two Segal categories $A$ and $B$ I will denote
by $\mathbb{R}\underline{Hom}(A,B)$ the derived internal $Hom$ between
$A$ and $B$ (see for example \cite[Cor. 2.0.5]{to-ve1}).
By definition, $Top$ is the Segal category of simplicial sets $Top:=LSSet$.

We will avoid using the expression \emph{functor} in the context of
Segal categories, and will reserve this expression for usual
category theory. Instead we will speak about
\emph{morphisms between Segal categories}, and such morphisms
will implicitly be assumed to live in some $\mathbb{R}\underline{Hom}$,
or in other words that we have performed some fibrant replacement.

\newpage

\section{Homotopy theories}

An abstract homotopy theory consists of a category of objects $C$, together
with a sub-category $W$ of equivalences and through which
objects of $C$ have to be identified. There exist several approaches
to deal with such a couple $(C,W)$, and the purpose of this section is
to present some of them and to explain their advantages as well as their limits.
What I would like to stress out is that it is sometimes necessary to use
more than one approach at the same time in order to have a
complete understanding of the situation. More precisely I would like
to explain through some pertinent examples that
the combination of model category theory and Segal category theory
is probably one of the most efficient way of doing homotopy theory. I hope this
will justify my choice of using model categories and Segal categories
as \textit{homotopical mathematics analogs} of categories.

\begin{center} \textit{Localization} \end{center}

The very first approach to abstract homotopy theory
is the localization construction which associates
to a pair $(C,W)$ a localized category $W^{-1}C$, obtained from
the category $C$ by formally inverting the morphisms of $W$ (see \cite{ga-zi}).
This construction has been of fundamental importance, for examples in order
to define the homotopy category of spaces as well as derived categories.
It is however very brutal, and
much too coarse in many contexts. A striking example is
when $C=Cat$ is the category of categories and $W$ consist
of all categorical equivalences. Then $W^{-1}C$ is the category
of categories and isomorphism classes of functors between them.
As functors can have non-trivial automorphisms, clearly
some important information (e.g. the automorphism groups
of the various functors) is lost when passing to $W^{-1}C$.
Another, deeper, example is the fact that Waldhausen $K$-theory is not an invariant
of homotopy categories, even when endowed with their natural triangulated
structures when they exist (see \cite{sch}). Other related troubles, as
the non functoriality of cones in derived categories, or the fact that
derived categories of sheaves in general do not form a stack (see e.g. the
introduction of \cite{hi-si}) also show
that the localization construction rapidly reaches its limits. \\

\begin{center} \textit{Model category} \end{center}

The major advance in abstract homotopy theory has been the notion
of model categories introduced by D. Quillen (who has been
followed by many different authors introducing various
modified versions). This notion had
an enormous impact and still today it seems difficult to
do serious homotopy theory without referring to model structures.
D. Quillen noticed that when a pair $(C,W)$ is furthermore structured and
admits a model structure then its homotopy theory on one hand becomes much more
easy to deal with and to describe, and on the other hand possesses several
additional important properties. First of all, the localized category
$W^{-1}C$, which in general is far from being easy to describe in
a useful way, possesses a very friendly presentation using
homotopy classes of maps between cofibrant and fibrant objects.
The second fundamental consequence of the existence of a model
structure is the existence of homotopy limits, homotopy colimits and of
mapping spaces, which are very important for many homotopical constructions.
Furthermore, homotopy limits, homotopy colimits and mapping spaces also
possess workable descriptions which allow to really deal with them.
These consequences of the existence of a model structure are definitely
non-trivial, and they show that the existence of a model structure
is never an easy nor a formal result.

Another very nice fact about model categories is that they tend to be
easily encountered and appear in many different contexts. This is probably due
to their good stability properties. For example, very often, presheaves (or even
sheaves) with values in a model category again form a model category. In the same
way, monoids objects, groups objects (or more general
abstract structures) in a model category again
form a model category. Starting from two fundamental examples of simplicial sets
and complexes over a ring, one constructs this way a lot of interesting
model categories.

However, model category theory also has its own limits. Indeed, all homotopical
invariants one is interested in a model category only depends on the
notion of equivalences, and not on the whole model structure. Model categories
are therefore too much structured objects, and tend to be a too much rigid notion
for certain purpose. For example, given two model categories $M$ and $N$, it
does not seem to exist a reasonable model category of functors
between $M$ and $N$. In other words, contrary to category theory, model category theory is not
an \textit{internal theory}, and this can cause troubles when one is interested
in model categories as objects in their own.
A good example of such a situation is when one considers a model category $M$,
a (fibrant) object $x \in M$, and the model category $M/x$ of objects over $x$.
Then, the simplicial monoid $aut(x)$ of self equivalences of $x$ is
expected to \textit{act} on the model category $M/x$, and this action does not seem
to be easy to describe inside the theory of model category itself.
This shows that
one is naturally lead to go beyond the theory of model category, though
every one seem to agree that model categories do include all
the examples one would like to study. From my point of view, the problem
encountered with model categories is not that they are too coarse
objects but on the contrary that they are over-structured objects. \\

\begin{center} \textit{D\'erivateurs} \end{center}

For a pair $(C,W)$ and a category $I$, one can construct
a new pair $(C^{I},W_{I})$, where $C^{I}$ is the category of functors
from $I$ to $C$ and $W_{I}$ the sub-category of natural transformations
in $C^{I}$ which levelwise belongs to $W$. Therefore, to an abstract homotopy theory
$(C,W)$ and a category $I$ one can construct a new abstract homotopy theory
of $I$-diagrams $(C^{I},W_{I})$.

The two categories
$(W^{-1}C)^{I}$ and $W_{I}^{-1}C^{I}$ are very different, and in general
the category $W^{-1}C$ alone does not determine the category
$W_{I}^{-1}C^{I}$. This is precisely one of the reason why
the localization construction $(C,W) \mapsto W^{-1}C$ is too crude, and why
the standard homotopical constructions (e.g. homotopy limits and colimits,
mappings spaces \dots) do depend on strictly more than  the
localized category $W^{-1}C$ alone.
The main idea of the theory of \textit{D\'erivateurs}, which seems
independently due to A. Heller and A. Grothendieck, is that from
an abstract homotopy theory $(C,W)$ one should not only consider
the localized category $W^{-1}C$, but
all the localized category of diagrams $W_{I}^{-1}C^{I}$ for various
index categories $I$. Precisely,  one should
consider the ($2$-)functor
$$\begin{array}{cccc}
\mathbb{D}_{(C,W)} : & \underline{Cat}^{op} & \longrightarrow & \underline{Cat} \\
 & I  & \mapsto & \mathbb{D}_{(C,W)}(I):=W_{I}^{-1}C^{I} \\
 & (f : I \rightarrow J) & \mapsto & f^{*},
\end{array}$$
from the (opposite) $2$-category of categories to itself.

According to A. Heller and A. Grothendieck the right context to do
homotopy theory is therefore the $2$-category of $2$-functors
$\underline{Cat}^{op} \longrightarrow \underline{Cat}$, called the
$2$-category of \textit{pr\'e-d\'erivateurs}, and denoted by
$\underline{PDer}$. \textit{D\'erivateurs} are then defined to be
pr\'e-d\'erivateurs satisfying certain additional property, as
for example the fact that the pull-backs functors $f^{*}$ possess
right and left adjoint (see \cite{he,gr2,ma1} for details). A fundamental
fact is that when a pair $(C,W)$ does admit a model structure
then the associated pr\'e-d\'erivateur $\mathbb{D}_{(C,W)}$ is a d\'erivateur
(this is essentially the existence of homotopy limits and colimits in
model categories).

The theory of d\'erivateurs has been quite successful
for many purposes. From a conceptual
point of view, the d\'erivateur associated with a model
category is a more intrinsic object than the model
category itself (e.g. it only depends on the notion
of equivalence and not on the whole model structure), which
furthermore contains a lot of the homotopy invariants
of the model category (as for example homotopy limits and colimits).
It has been used for example in order to state some kind of universal
properties that were lacking for derived categories and for the homotopy theory of spaces
(see \cite{ma1} for historical references).
The theory also has the advantage of solving many of the problems encountered
with the localization construction, as for example the non functoriality of
cones in derived categories (this was apparently one
of the motivations to introduce them).
Finally, on the contrary with model categories it has the advantage
of being an internal theory, in the sense that given two d\'erivateurs one can define
a reasonable (pr\'e-)d\'erivateur of morphisms between them.

Because of all of these nice properties the theory of d\'erivateurs
seems at first sight to be the \textit{right context for doing homotopy theory}. However,
there exist homotopical constructions that does not factor through
the theory of d\'erivateurs, and the fundamental reason is that d\'erivateurs
form a $2$-category, which for many purposes is a too coarse structure
missing some important higher homotopical information.
One can make for instance the same kind of remark as for model categories.
Given a model category $M$ and $x$ a (fibrant) object in $M$, the simplicial
monoid $aut(x)$ is expected to act on the d\'erivateurs associated
to $M/x$. However, as this d\'erivateur lives in a $2$-category this action
automatically factors through an action of the $1$-truncation $\tau_{1}aut(x)$, and therefore
one sees that the d\'erivateurs associated to $M/x$ does not see
the higher homotopical information encoded in the action of the whole
space $aut(x)$. This fact implies for example that the
$K$-theory functor (in the sense Waldhausen)
can not be reasonably defined on the level of pr\'e-d\'erivateurs (see Prop. \ref{p2}
and Cor. \ref{c3}).
As explained by the result Thm. \ref{t6}, the theory of pr\'e-d\'erivateurs
is only a approximation \textit{up to $2$-homotopies} of a more
complex object encoding higher homotopical data. \\

\begin{center} \textit{Simplicial localization} \end{center}

In \cite{dk1} B. Dwyer and D. Kan introduced a refined version of the
localization construction $(C,W) \mapsto W^{-1}C$, which associates
to a pair $(C,W)$ an $S$-category (i.e. a category enriched over the
category of simplicial sets) $L(C,W)$ whose category of connected
component $Ho(L(C,W))$ is naturally isomorphic to $W^{-1}C$.
As the localized category $W^{-1}C$ satisfies a universal property in the category
of categories, the $S$-category $L(C,W)$ satisfies some universal
property, \textit{up to equivalence}, in the category of $S$-categories.
One of the main result proved by B. Dwyer and D. Kan states that
when the pair $(C,W)$ has a model structure then
the simplicial localization $L(C,W)$ can be described using the
mapping spaces defined in terms of fibrant-cofibrant resolutions (see \cite{dk2}).
This last result is an extension to higher homotopies of the
well known description of the category $W^{-1}C$ in terms
of homotopy classes of maps between fibrant and cofibrant objects.

The very nice property of the $S$-category $L(C,W)$ is that it
seems to contain all of the interesting homotopical information
encoded in the pair $(C,W)$. For example, when $(C,W)$ is endowed
with a model structure, the mapping spaces as well as the homotopy limits and
colimits can all be reconstructed from $L(C,W)$. Also, as shown
by Thm. \ref{t1}, the $K$-theory functor does factor through the theory of
$S$-categories, which shows that $L(C,W)$ contains definitely more
information than the localized category $W^{-1}C$ (even when endowed
with its additional triangulated structure, when it exists).
Actually, the $S$-category $L(C,W)$ almost reconstructs, in a sense I will not
precise here, the model category $(C,W)$.

However, once again the theory of $S$-categories is not
well behaved with respect to categories of functors
and suffers the same troubles than the theory of model categories. Indeed,
it does not seem so easy to define for two $S$-categories
a reasonable $S$-category of morphisms between them. Of course
there exists a natural $S$-category of morphisms between two $S$-categories
but one can notice very easily that it does not have the right homotopy type
(for example it is not invariant under equivalences of $S$-categories).
Actually, there exist well known conceptual reasons why the category of $S$-categories
can not be enriched over itself in a homotopical meaningful manner (see \cite{hi-si}, the remark after Problem 7.2).
The situation is therefore very similar to the case of model categories,
as $S$-categories seem to model all what we want but does not
provide and internal theory. In any case, one advantage of the notion $S$-categories
compare to model categories is that it is more intrinsic, and the $S$-category
associated to a model category only depends on its notion of equivalence
and not on the whole model structure. The theory of $S$-categories is
also better suited than the theory of d\'erivateur as
there does not seem to exists homotopical invariants of model categories that
can not be reconstructed from its simplicial localization. \\

\begin{center} \textit{Segal categories} \end{center}

Segal categories are weak form of $S$-categories, in which
composition is only defined \textit{up to a coherent system of
equivalences} (see \cite{hi-si,pe} for details). An $S$-category is in an
obvious way a Segal category, and any Segal category is equivalent
to an $S$-category. More generally, the homotopy theory of
$S$-categories and of Segal categories are equivalent (see
\cite[\S 2]{si1}). However, the main advantage of Segal categories is that
they do form an internal theory. Given two Segal categories $A$
and $B$, there exists a Segal category of morphisms
$\mathbb{R}\underline{Hom}(A,B)$ satisfying all the required properties
(as the usual adjunction rule, invariance by equivalences \dots). When
applied to two $S$-categories $A$ and $B$, the Segal category
$\mathbb{R}\underline{Hom}(A,B)$ should be interpreted as
the Segal category of \textit{lax simplicial functors from $A$ to $B$}.
In conclusion, the theory of Segal categories is equivalent to the
theory of $S$-categories, but does behave well with respect to
categories of functors.

For any abstract homotopy theory $(C,W)$ one has
the $S$-category $L(C,W)$, which can be considered
as a Segal category.
In particular, for two model categories $M$ and $N$, one can consider
their simplicial localizations $LM$ and $LN$ (along the
sub-categories of equivalences), and then consider
the Segal category $\mathbb{R}\underline{Hom}(LM,LN)$.
The Segal category $\mathbb{R}\underline{Hom}(LM,LN)$
precisely plays the role of the non-existing model
categories of morphisms from $M$ to $N$. For a (fibrant) object
$x$ in $M$ one can also make sense of the action of
the simplicial monoid $aut(x)$ on the Segal category $L(M/x)$
(the model category of Segal category is a simplicial model category, and
therefore an action of a simplicial monoid on a Segal category makes
perfect sense).

The conclusion is that the theory of Segal categories can be used
in order to make constructions with model categories that can not
be done inside the theory of model categories itself. Several
examples of such constructions are given in this m\'emoire. It
worth also mentioning that model categories and Segal categories
essentially model the same objects. Of course, it is not true that
any Segal category can be written as some $LM$ for a model
category $M$. However, any Segal category can be written as
$L(C,W)$, where $C$ is a full sub-category of a model category
which is closed by equivalences, and $W$ is the restriction of
equivalences for the ambient model structure. This remarks shows
that, at least on the level of objects, Segal categories and model
categories are essentially the same thing. What Segal categories
really bring are new, much more flexible and powerful,
well behaved functoriality properties. \\

\begin{center} \textit{Other approaches} \end{center}

There exist other theories that I will not consider in this text, but which
definitely behaves in very similar manners than Segal categories, and for which
one could make the same analysis. I am thinking in particular to
the theory of Quasi-categories of A. Joyal (see \cite{jo}), and the
theory of complete Segal spaces of C. Rezk (see \cite{re}). Some comparison functors
between these theories and the theories of $S$-categories and Segal categories
are given in appendix B. \\

\bigskip

What shows the above overview of the various ways of doing
homotopy theory is that the combination of model category theory
and Segal category theory seems the most suited for many purposes.
Of course, one could also try to use only Segal category theory,
but it is sometimes useful to use model categories which in
practice provide a much more friendly setting than Segal
categories. Roughly speaking, model categories are used in order
to do explicit computations and Segal categories are rather used
in order to produce abstract constructions
(though this way of thinking is rather artificial). A key result which
allows to really do this in practice is the so-called \textit{strictification
theorem} of \cite[Thm. 4.2.1]{to-ve1}.
By analogy, one could compare the relations between model categories
and Segal categories with the existing relations between Grothendieck sites
and topoi (passing from a site to its associated topos being the
analog of the simplicial localization construction).
I personally like to think that writing a Segal category $A$ as some simplicial localization
$LM$ for a model category $M$ is very much like \textit{choosing coordinates
on $A$}, in the same way that one choose coordinates on an abstract
manifold when one does local computations. The relationship between
Segal categories and model categories will be described with more
details in the next section.

\section{Segal categories}

In this first section I will present some results and constructions about Segal categories, or
that use Segal categories in an essential way. They are
only few examples of the kind of results one can obtain using
Segal categories, and I am convinced that more interesting results
could also be proved in the future.

\subsection{Segal categories and model categories}

An $S$-category is by definition a category enriched over the category
of simplicial sets. Segal categories are weak form of $S$-categories
were the composition is only defined \textit{up to a coherent system
of equivalences}. For the details on the theory of $S$-categories and Segal categories
we refer to \cite{dk1,dk2,hi-si,pe}, and to \cite{to-ve1} for an short overview and notations.
From a purely esthetic point of view, it is also useful to think of Segal categories as $\infty$-categories
in which $i$-morphisms are invertible (up to $i+1$-isomorphisms) for any $i>1$.

The reader should keep in mind that Segal category theory works in a very similar
manner than usual category theory and that most (if not all) of
the standard categorical notions can be reasonably defined in the Segal setting.
Here follows a sample of examples of such (once again we refer to
the overview \cite{to-ve1} for more details).

\begin{itemize}

\item As any $S$-category, any Segal category $A$ possesses a homotopy
category $Ho(A)$ (which is a category in the usual sense),
having the same objects as $A$ and homotopy classes
of morphisms of $A$ as morphisms between them. We also recall that
given two objects $a$ and $b$ in a Segal category, morphisms
between $a$ and $b$ in $A$ form a simplicial set denoted by
$A_{(a,b)}$.

For a morphism of Segal categories $f : A \longrightarrow B$, one says that
$f$ is essentially surjective (resp. fully faithful) if
the induced functor $Ho(f) : Ho(A) \longrightarrow Ho(B)$ is essentially surjective (resp.
if for any two objects $a$ and $b$ in $A$ the induced morphism
$f_{(a,b)} : A_{(a,b)} \longrightarrow B_{(f(a),f(b))}$ is an
equivalence of simplicial sets). Of course one says that $f$ is an equivalence if it is
both fully faithful and essentially surjective.

\item The foundational result about Segal categories is the existence
of a model structure (see \cite{hi-si,pe}). Segal categories are particular cases of
Segal pre-categories, and the category of Segal pre-categories
is endowed with a cofibrantely generated model structure. Every object
is cofibrant, and the fibrant objects
for this model structure are Segal categories, but not all Segal category is
a fibrant object, and in general fibrant objects are quite difficult to
describe. The model category is furthermore enriched over itself (i.e.
is an internal model category in the sense of
\cite[\S 11]{hi-si}). This implies that
given two Segal categories $A$ and $B$ one can associate a Segal category
of morphisms
$$\mathbb{R}\underline{Hom}(A,B):=\underline{Hom}(A,RB),$$
where $RB$ is a fibrant model for $B$ and $\underline{Hom}$ denote the
internal $Hom$'s in the category of Segal pre-categories.
From the point of view of $\infty$-categories, $\mathbb{R}\underline{Hom}(A,B)$
is a model for the $\infty$-category of (lax) $\infty$-functors
from $A$ to $B$. In general, the expression \textit{$f : A \longrightarrow B$ is
a morphism of Segal categories} will mean that $f$ is an object
in $\mathbb{R}\underline{Hom}(A,B)$. In other words we implicitly allows
ourselves to first take a fibrant replacement of $B$ before considering
morphisms into $B$.

\item There is a notion of \emph{Segal groupoid}, which is to Segal category
theory what groupoids are for category theory. By definition,
a Segal category $A$ is a Segal groupoid if its homotopy category
$Ho(A)$ is a groupoid in the usual sense (or equivalently, if any morphism
in $A$ is a homotopy equivalence).

Furthermore for any Segal category $A$, one can define
its geometric realization $|A|$, which is the diagonal simplicial
set of the underlying bi-simplicial of $A$ (see \cite[\S 2]{hi-si}, where
$|A|$ is denoted by $\mathcal{R}_{geq 0}(A)$). The construction
$A \mapsto |A|$ has a right adjoint, sending a simplicial
set $X$ to its fundamental Segal groupoid $\Pi_{\infty}(X)$
(this one is denoted by $\Pi_{1,se}(X)$ in
\cite[\S 2]{hi-si}).
By definition, the set of object of $\Pi_{\infty}(X)$
is the set of $0$-simplicies in $X$, and
for $(x_{0},\dots,x_{n}) \in X_{0}^{n+1}$ the simplicial set
$\Pi_{\infty}(X)_{n}$ is the sub-simplicial set of
$X^{\Delta^{n}}$ sending the $i$-th vertex of $\Delta^{n}$ to
$x_{i}\in X_{0}$. A fundamental theorem states that
the constructions $A \mapsto |A|$ and
$X \mapsto \Pi_{\infty}(X)$ provide an equivalence
between the homotopy theories of
Segal groupoids and of simplicial sets (see \cite[\S 6.3]{pe}). This last
equivalence is a Segal version of
the well known equivalence between the homotopy theories of
$1$-truncated homotopy types and of groupoids.

\item Given a Segal category $A$ and a set of morphisms $S$ in $Ho(A)$ (we recall
that $Ho(A)$ is the homotopy category of $A$), one can construct a
Segal category $L(A,S)$ by formally inverting the arrows of $S$. This construction
is the Segal analog of the Gabriel-Zisman localization for categories. By definition,
the Segal category $L(A,S)$ comes with a localization morphism $l : A \longrightarrow L(A,S)$
satisfying the following universal property: for any Segal category $B$, the
induced morphism
$$l^* : \mathbb{R}\underline{Hom}(L(A,S),B) \longrightarrow \mathbb{R}\underline{Hom}(A,B)$$
is fully faithful, and its essential image consists of
morphisms $A \longrightarrow B$ sending morphisms of $S$ into equivalences in $B$ (i.e.
isomorphisms in $Ho(B)$).

When applied to the case where $C$ is a category considered as a Segal category, the construction
$L(A,S)$ described above coincides, up to an equivalence, with the simplicial localization
construction of \cite{dk1}. From the $\infty$-category point of view this means that
$L(C,S)$ is the \textit{$\infty$-category obtained from $C$ by formally inverting
the arrows in $S$}. The fact that this localization procedure produces $\infty$-categories
instead of categories is, from my point of view,
the deep heart of the relations between homotopy theory and higher category theory.

\item Given a model category $M$, one can construct a Segal category $LM:=L(M,W)$ by localizing
$M$, in the Segal category sense, along its sub-category of equivalences $W$.
This gives a lot of examples of Segal categories. Using the main result of
\cite{dk3} the Segal categories $LM$ can be explicitly describe in terms
of mapping spaces in $M$. In particular, when $M$ is a simplicial model category
$LM$ is equivalent to the simplicial category of fibrant-cofibrant objects in $M$.
For the
model category of simplicial sets we will use the notation $Top:=LSSet$. The Segal
category $Top$ is as fundamental as the category of sets in category theory.

Given a Segal category $A$ one can construct a Yoneda embedding
morphism
$$h : A \longrightarrow \mathbb{R}\underline{Hom}(A^{op},Top),$$
which is know to be fully faithful (this is the Segal version of
the Yoneda lemma). Any morphism $A^{op} \longrightarrow Top$ in the essential image of this
morphism is called \textit{representable}. Dually, one has a notion of
\textit{co-representable morphism}.

\item Given a morphism of Segal categories $f : A \longrightarrow B$, one says that
$f$ has a right adjoint if there exists a morphism $g : B \longrightarrow A$
and a natural transformation $h \in \mathbb{R}\underline{Hom}(A,A)_{(Id,gf)}$,
such that for any two objects $a \in A$ and $b\in B$ the natural
morphism induced by $h$
$$\xymatrix{A_{(f(a),b)} \ar[r]^-{g_{*}} & A_{(gf(a),g(b))} \ar[r]^-{h^{*}} &  A_{(a,g(b))}}$$
is an equivalence of simplicial sets.
This definition allows one to talk about adjunction between Segal categories.
An important fact is that a Quillen adjunction
between model categories
$$f : M \longrightarrow N \qquad M \longleftarrow N : g$$
gives rise to a natural adjunction of Segal categories
$$Lf : LM \longrightarrow LN \qquad LM \longleftarrow LN : Lg.$$

\item Given two Segal categories $A$ and $I$, one says that $A$ has limits (resp. colimits) along $I$ if the
constant diagram morphism $A \longrightarrow \mathbb{R}\underline{Hom}(I,A)$
has a right adjoint (resp. left adjoint). This allows one to talk about
Segal categories having (small) limits (resp. colimits), or finite limits (resp. colimits).
In particular one can talk about fibered and cofibered square, final and initial
objects, left and right exactness \dots.

\item Existence of Segal categories of morphisms also permits to define
notions of algebraic structures in a Segal category. For example, if $A$ is
a Segal category with finite limits, the Segal category of monoids in $A$ is
the full sub-Segal category of $\mathbb{R}\underline{Hom}(\Delta^{op},A)$
consisting of morphisms $F : \Delta^{op} \longrightarrow A$ such that
$F([0])=*$ and such that the Segal morphisms $F([n]) \rightarrow F([1])^{n}$
are equivalences in $A$. One can also defines this way groups, groupoids, rings, categories
\dots in $A$.

\item More advanced notions for Segal categories, as topologies, stacks and topos
theory, or monoidal structures will be given in \S 3 and \S 5.

\end{itemize}

The previous list (highly non-exhaustive) of examples of standard constructions
one can do with Segal categories is very much useful in practice as it
allows to use Segal category theory as the category theory we have learned at school.
However, as fibrant objects and fibrant resolutions are extremely difficult to
describe these constructions turn out to be quite hard to manipulate
in concrete terms. This difficult is solved by the so-called \textit{strictification theorems},
which stipulate that when the Segal categories involved are of the form
$LM$ for $M$ a model category, all of these categorical constructions
can be expressed in terms of standard categorical construction inside
the well known world of model categories. The most important strictification theorem,
concerning categories of diagrams, is the following.

Let $M$ be a simplicial model category (in the sense of \cite{ho}), and $T$ be an $S$-category.
Let $M^{T}$ be the category of simplicial functors from $T$ to $M$. A morphism
$f : F \rightarrow G$ in $M^{T}$ will be
called an equivalence if for any object $t \in T$ the induced morphism
$f_t : F(t) \rightarrow G(t)$ is an equivalence in $M$. By the universal property
of the localization construction, one defines a natural morphism of Segal categories
$$L(M^{T}) \longrightarrow \mathbb{R}\underline{Hom}(T,LM),$$
where the localization on the left is perform with respect to
the above notion of equivalences in $M^{T}$.

\begin{thm}{(Hirschowitz-Simpson, \cite[Thm. 18.6]{to-ve1})}\label{tstrict}
Under the previous assumption, the natural morphism
$$L(M^{T}) \longrightarrow \mathbb{R}\underline{Hom}(T,LM)$$
is an equivalence of Segal categories.
\end{thm}

Theorem \ref{tstrict} has many important consequences. First of all, it is
the key argument in the proof of the Yoneda lemma for Segal categories. Furthermore, the
Yoneda lemma implies that any Segal category $A$ is equivalent to a full sub-Segal category
of $L(SSet^{T})$, where $T$ is an $S$-category equivalent to $A^{op}$. This implies that
any Segal category can be represented up to equivalences by a full sub-$S$-category
of fibrant-cofibrant simplicial presheaves on a category. In particular, Segal categories
and model categories are very close, and from my point of view
are essentially the same kind of objects (i.e. $\infty$-categories
where $i$-arrows are invertible for $i>1$).
Another, very important,
consequence of theorem \ref{tstrict} is that for a model category $M$ the Segal
category $LM$ has all limits and colimits, and these can be concretely computed
in terms of homotopy limits and homotopy colimits in $M$. \\

In conclusion, the localization functor $L$ permits to pass from
model categories to Segal categories for which many interesting
and abstract categorical constructions are available (e.g.
Segal categories of functors). The strictification theorem
then stipulates that these constructions on the level of Segal categories
do have model category interpretations which in practice allows
to reduce problems and computations to model category theory
(and therefore to standard category theory). This principle will
be highly used all along this work but in an rather implicit manner, and corresponding constructions
in model category and Segal categories will always be identified in some sense:
for a given
situation I will simply used the most appropriate theory.
I personally like to
deal as much as possible with Segal categories for general constructions, but
use model
categories to provide explicit descriptions. This does not
quite follow the general philosophy, which I agree with,
that Segal category theory (or any equivalent theory) should at some point completely replace
model category theory and provide a much more powerful and friendly setting.
My personal
feeling is that the theory of Segal category is very
much at its starting point and still we are
not totally confident with the kind of arguments and constructions that one
is allowed to use. Model categories on the other hand as been highly studied since
many years, and I have the feeling that keeping an eye on model category theory while
dealing with Segal categories will help us to learn how the theory really works.
In particular
this way of proceeding should
provide a whole list of arguments, manipulations and constructions one can safely used,
which I hope will be part of the standard mathematical knowledge in the future much in the
same way as category theory is today. \\

To finish this paragraph on Segal category theory let me mention
the existence of higher Segal categories. I will not use this higher
notion very often, but it will happen that the notion of $2$-Segal categories
is needed. Also, I will implicitly use the \emph{change of $n$} constructions
given in \cite[\S 2]{hi-si}, and always consider that
a Segal category is also in a natural way a $2$-Segal category. For all
details I refer to \cite{hi-si} and \cite{pe}. \\

\subsection{$K$-Theory}

In the sixties A. Grothendieck asked whether or not $K$-theory
is an invariant of triangulated categories. This question has
been studied by several authors (R. Thomason, A. Neeman \dots) and we now know that
the answer is negative: there is no reasonable $K$-theory functor
defined on the level of triangulated categories (see \cite{sch}).
In other words, given a Waldhausen category $C$, the homotopy category
$Ho(C)$, even endowed with its triangulated structure when it exists,
is not sufficient to recover the $K$-theory spectra $K(C)$ (or even
the $K$-theory groups $K_{n}(C)$).

However, for a Waldhausen category $C$,
one can consider
its simplicial localization $LC$ defined by Dwyer and Kan and which
is a refinement of the homotopy category $Ho(C)$ (see \cite{dk1}). It is an $S$-category
that can be considered as a Segal category, and under some conditions
on $C$ (e.g. when it is \textit{good} in the sense
of \cite{to-ve2}, though it seems the result stays correct under much weaker
assumptions) $LC$ is enough to recover the $K$-theory spectrum $K(C)$.

\begin{thm}{(To\"en-Vezzosi, \cite{to-ve2})}\label{t1}
The $K$-theory spectrum $K(C)$ of a \emph{good} Waldhausen
category $C$ can be recovered (up to equivalence) functorially from
the Segal category $LC$.
\end{thm}

An important non-trivial corollary of the previous theorem is the
following. For a Segal category $A$ we let $aut(A)$ its
simplicial monoid of self-equivalences (as defined for example
as in \cite{dk1}). Concretely, $aut(A)$ is the classifying space
of the maximal sub-Segal groupoid $\mathbb{R}\underline{End}(A)^{int}$
of the Segal category of endomorphisms $\mathbb{R}\underline{End}(A)$.

\begin{cor}\label{c1}
For a \emph{good} Waldhausen category $C$ the simplicial monoid
$aut(LC)$ acts naturally on the
Waldhausen $K$-theory spectrum $K(C)$.
\end{cor}

The proof of the above results given in \cite{to-ve2} is direct and does not involve
Segal category techniques. However, as mentioned at the end of \cite{to-ve2}, one
could also prove theorem \ref{t1} by first defining a $K$-theory functor
on the level of Segal categories having finite limits,
and then proving
that when applied to $LC$ for a (good) Waldhausen category $C$ the two constructions
coincide. Without going into too technical details let us just mention
the following unpublished result which follows from the results of \cite{to-ve2}, the construction
sketched at the end of \cite{to-ve2} and the strictification theorem of \ref{tstrict}. For this,
we let $GWCat$ be the category of good Waldhausen categories and exact
functors. In $GWCat$, an exact functor $F : C \longrightarrow D$
is called an $L$-equivalence if the induced morphism $LF : LC \longrightarrow LD$
is an equivalence of Segal categories. We let
$LGWCat$ be the Segal category obtained from $GWCat$ by
applying the simplicial localization functor with respect to the
$L$-equivalences. In the same way, let $LSeCat_{*}$ be the Segal category
obtained from the category of pointed Segal categories
\footnote{A Segal category is pointed if it has an initial
object which is also a final object.} by applying the simplicial
localization functor with respect to equivalences. We denote by
$LSeCat_{*}^{fl}$ the sub-Segal category of $LSeCat_{*}$ consisting
of Segal categories having finite limits and left exact functors.
Finally, let $LSp$ the
simplicial localization of the model category of spectra. The Waldhausen
$K$-theory functor induces a well defined morphism of Segal categories
$$K^{Wal} : LGWCat \longrightarrow LSp.$$
In the same way, the simplicial localization functor induces a morphism
$$L : LGWCat \longrightarrow LSeCat_{*}^{fl}.$$

\begin{thm}\label{t2}
There exists a commutative diagram in the homotopy category of Segal categories
$$\xymatrix{
 LGWCat \ar[dd]_-{L}\ar[rr]^-{K^{Wal}} & & LSp \\
& \\
LSeCat_{*}^{fl}. \ar[rruu]_-{K^{Segal}} & &  }$$
\end{thm}

The conclusion of theorem \ref{t2} really is: \\

\textit{$K$-theory is an invariant of Segal categories.} \\

\subsection{Stable Segal categories}

In the last part, we have seen that the $K$-theory spectrum of a
Waldhausen category can be recovered from
its simplicial localization $LC$ without any additional structures. At first
sight, this might look surprising as several works around this
type of questions involves triangulated structures (see \cite{du-sh,ne}).
In fact, the $S$-category $LC$ completely determines the
notion of fiber and cofiber sequences in the homotopy category $Ho(C)$.
In particular, the category
$Ho(C)$ together with its triangulated structure (when it exists)
is completely determined by the Segal category $LC$.

This observation has led A. Hirschowitz, C. Simpson and myself to
introduce a notion of \textit{stable Segal categories}. This notion
clearly is very close to the notions of \textit{enhanced
triangulated categories} (see \cite{bo-ka}), of
\textit{triangulated $A_{\infty}$-categories} of M. Kontsevich, and of
\textit{stable model categories} (see \cite[\S 7]{ho}).

\begin{df}{(Hirschowitz-Simpson-To\"en)}\label{d1}
A Segal category $A$ is \emph{stable} if it satisfies the following three conditions.
\begin{enumerate}
\item The Segal category $A$ possesses finite limits and colimits (in particular
its a final object and an initial object).

\item The final and initial object in $A$ are equivalent.

\item The suspension functor
$$\begin{array}{cccc}
S : & Ho(A) & \longrightarrow & Ho(A) \\
 & x & \mapsto & *\coprod_{x}*
\end{array}$$
is an equivalence of categories.
\end{enumerate}
\end{df}

\begin{rmk}\label{r2}
\emph{Though the notion of stable Segal category is not
strictly speaking a generalization of the notion of abelian categories
(a category which is stable in the Segal sense is trivial),
stable Segal categories really play the role of abelian categories in the Segal setting.}
\end{rmk}

The main properties of stable Segal categories are gathered in the following unpublished theorem.

\begin{thm}\label{t3}
\begin{enumerate}
\item The homotopy category $Ho(A)$ of a stable Segal category $A$
has a natural triangulated structure, for which triangles are induced
by the images of fiber sequences in $A$. Any exact morphism $f : A \longrightarrow B$
between stable Segal categories induces a functorial triangulated functor
$Ho(f) : Ho(A) \longrightarrow Ho(B)$.

\item A morphism between stable Segal category is left exact if and only
if it is right exact.

\item For any stable Segal category $A$ and any Segal category $B$
the Segal category $\mathbb{R}\underline{Hom}(B,A)$
is stable.

\item If $C$ is a full sub-category of a stable model category $M$ (in the sense
of \cite[\S 7]{ho}) which is closed by equivalences $M$ and
homotopy fibers and contains the initial-final object,
then the Segal category $LC$ is stable.

\item For any stable Segal category $A$, there exists a full sub-category $C$
of a stable model category and which is closed by equivalences in $M$ and
homotopy fibers
such that $LC$ is equivalent to $A$. Furthermore, one can chose
$M$ to be a model category of presheaves of spectra over some category.

\end{enumerate}
\end{thm}

The previous theorem clearly shows that stable Segal categories are quite close
to triangulated categories. However, the additional structure encoded
in Segal categories allows one to have the fundamental property $(3)$, which
is violated for triangulated categories. Furthermore,
the theorem \ref{t1} together with the counter-example given in \cite{sch}, show that
there exists two non-equivalent stable Segal categories whose homotopy categories
are equivalent as triangulated categories. This of course implies that
a stable Segal category $A$ contains strictly more information than
its triangulated homotopy category $Ho(A)$. Because of all these
reasons, we propose the notion of stable Segal categories
as an alternative to the notion of triangulated categories, and we think that
several troubles classically encountered with triangulated categories can be solved
this way. This is a reasonable thing to do as
all of the interesting triangulated categories we are award of
(in particular all of triangulated categories of geometric
origin, as for example derived categories of sheaves)
are of the form $Ho(A)$ for some stable Segal category $A$.
I personably tend to think that triangulated
categories which are not of the form $Ho(A)$ for some stable
Segal categories $A$ are unreasonable object and should not be
considered at all.

To finish, let me mention that
the notion of stable Segal categories has been already used
in several contexts, as for examples the Tannakian formalism for
Segal categories (see \S 5) and the stacks point of view
on Grothendieck's duality theory (Hirschowitz-Simpson-To\"en, unpublished). \\

\subsection{Hochschild cohomology of Segal categories}

In this part I present some thoughts about the notion of
Hochschild cohomology of Segal categories. My main objective was
to understand in which sense \textit{Hochschild cohomology is the
space of endomorphisms of an identity functor}, as stated by several
authors (see e.g. \cite{sei}). The point of view taken in this paragraph
is non-linear, and is concerned with the discrete version of
the Hochschild cohomology. Taking linear structures into account is
probably more complicated and actually a bit tricky, but I guess it should
be doable.

All the material of this section is not written up, except the last result
proved in collaboration with J. Kock (see \cite{ko-to}).\\

For a Segal category $A$ one has its Segal category
of endomorphisms $\mathbb{R}\underline{End}(A)$. One consider the
object $Id \in \mathbb{R}\underline{End}(A)$ and its
simplicial set of endomorphisms
$\mathbb{R}\underline{End}(Id_{A}):=\mathbb{R}\underline{End}(A)_{(Id_{A},Id_{A})}$.

\begin{df}\label{d2}
The \emph{Hochschild cohomology} of a Segal category $A$ (also called
the \emph{center} of $A$) is the
simplicial set
$$\mathbb{HH(A)}:=\mathbb{R}\underline{End}(Id_{A}).$$
\end{df}

\begin{rmk}\label{r1}
\emph{When applied to the Segal category $LM$ for a model category
$M$, definition \ref{d2} gives a notion of the} Hochschild cohomology of
the model category $M$.
\end{rmk}

This definition is a direct generalization of the center of a category
$C$, defined as the monoid of endomorphisms of the identity functor
of $C$. For a category
$C$, its center is always a commutative monoid. This follows from a very standard
argument as the center is always endowed with two compatible unital and associative
composition laws. For a Segal category $A$, we will see that $\mathbb{HH}(A)$ is not
quite commutative, but is a \textit{$2$-Segal monoid}, or in other words is endowed with
two compatible unital and associative \textit{weak} composition laws. Therefore,
$\mathbb{HH}(A)$ looks very much like a $2$-fold loop space but for which
the composition laws are not necessarily invertible. In order to state
precise results let me start by some definitions (the reader could
consult \cite{ko-to} for details). \\

We let $C$ be a category with a notion of equivalences
(e.g. a model category or a Waldhausen category) and finite products, and such that
finite products preserve equivalences. A \emph{Segal monoid} in $C$ is a functor
$$\begin{array}{cccc}
H : & \Delta^{op}  & \longrightarrow & C \\
& [n] & \mapsto & H_{n},
\end{array}$$
such that
\begin{enumerate}
\item $H_{0}=*$

\item For any $n\geq 1$ the Segal morphism (see \cite{hi-si} or \cite[\S 2]{to-ve1})
$$H_{n} \longrightarrow H_{1}^{n}$$
is an equivalence.

\end{enumerate}

A morphism between Segal monoids is simply a natural transformation, and
a morphism is an equivalence if it is an equivalence levelwise (or equivalently
on the image of $[1]$). Clearly, Segal monoids in $C$ do form a category
$SeMon(C)$, again with a notion of equivalences. Furthermore, $SeMon(C)$
has finite products which again preserve equivalences. The construction can
therefore by iterated.

\begin{df}\label{d4}
Let $C$ be a category with a notion of equivalences and finite products, and such that
finite products preserve equivalences.
\begin{enumerate}
\item The category of \emph{$d$-Segal monoids} in $C$ is defined inductively by
$$0-SeMon(C)=C \qquad d-SeMon(C):=SeMon((d-1)-SeMon(C)).$$
$d$-Segal monoids in simplicial sets are simply called $d$-Segal monoids.

\item The underlying object of a Segal monoid $A$ in $C$ is
$A_{1} \in C$. Inductively, the underlying object of a $d$-Segal monoid
$A$ in $C$ is the underlying object of the $(d-1)$-Segal monoid
$A_{1}\in (d-1)-SeMon(C)$. The underlying object
of a $d$-Segal monoid $A$ is again denoted by $A$.

\item The category of \emph{$d$-fold monoidal Segal categories}
is the category of $d$-Segal monoids in the category of Segal categories.

\end{enumerate}
\end{df}

Almost by definition, one has the following elementary proposition.

\begin{prop}\label{p1}
Let $A$ be a $d$-fold monoidal Segal category, and still denote by $A$ the
underlying Segal category. The Hochschild cohomology of $A$, $\mathbb{HH}(A)$,
has a natural structure of a $(d+2)$-Segal monoid.
\end{prop}

\begin{rmk}\label{rp1}
\begin{enumerate}
\item
\emph{The expression} has a natural structure of a $(d+2)$-Segal monoid
\emph{is a bit ambiguous. It as of clearly to be understood
up to equivalence. More precisely, this means that
there exists a natural $(d+2)$-Segal monoid $M$ whose underlying object
is naturally equivalent to $\mathbb{HH}(A)$.}
\item
\emph{Of course, the $(d+2)$-Segal monoid structure
on $\mathbb{HH}(A)$ depends on
the $d$-fold monoidal structure given on $A$.}
\end{enumerate}
\end{rmk}

The above proposition becomes really interesting with the following theorem which
relates the definition of Hochschild cohomology of a Segal category to
a more usual one. Let us recall first that for a simplicial monoid
$H$ one has a notion of $H$-modules, or equivalently of simplicial sets
with an action of $H$. The category of $H$-module is known to be a model
category for which equivalences are defined on the underlying simplicial set.
The main theorem, which is
a direct consequence of the strictification theorem \ref{tstrict}, is the
following.

\begin{thm}\label{t4}
Let $H$ be a simplicial monoid and $A=BH$ be the Segal category with one
object and $H$ as endomorphisms of this object. Then, one has a natural equivalence
of simplicial sets
$$\mathbb{HH}(BH)\sim Map_{H\times H^{op}-Mod}(H,H),$$
where the right hand side is the mapping space computed in the model
category of $H\times H^{op}$-modules.
\end{thm}

An important corollary is a non-linear analog of Deligne's conjecture
on the complex of Hochschild cohomology of an associative differential graded algebra
(see \cite{ko-to} for detailed references).
It follows from proposition \ref{p1} and theorem \ref{t4} and by the observation that
if $H$ is a $d$-Segal monoid then $BH$ is a naturally a $(d-1)$-fold monoidal Segal category
(the notions of modules over simplicial monoids extend naturally
to Segal monoids, using for example that any Segal monoid is
naturally equivalent to a simplicial monoid).

\begin{cor}\label{c2}
Let $H$ be a $d$-Segal monoid (in simplicial sets). Then the
simplicial set $Map_{H\times H^{op}-Mod}(H,H)$ has
a natural structure of a $(d+1)$-Segal monoid.
\end{cor}

By different techniques, which use in an essential way
the simplicial localization functor of Dwyer and Kan one can also
prove the following related result. It does imply corollary
\ref{c2} for $d=1$, but also has its own interest. It is a model category analog
of the fact that the endomorphism of the unit in a monoidal category is
a commutative monoid.

\begin{thm}{(Kock-To\"en, \cite{ko-to})}\label{t5}
Let $M$ be a monoidal model category in the sense of \cite[\S 4.3]{ho} and
$1_{M}$ be its unit. Then the simplicial set
$Map_{M}(1_{M},1_{M})$ has a natural structure of a
$2$-monoid.
\end{thm}

\subsection{Segal categories and d\'erivateurs}

In this last paragraph I will compare the theory of \textit{d\'erivateurs}
of A. Heller and A. Grothendieck with the theory of Segal categories. The main theorem
states that the theory of d\'erivateurs is essentially a $2$-truncation of the
theory of Segal categories, and so the two theories are more or less equivalent
\textit{up to $2$-homotopies}. The results of this paragraph have not been written
up. \\

We denote by $\underline{PDer}$ the $2$-category of pr\'e-d\'erivateurs in the sense
of \cite{ma1}. Recall that $\underline{PDer}$ is the $2$-category of
$2$-functors $\underline{Cat}^{op} \longrightarrow \underline{Cat}$.
For a Segal category $A$ we define an object $\mathbb{D}_{A} \in \underline{PDer}$
in the following way
$$\begin{array}{cccc}
\mathbb{D}_{A} : & \underline{Cat}^{op} & \longrightarrow & \underline{Cat} \\
 & I & \mapsto & Ho(\mathbb{R}\underline{Hom}(I,A)).
\end{array}$$
The construction $A \mapsto \mathbb{D}_{A}$ clearly defines a
morphism of $2$-Segal categories
$$\mathbb{D} : \underline{SeCat} \longrightarrow \underline{PDer},$$
where the left hand side is the $2$-Segal category of Segal categories
as defined in \cite[\ S]{hi-si}.
As the $2$-Segal category $\underline{PDer}$ is a $2$-category, this morphism factors as a
$2$-functor between $2$-categories
$$\mathbb{D} : \tau_{\leq 2}\underline{SeCat} \longrightarrow \underline{PDer},$$
where we have denoted by $\tau_{\leq 2}B$ the $2$-category obtained from
a $2$-Segal category by replacing all $1$-Segal categories of morphisms
in $B$ by their homotopy categories ($\tau_{\leq 2}B$ is the $2$-Segal analog
of the homotopy category of Segal categories). The following theorem has been proved
in collaboration with D.-C. Cisinski.

\begin{thm}{(Cisinski-To\"en)}\label{t6}
The above $2$-functor
$$\mathbb{D} : \tau_{\leq 2}\underline{SeCat} \longrightarrow \underline{PDer}$$
is fully faithful (in the sense of $2$-categories).
\end{thm}

The above result implies that the theory of pr\'e-d\'erivateurs is
an approximation, \textit{up to $2$-homotopies}, of the theory of Segal categories.
However, the functor $\mathbb{D}$ is surely not essentially surjective
and I personably think that pr\'e-d\'erivateurs not in the essential
image of $\mathbb{D}$ are very unnatural objects which should not
be considered at all.

I would like to finish this part by the proposition below showing that
the higher homotopies that are not taken into account in the theory
of pr\'e-d\'erivateurs are of some importance. For this, we recall that
there is a $K$-theory functor $K^{Wal} : LGWCat \longrightarrow LSp$,
from the Segal category of good Waldhausen categories to the Segal
categories of spectra (see \S 2.1).

\begin{prop}\label{p2}
Let $n$ be any integer.
The morphism of Segal categories
$$K^{Wal} : LGWCat \longrightarrow LSp$$
does not factor, in the homotopy category of Segal categories, through
any Segal category whose simplicial sets of morphisms are $n$-truncated.
\end{prop}

The proof of this proposition relies on the fact that
for a space $X$, the simplicial monoid of auto-equivalences
$aut(X)$ acts naturally on its $K$-theory spectrum $K(X)$
(here $K(X)$ is the space of algebraic $K$-theory of $X$, as defined by
Waldhausen),
and in general this action does not factor through
the $n$-truncation of $aut(X)$. One could also use the same
kind of arguments using the action of the simplicial monoid
of auto-equivalences $aut(A)$ of a simplicial ring $A$ on its
$K$-theory spectra $K(A)$.

An important consequence of the last proposition is the following corollary.

\begin{cor}\label{c3}
Waldhausen $K$-theory functor can not factor, up to a natural
equivalence, through any full sub-$2$-category of $\underline{PDer}$.
\end{cor}

As a consequence of the corollary \ref{c3} we see
that Waldhausen $K$-theory can not be reasonably defined on the level
of \textit{d\'erivateurs triangul\'es} of \cite{ma2}. However,
this does not give a counter example to conjecture $1$
of \cite{ma2} as it is only stated for the $K$-theory of an exact category and not
for a larger class of Waldhausen categories (e.g. including
Waldhausen categories computing $K$-theory of spaces or
of simplicial rings).
In the same way,
Proposition \ref{p2} also implies that there is no reasonable
$K$-theory functor defined on the level of
triangulated categories. \\

The conclusion of proposition \ref{p2} really is: \\

\textit{Waldhausen $K$-theory is not an invariant of $n$-categories for any $n<\infty$.} \\

\section{Segal categories, stacks and homotopy theory}

Segal categories are generalizations of categories, and \textit{Segal topoi}
are to Segal categories what Grothendieck topoi are to categories. The basic notions and
results of the
theory are presented in \cite{to-ve1} and \cite{to-ve3}.
The notion
of Segal topoi appeared naturally at the very beginning of my joint work with
Gabriele Vezzosi, in
our investigation of \textit{homotopical algebraic geometry} (``HAG'' for short). Indeed,
a very natural setting for algebraic geometry is the category of sheaves
of sets on the site of affine schemes, or more generally of $1$-stacks,
$2$-stacks or even $\infty$-stacks (simply called \textit{stacks} in this text).
While developing the basic theory of HAG we discovered that the notion of stacks
over Grothendieck sites is too restrictive for our purposes. Instead,
a notion of \textit{stacks over Segal sites} (i.e. a Segal category
endowed with a suitable notion of Grothendieck topologies) were needed. As topoi
are categories of sheaves, \textit{Segal topoi are Segal categories of stacks
over Segal sites}.

Despite its conceptual interest, the notion of Segal topoi turned out to appear in several
contexts, and seem to be a natural and useful notion (see
\cite{la} for a surprising context of apparition). As an example of
application we have investigated (still with Gabriele Vezzosi)
a reinterpretation and a generalization of Artin-Mazur's homotopy type, which
appear now as part of \textit{higher topos theory}, in the same
spirit as A. Grothendieck's considerations on homotopy types of
topoi found in his letter to L. Breen (see \cite{gr1}).

\subsection{Segal topoi}

\begin{df}\label{d5}
\begin{enumerate}
\item Let $A$ be a Segal category. A Segal category $B$ is a \emph{left exact localization of $A$}
if it is equivalent to a full sub-Segal category $B'$ of  $A$ such that
the inclusion functor $A' \hookrightarrow B$ has a left exact left adjoint.

\item A Segal category $A$ is a \emph{Segal topos} if there exists
a Segal category $T$ such that $A$ is a left exact localization of
$\hat{T}:=\mathbb{R}\underline{Hom}(T^{op},Top)$.

\end{enumerate}
\end{df}

The above definition is based on the fact that a Grothendieck topos is a category which is
a left exact localization of a category of presheaves of sets
(see e.g. \cite{mac-moe,schu}). Let me mention immediately that there exists
Segal topoi which are not left exact localization of $\hat{C}$ for a category
$C$ (see \cite[Rem. 2.0.7]{to-ve1} for a counter-example).
This implies that the fact that $T$ is a Segal category in definition \ref{d5} (2)
and not just a category is crucial. Actually, the Segal topoi
appearing in HAG are not exact localization of $\hat{C}$
for a category $C$. \\

\begin{df}\label{d6}
A \emph{Segal topology on a Segal category} $T$ is a Grothendieck
topology on the homotopy category $Ho(T)$.
A Segal category together with a Segal topology is called a \emph{Segal site}.
\end{df}

When a Segal category $T$ is endowed with a topology $\tau$ one can define
a notion of hyper-coverings which generalizes the usual notion (see \cite[Def. 3.3.2 (1)]{to-ve1}). More precisely, one says that
a morphism $f : F \longrightarrow F'$ in $\hat{T}$ is
a ($\tau$-)epimorphism if for any object $t \in T$ and
any $x\in \pi_{0}(F'(x))$, there is a covering sieve
$S$ of $t \in Ho(T)$ such that for any $u\rightarrow x$
belonging to $S$ there exists $y \in \pi_{0}(F(u))$
with $f(y)=x_{|u}$. Now, a morphism $f : F \longrightarrow F'$
will be called a ($\tau$-)hypercovering if
for any integer $n\geq 0$ the natural morphism
$$F \longrightarrow
F'\times_{(F')^{\partial \Delta^{n}}}F^{\partial \Delta^{n}}$$
is a ($\tau-$)epimorphism in $\hat{T}$. Here, $F^{K}$ denotes the
exponentiation of $F$ by an object $K \in Top$, which is uniquely determined by
the usual adjunction formula\footnote{Such exponentiation exists
in any Segal category with arbitrary limits.}
$$\hat{T}_{(G,F^{K})} \simeq Top_{(K,\hat{T}_{(G,F)})}.$$

For an object $G \in \hat{T}$, we will say that $F$ satisfies the
descent condition for hypercoverings if for any
hypercovering $F \longrightarrow F'$ in $\hat{T}$ the natural
morphism
$$\hat{T}_{(F',G)} \longrightarrow \hat{T}_{(F,G)}$$
is an equivalence in $Top$. This decent condition is the
Segal analog of the usual sheaf and the stack conditions.

Concerning the terminology,
the Segal category of pre-stacks on a Segal site $(T,\tau)$ is
the Segal category $\hat{T}=\mathbb{R}\underline{Hom}(T^{op},Top)$, and the
Segal category of stacks on $(T,\tau)$ is the full sub-Segal category
of $\hat{T}$ consisting of pre-stacks having the descent property for
hyper-coverings (see \cite[Def. 3.3.2 (2)]{to-ve1}). The Segal category of stacks is
denoted by $T^{\sim,\tau}$.

An important result states that the Segal category $T^{\sim,\tau}$ is the localization
of a natural model category of stacks. Indeed, in \cite{to-ve3} is constructed
a model category of stacks on the Segal site $(T,\tau)$, which is denoted by
$SPr_{\tau}(T)$. It is a direct consequence of the strictification
theorem \ref{tstrict} that $LSPr_{\tau}(T)$ is equivalent to $T^{\sim,\tau}$.
In particular, one sees that the Segal category $T^{\sim,\tau}$ possesses all
kind of limits and colimits. In practice the existence of this
model structure provides a more friendly setting for stacks over Segal sites. \\

The main theorem of Segal topos theory is the following statement that
relates Segal topoi and Segal categories of stacks. It is a generalization
of the classical correspondence between Grothendieck topologies and
exact localizations of categories of pre-sheaves. For this we need
the extra notion of \textit{t-complete Segal categories} for which we refer
to \cite[Def. 3.3.6]{to-ve1} for a detailed definition (a Segal topos
is $t$-complete if every hyper-covering is contractible).

\begin{thm}{(To\"en-Vezzosi, \cite[Thm. 3.3.8]{to-ve1})}\label{t7}
Let $T$ be a Segal category.
\begin{enumerate}
\item For any topology $\tau$ on $T$, the inclusion morphism $T^{\sim,\tau} \hookrightarrow \hat{T}$
possesses a left exact left adjoint.

\item The map $\tau \mapsto T^{\sim,\tau}$, which associates to a topology $\tau$
the full sub-category of stacks on $(T,\tau)$ induces a bijective correspondence
between topologies on $T$ and t-complete full sub-Segal categories $A$ on $\hat{T}$ which
inclusion functor $A \hookrightarrow \hat{T}$ possesses a left exact left adjoint.

\item The t-complete Segal topoi as precisely the Segal categories of stacks
over a Segal site.

\end{enumerate}
\end{thm}

The interest of this last result lies in the fact that it justifies our
notion of topologies on Segal categories, at least when one is
dealing with t-complete Segal topoi. It seems possible
however to drop the t-complete assumption by replacing the notion of topologies
by a weaker notion of \textit{hyper-topology}. Informally, in a topology one
fixes the data of coverings whereas in an hyper-topology one fixes
the data of hyper-coverings (it is worth mentioning here that even in the case
of a category, the two notions of topologies and hyper-topologies do not
coincide). We have not investigate this notion
seriously as we did not find any reasons to work with Segal topoi
which are not $t$-complete, and it is not clear at all that this
new notion of hyper-topology has any interest besides
a conceptual one.

Of course theorem \ref{t7} is only the starting point of the whole theory,
and much work has to been done in order to have a workable and powerful
theory of Segal topoi. Unfortunately we have not done much more, except
stating a conjectural Giraud's style characterization of Segal topoi (see \cite[Conj. 5.1.1]{to-ve1}) for which
some recent progress have been made by J. Lurie (see \cite{lu}). \\

To finish this paragraph, let me mention the unpublished work
of C. Rezk on the notion of \textit{homotopy topos} (I personally
prefer the expression \textit{model topos}) which is a model
category analog of the notion of Segal topos (see \cite[Def. 3.8.1]{to-ve3}). By definition,
if $M$ is a model topos then the Segal category $LM$ is
a Segal topos, and furthermore any Segal topos is obtained this
way. However, the main advantage of Segal topoi compare to model
topoi is the existence of a good notion of morphisms.
In other words, given two Segal topoi $T$ and $T'$ there exist
a Segal category of geometric morphisms $\mathbb{R}\underline{Hom}^{geom}(T,T')$,
defined as the full sub-Segal category of $\mathbb{R}\underline{Hom}(T',T)$
consisting of exact morphisms which admit a right adjoint,
which allows one to consider all Segal topoi together assembled in a
$2$-Segal category. These Segal categories of geometric morphisms will be
used in a an essential way in the next paragraph.

\subsection{Homotopy type of Segal topoi}

The starting point of homotopy theory of Segal topoi has been
the following Grothendieck's style interpretation of
homotopy types of spaces. It is also a possible answer to
some conjecture of Grothendieck that appear in
one of his letter to L. Breen (see \cite{to1} for comments on this
conjecture). It uses the notion of
Segal topoi and Segal categories of geometric morphisms between them.
Before stating this result let me recall that given
two Segal topoi $T$ and $T'$ there exists a Segal category
of geometric morphisms $\mathbb{R}\underline{Hom}^{geom}(T,T')$ (it is defined
to be the full sub-Segal category of $\mathbb{R}\underline{Hom}(T',T)$
consisting of exact morphisms which admit a right adjoint). For an
object $p : T \longrightarrow T'$ in $\mathbb{R}\underline{Hom}^{geom}(T,T')$
we will denote by $p^{*}$ the corresponding object in
$\mathbb{R}\underline{Hom}(T',T)$ and by
$p_{*} \in \mathbb{R}\underline{Hom}(T,T')$ its right adjoint.

In order to state the theorem, let us call a stack $F\in T$
in a Segal topos $T$
\emph{constant} if it is in the essential image of $\pi^{*}$, the
inverse image of the unique morphism of Segal topoi
$\pi : T \longrightarrow Top$. A stack
$F$ will be called \emph{weakly locally constant}
\footnote{The expression \emph{locally constant} will
be reserved for a stronger notion that will be introduced
later in the text.} if there
exists an epimorphism $X \longrightarrow *$ in $T$
such that $F\times X \longrightarrow X$ is a constant
stack in the Segal topos $T/X$. When $T=St(X)$ is the topos
of stacks on a topological space $X$, then a stack
$F \in St(X)$ is weakly locally constant if and only if
it is locally equivalent to a constant simplicial presheaves.

\begin{thm}{(To\"en-Vezzosi, \cite[Thm. 5.2.1]{to-ve1})}\label{t8}
\begin{enumerate}
\item For any CW complex $X$ let $St(X)$ be the Segal category of
stacks on $X$. The full sub-Segal category
$Loc(X)$ of $St(X)$, consisting of locally constant stacks on $X$ is a t-complete Segal topos.

\item For two CW complexes $X$ and $Y$, the Segal category of geometric morphisms
$\mathbb{R}\underline{Hom}^{geom}(Loc(X),Loc(Y))$ is a Segal groupoid. Furthermore,
there is a natural equivalence of simplicial sets
$$Map(X,Y)\simeq |\mathbb{R}\underline{Hom}^{geom}(Loc(X),Loc(Y))|,$$
where the right hand side is the nerve of the Segal groupoid $\mathbb{R}\underline{Hom}^{geom}(Loc(X),Loc(Y))$.

\item The morphism of $2$-Segal categories
$$Top \longrightarrow \{Segal \; Topoi\},$$
sending $X$ to the Segal topos $Loc(X)$ is fully faithful.

\end{enumerate}
\end{thm}

The theorem has the following important consequence. It is
a Segal category version of \emph{hypoth\`ese
inspiratrice} of \cite{gr2}, stating that
the homotopy category of spaces does not have any
non-trivial auto-equivalences.

\begin{cor}\label{ct8}
Let $\mathbb{R}\underline{Aut}(A)$ be the full sub-Segal category
of $\mathbb{R}\underline{Hom}(A,A)$ consisting of
equivalences. Then, one has $\mathbb{R}\underline{Aut}(Top)=simeq *$.
\end{cor}

Based on the previous theorem, for any Segal topos $T$ we define its \textit{homotopy shape}
to be the morphism of Segal categories
$$H_{T} : Top \longrightarrow Top,$$
sending $X$ to $|\mathbb{R}\underline{Hom}^{geom}(T,loc(X))|$. The key observation,
base on theorem \ref{t8}, is that
if $T=Loc(Y)$ for a CW complex $Y$ then $H_{T}$ is co-represented (in the sense of
Segal categories) by the
homotopy type of $Y$. In the general situation the functor $H_{T}$ is
only co-representable by a pro-object in the Segal category $Top$. Precisely, one can
prove the following pro-representability result. Its proof essentially
relies on the fact that the homotopy shape $H_{T}$ is a
left exact morphism of Segal categories. Because of some
technical difficulties we will assume that $T$ is a $t$-complete Segal
topos, but I expect the proposition to be correct in general (this complication
is related with the problem of defining a reasonable notion of
hyper-topologies as explained after theorem \ref{t7}).

\begin{prop}{(To\"en-Vezzosi, see \cite{to-ve1})}\label{p3}
Let $T$ be a $t$-complete Segal topos and $H_{T}$ its homotopy shape as defined above.
There exists a left filtered Segal category $A$ (e.g. which possesses finite limits) and
a morphism $K_{T} : A \longrightarrow Top$, which co-represents the morphism $H_{T}$.
In other words, the following two endomorphisms of the Segal category $Top$
$$X \mapsto H_{T}(X) \qquad X \mapsto Hocolim_{a\in A^{op}}Map(K_{T}(a),X)$$
are equivalent.
\end{prop}

In the last proposition there is a universe issue that is not
mentioned. In fact, if $T$ is a $\mathbb{U}$-Segal topos (i.e. the
Segal category $T$ in definition \ref{d5} can be chosen to be
$\mathbb{U}$-small) then the Segal category $A$ can be chosen to be
$\mathbb{U}$-small.

\begin{df}\label{d7}
The morphism of Segal categories $K_{T} : A \longrightarrow Top$ of proposition \ref{p3} is called
the \emph{pro-homotopy type of the Segal topos $T$}.
\end{df}

The pro-homotopy type $K_{T}$ of a Segal topos $T$ is not a
pro-object in the category of simplicial set. Indeed,
it is only a pro-object in the sense of Segal categories (i.e.
the category of indices $A$ is Segal category), but it
seems however that this notion is not strictly more general than
the usual notion of pro-simplicial sets (see \cite[App. B]{lu}). \\

The Segal category $A$ of definition \ref{d7} is left filtered, and
it follows that so is its homotopy category $Ho(A)$. In particular,
the morphism $K_{T} : A \longrightarrow Top$ gives rise to a pro-object
$Ho(K_{T}) : Ho(A) \longrightarrow Ho(Top)$, which is a pro-object in the homotopy
category of spaces. When the Segal topos $T$ is in fact
the Segal category of stacks over a locally connected Grothendieck site $C$,
I suspect that the pro-object $Ho(K_{T}) : Ho(A) \longrightarrow Ho(Top)$ is
isomorphic to the Artin-Mazur's pro-homotopy type of $C$ as defined
in \cite{ar-ma} (at least after some
$\sharp$-construction). The pro-homotopy type $K_{T}$ of definition \ref{d7} is therefore
a refinement of Artin-Mazur's construction, already for the case
of locally connected Grothendieck sites.

The fundamental property of the pro-homotopy type $K_{T}$ of a Segal topos
$T$ is the following.
In order to state it let us first start by some
general notions. We fix a $t$-complete Segal topos $T$ and we simply
let $K$ be its
pro-homotopy type of definition \ref{d7}.

\begin{itemize}

\item Let $\underline{Top}$ be the constant Segal stack over $T$
associated with the Segal category $Top$ (the theory of Segal
stacks of \cite{hi-si} generalizes in an obvious way to Segal
stacks over Segal topoi). The \emph{Segal category of locally constant
stacks on $T$} is defined to be
$$Loc(T):=\mathbb{R}\underline{Hom}(*,\underline{Top}),$$
where the right hand side is the Segal category of
(derived) morphisms of Segal stacks over $T$. I warn the reader that
the natural morphism $Loc(T) \longrightarrow T$, given by
descent theory, is not fully faithful in general. The Segal category $Loc(T)$
is therefore not a full sub-Segal category of $T$ in general, and its
objects consist of objects of $T$ endowed with certain additional
structures. Also, the notion of locally constant stack
on $T$ is clearly different from the notion of
weakly locally constant stacks used before. Not even every weakly locally
constant stack in $T$ lies in the essential image of
the morphism $Loc(T) \longrightarrow T$. However, when $T$ is
locally contractible
(in particular for $T=St(X)$ for a CW complex $X$),
then $Loc(T) \longrightarrow T$ is fully faithful
and its image consists precisely of all weakly locally constant stacks.

\item Let $K : A \longrightarrow Top$ be the pro-homotopy type
of the Segal topos $T$. Composing with the morphism $X \mapsto Loc(X)$
one gets an ind-object in the $2$-Segal category of Segal categories
$$\xymatrix{A^{op} \ar[r]^-{K} & Top^{op} \ar[r]^-{Loc(-)} & \underline{SeCat}.}$$
The colimit of this morphism (in the $2$-Segal category $\underline{SeCat}$)
is denoted by
$$Loc(K):=Colim_{a\in A^{op}}Loc(K_{a}) \in \underline{SeCat},$$
and is called the \emph{Segal category of locally constant
stacks over $K$}.

\end{itemize}

The above unpublished theorem is the universal property of the pro-homotopy type
$K$ of the Segal topos $T$.

\begin{thm}\label{p3'}
With the above notations,
there exists a natural equivalence of Segal categories
$$Loc(T)\simeq Loc(K).$$
\end{thm}

The theorem \ref{p3'} gives a universal property of the
pro-homotopy type associated to a Segal topos, which is
a generalization of the well known universal property of
the fundamental groupoid of a locally connected topos.
When $T$ is the Segal topos of stacks over a locally connected
Grothendieck site, this universal property were not known
(and is probably uneasy or even impossible to state) for
Artin-Mazur's homotopy type, as the definition of $Loc(K)$ really uses
the fact that $K$ is a pro-object in $Top$ and not only in its
homotopy category.

\section{Homotopical algebraic geometry}

The main references for this section are
\cite{to-ve1,to-ve3,to-ve4,to-ve5,to-ve6}. \\

Developing homotopical algebraic geometry (``HAG'' for short) is a project
we started together
with G. Vezzosi during the fall 2000, and which is still in progress.
The main goal of HAG is to provide a mathematical
setting in which one can talk about schemes in a context were
affine objects are modelled by homotopy-ring-like objects (e.g. commutative
differential graded algebras, commutative ring spectra, symmetric monoidal
categories \dots). We already know that algebraic geometry possesses
generalization to a relative setting for which affine objects are
modelled by commutative rings in a general ringed topos or
in a Tannakian category (see \cite{de,hak}). The new feature
appearing in HAG is the fact that the category of models for
affine objects comes with a non-trivial homotopy theory
(e.g. a model category structure) that have to be taken into account.
Clearly, \textit{homotopical algebraic geometry
is to algebraic geometry what homotopical mathematics are to
mathematics}.

The original motivations for starting such a project were various. As an indication
for the reader, I remember below three of them (see also \cite{to-ve4,to-ve5} for more details).

\begin{itemize}
\item \textit{Tannakian duality:} Form the algebraic point of view the Tannakian dual of
a (neutral) Tannakian category (over a field $k$) is a commutative Hopf $k$-algebra, which from the
geometrical point of view corresponds to an affine group scheme over $k$.
In the (conjectural) Tannakian formalism for Segal categories
(see \S 5) the base monoidal category
of $k$-vector spaces is replaced by the monoidal Segal category of complexes of $k$-vector
spaces. The Tannakian dual of a (neutral) Tannakian Segal  category is therefore
expected to be \textit{a commutative Hopf algebra in the monoidal Segal category of complexes}
(i.e. some kind of commutative dg-Hopf $k$-algebra), and geometrically one would like
to consider this Hopf algebra as an \textit{affine group scheme over the Segal category
of complexes of $k$-vector spaces}. My original interests in this
notion of Tannakian Segal category were based on the observation that
several interesting derived categories appearing in algebraic
geometry (e.g. derived category of perfect complexes
of local systems or l-adic sheaves, of perfect complexes with
flat connections, of perfect complexes of F-isocrystals \dots)
are in fact the homotopy categories of natural
Tannakian Segal categories, whose dual can be considered
as certain homotopy types (in the same way as the dual of
a neutral Tannakian category is considered
as a fundamental group, or more generally as a $1$-truncated homotopy type).
This point of view, which can be avoided to actually construct
these homotopy types (see \S 6), seems to me
very powerful for the study of homotopy types in algebraic geometry
(see \S 6.4).

\item \textit{Derived algebraic geometry:} There are essentially two kind of
general constructions in the category of schemes, colimits (e.g. quotients)
and limits (e.g. fibered products). These two constructions are not
\textit{exact} in some sense and according to a very general philosophy they
should therefore be \textit{derived}. Stacks and algebraic stacks theory
has been introduced in order to be able to make \textit{derived quotients}.
More generally, higher stacks provide a theory in which one can do arbitrary
\textit{derived colimits}. On the other side,
the notion of \textit{dg-schemes} (see e.g. \cite{cio-kap1,cio-kap2}) have been
introduced in order to be able to do \textit{derived fibered products}, and more
generally \textit{derived limits}. However, this approach have encountered
two major problems, already identified
in \cite[0.3]{cio-kap2}.

\begin{enumerate}

\item The definition of dg-schemes and dg-stacks seems too
rigid for certain purposes. By definition, a dg-scheme is a space obtained by
\textit{gluing commutative differential graded algebras for the Zariski topology}. It seems
however that certain constructions really require a weaker notion of gluing, as for example
\textit{gluing differential graded algebras up to quasi-isomorphisms} (and a weaker topology).

\item The notion of dg-schemes is not very well suited with respect to the functorial point of view,
as representable functors would have to be defined on the derived category of dg-schemes
(i.e. the category obtained by formally inverting quasi-isomorphisms of dg-schemes),
which seems difficult to describe and to work with.
As a consequence, the derived moduli spaces constructed in \cite{kap1,cio-kap1,cio-kap2}
do not arise as solution to natural
\textit{derived moduli problems}, and are constructed in a rather ad-hoc way.

\end{enumerate}

The main idea to solve these two problems was to interpret dg-schemes
as \textit{schemes over the category of complexes}. Therefore,
it appeared to us that the theory of dg-schemes should
be only an approximation of what
algebraic geometry over the category of complexes is. Such a theory
actually did provide to a us a context in which wa have been
able to construct many new derived moduli spaces that were not
constructed (and probably could not be constructed) as dg-schemes
or dg-stacks.

\item \textit{Moduli spaces of multiplicative structures:} Fixing a
finite dimensional vector space $V$, one can define the classifying stack
of algebra structures on $V$, $Alg_{V}$. This is an algebraic stack in the sense
of Artin, which is a solution to the classification problem
of algebra structures on $V$.

In algebraic topology, the classification problem of ring structure
on a given spectra $X$ appears naturally and in several interesting
contexts (see \cite{goe-hop,laz}). However, this classification problem has, until now,
only been solved in a rather crude way by using classifying spaces
as in \cite{re}. These classifying spaces are homotopy types, and therefore
are discrete invariants (in our previous example they would correspond
to the space of global sections of the stack $Alg_{V}$ alone, thus it
looks like the underlying set of points of an algebraic variety).
It seems therefore
very natural to look for additional algebraic structures on the classifying
space of multiplicative structures on a given spectra $X$, reflecting
some global geometry. The main idea
is be to define a \textit{classifying stack $Alg_{X}$}, of algebra structures
on $X$, in a very similar manner than the stack $Alg_{V}$ is defined. Our
point of view was that the stack $Alg_{X}$ only exist in a reasonable sense
as a \textit{stack over the category of spectra}, and therefore belongs
to \textit{algebraic geometry over spectra}. As we will see
below, HAG actually provide
a context in which this construction, and several others,
makes sense.

\end{itemize}

In the following paragraph I will present the general theory, as well
as some examples of constructions of moduli spaces in this new context.
Details can be found
in \cite{to-ve3,to-ve4,to-ve5}.

\subsection{HAG: The general theory}

The general formalism was inspired to us by the work
of C. Simpson around the notion of $n$-geometric stacks, as
exposed in \cite{si6}. The only new difficulty here is to take into
account correctly the homotopy theory of the category of affine objects, which
is possible thanks to our work on Segal sites, stacks and Segal topoi
presented before. \\

The starting point is a symmetric monoidal model category $M$ which will be
the base category of the theory
(se also remark \ref{rd9}). One considers the category
$Comm(M)$, of $E_{\infty}$-algebras in $M$, which is very often
a model category for which equivalences are simply defined
on the underlying objects in $M$
(I assume that it is for simplicity, see \cite{sp} for details). The
model category of \textit{affine schemes over $M$} is defined
to be the opposite model category $Aff_{M}:=(Comm(M))^{op}$, and
its Dwyer-Kan localization $LAff_{M}$ is called the \emph{Segal category of
affine schemes over $M$}.
We assume that one is given a topology $\tau$ on the Segal
category $LAff_{M}$ (i.e. a Grothendieck topology on the
homotopy category $Ho(Aff_{M})$). Of course, finding interesting
topologies on $LAff_{M}$ very much depends on the context and could be sometimes
not so easy (examples will be given in the next paragraphes).
The Segal category $LAff_{M}$ together with the topology
$\tau$ is a Segal site in the sense of Def. \ref{d6}, and we can therefore apply the
general theory of Segal topoi and stacks  in order to
produce the Segal topos of stacks for the topology $\tau$, denoted by $LAff_{M}^{\sim,\tau}$.
The Segal topos $LAff_{M}^{\sim,\tau}$
is the fundamental object in order to \emph{do algebraic geometry
over $M$}, and plays exactly the same role as the
category of sheaves on the big site of schemes in classical
algebraic geometry.
As explained in \cite[\S 4]{to-ve3}, there is a natural model category
of stacks $Aff_{M}^{\sim,\tau}$ whose simplicial localization
is equivalent to $LAff_{M}^{\sim,\tau}$. In practice the existence
of the model category $Aff_{M}^{\sim,\tau}$ is very helpful as it
allows to work within a model category rather than in a Segal category
and reduces statements for Segal category theory to usual category theory.
Because of this, very often, objects in $LAff_{M}^{\sim,\tau}$ will be
implicitly considered as objects of the
model category $Aff_{M}^{\sim,\tau}$. This has the advantage that stacks will be represented
by actual functors (in the usual sense)
$$F : Comm(M) \longrightarrow SSet,$$
from the category of $E_{\infty}$-algebras in $M$ to the category of simplicial sets,
and which satisfies the following two conditions:
\begin{itemize}
\item The functor $F$ preserves equivalences.

\item The functor $F$ has the descent property with respect to
hyper-coverings (see \cite[\S 4.4]{to-ve3} for details).
\end{itemize}

As for any Segal category, one has
a Yoneda embedding morphism
$$h : LAff_{M} \hookrightarrow \hat{LAff_{M}}.$$
For an $E_{\infty}$-algebra $R \in Comm(M)$, its image by the morphism $h$
is usually denoted by $\mathbb{R}Spec\, R$.
We will assume that the topology is \textit{sub-canonical}, or equivalently
that the Yoneda embedding factors through the sub-Segal category of stacks
$$h : LAff_{M} \hookrightarrow LAff_{M}^{\sim,\tau}\subset \hat{LAff_{M}}.$$

Objects of $LAff_{M}^{\sim,\tau}$ lying in the essential image of $h$ will be called
\textit{representable stacks}, and will play the role of affine schemes in the whole
theory. Directly related to this, one defines the notion of representable
morphisms between stacks. A morphism $f : F \longrightarrow G$ in
$LAff_{M}^{\sim,\tau}$ is \textit{representable} if
for any representable stack $H$, and any morphism $H \rightarrow G$, the
fiber product (in the Segal category $LAff_{M}^{\sim,\tau}$, or equivalently
the homotopy fibered product in the model category
$Aff_{M}^{\sim,\tau}$) $F\times_{G}H$ is a representable stack.

Let us now assume that one has a notion $P$ of morphisms
in the Segal category $LAff_{M}$. We suppose that morphisms in $P$ are stable
by compositions and pull-backs. We also assume some compatibilities
between morphisms in $P$ and the topology $\tau$ that I will not
explicitly state here. For example,
we assume that a morphism which locally (for the topology $\tau$) is in $P$ lies itself
in $P$. I will not describe all of these conditions and the reader is
advised to make the comparison with the case were $\tau$ is the usual
\'etale topology on affine schemes and $P$ consists of all smooth
or flat morphisms.
As usual, the notion of being in $P$ can be extended from
morphisms in $LAff_{M}$ to representable morphisms in $LAff_{M}^{\sim,\tau}$.
Precisely, a representable morphism $f : F \longrightarrow G$ is in $P$
if for any representable stack $H$ the induced morphism
$F\times_{G}H \longrightarrow H$ is in $P$ (this makes sense
since both stacks are representable and therefore this last
morphism can be considered as a morphism in $LAff_{M}$. We use
here the fact that the topology $\tau$ is sub-canonical.).

The main definition of HAG is the following. It is the obvious extension of the
notion of algebraic stacks to our general setting.

\begin{df}\label{d9}
A stack $F \in LAff_{M}^{\sim,\tau}$ is called \emph{$P$-geometric}
if it satisfies the following conditions.
\begin{enumerate}
\item The diagonal morphism $\Delta : F \longrightarrow F\times F$
is representable.

\item There exists a family of representable stacks $\{U_{i}\}$,
and an epimorphism
$$p : \coprod_{i} U_{i} \longrightarrow F$$
such that each $p_{i} : U_{i} \longrightarrow F$ is in $P$ (note that
each $p_{i}$ is
automatically representable by condition (1)).

\end{enumerate}
\end{df}

Without going into details (it would be too long) let me mention that
the above definition can be iterated in order to
define the notion of \textit{$n-P$-geometric stacks} for an integer $n$.
I refer to \cite{si6} and \cite[\S 3.3]{to-ve5} for some details
on this notion.
The very general definition will also appear in \cite{to-ve6}. \\

The general theory of geometric stacks can then be pursued
in parallel with the usual theory of schemes and
of algebraic stacks. One can for example
define notions of sheaves of modules, cohomology, $K$-theory \dots. Almost all
general notions available for algebraic stacks are also available for
geometric stacks in the sense of definition \ref{d9}. Of course, all of the
general theory is very formal and non-trivial mathematical statements and
constructions really start when one specializes the base model category
$M$. In the following I will gives three example, corresponding to the case
were $M$ is the model category of negatively graded complexes,
of unbounded complexes and of spectra. Most of the proofs will appear in
\cite{to-ve6}.

\begin{rmk}\label{rd9}
\emph{The starting point in the theory presented above were
a base symmetric monoidal model category $M$. There also exists
a} model category free \emph{approach, for which the starting point
is a symmetric monoidal Segal category $S$ (see Def. \ref{dtan1} for definitions).
The two point of views can be compared using the fact that
$LM$ is a symmetric monoidal Segal category when $M$ is
a symmetric monoidal model category.
This is of course a much more general theory but which requires
a lot of technology (essentially linear algebra in
monoidal Segal categories, as algebras, modules, \dots) in order to be
done. Ultimately, working entirely in the Segal setting and not refereing
to model categories would have great advantages, but it seems
to that Segal category theory is not well developed enough
in order for such an approach to be reasonable.}
\end{rmk}

\subsection{DAG: Derived algebraic geometry}

For introduction to derived algebraic geometry I refer to \cite{to-ve4}, in
which we have tried to exposed the general motivations, the
main philosophy and part of the history of the subject.  \\

In this part we specialize $M$ to be the model category
$C^{-}(k)$ of negatively graded (with increasing differential)
complexes of $k$-modules, for $k$ a commutative ring. It is a
symmetric monoidal model category for the projective model
structure of \cite[\S 2.3]{ho}, and one can therefore apply the general
constructions presented in the last part. The theory presented
below gives an alternative to the theory of dg-schemes of \cite{cio-kap1,cio-kap2},
having the advantage of providing a functorial point of view
which seems difficult to deal with inside the theory of dg-schemes.
It is also general enough in order to deal with objects as \textit{dg-stacks} and
\textit{higher dg-stacks} for which the approach of dg-schemes does not seem
to be always well suited. As a consequence, we are able to
consider moduli problems that does not seem to have been
constructed before, and in particular we construct
global counter parts of the formal moduli spaces
considered in \cite{kon,kon-soi}. \\

First of all, one defines an \'etale topology on the Segal category
$LAff_{C^{-}(k)}$ of affine schemes over $C^{-}(k)$ in the following way
(recall that the category $Aff_{C^{-}(k)}$ has been defined
as the opposite of the category of $E_{\infty}$-algebras in
the model category $C^{-}(k)$, or in other words of
non-positively graded $E_{\infty}$-algebras over $k$).
A morphism $f : A \longrightarrow B$ of $E_{\infty}$-algebras (over $k$)
is \textit{\'etale} (resp. \textit{strongly smooth},
resp. \textit{strongly flat}) if it
satisfies the following three conditions.
\begin{itemize}
\item The morphism $A \longrightarrow B$ is finitely presented in the
sense of Segal categories\footnote{
Recall that an object $x$ in a Segal category $A$ is finitely
presented if the morphism $A_{(x,-)} : A \longrightarrow Top$
commutes (up to equivalence)  with filtered colimits. A morphism
$x \rightarrow y$ in $A$ is finitely presented if
$y$ is a finitely presented as an object in the coma
Segal category $x/A$.}. Equivalently, this means that
$B$ is equivalent to a retract of a finite cell $A$-algebra
in the sense of \cite{ekmm} for example.

\item The induced morphism $H^{0}(A) \longrightarrow H^{0}(B)$ is an
\'etale (resp. smooth, resp. flat) morphism of rings.

\item The natural morphism $H^{*}(A)\otimes_{H^{0}(A)}H^{0}(B) \longrightarrow H^{*}(B)$
is an isomorphism.
\end{itemize}

Now, a family of morphisms of $E_{\infty}$-algebras $\{A \rightarrow B_{i}\}_{i\in I}$
is an \textit{\'etale covering} if each morphism $A \rightarrow B_{i}$ is \'etale,
and if the family $\{Spec\, H^{0}(B_{i}) \rightarrow Spec\, H^{0}(A)\}_{i \in I}$
is an \'etale covering of affine schemes. This defines an \'etale topology
on the Segal category $LAff_{C^{-}(k)}$, and therefore
a Segal topos of stacks $LAff_{C^{-}(k)}^{\sim,et}$. Furthermore, by
taking $P$ to be the set of strongly smooth morphisms in the sense above, one
obtains the notion of \textit{strongly geometric stacks in $LAff_{C^{-}(k)}^{\sim,et}$}
(these are called \textit{strongly geometric $D$-stacks} in \cite{to-ve4}).

I will not recall here all of the properties of stacks and
strongly geometric stacks. Let me however mention two
important facts.

\begin{itemize}

\item Let $LAff_{k}^{\sim,et}$ be the Segal category of stacks
on the Grothendieck site of affine $k$-schemes with the \'etale topology.
The natural inclusion functor from the category
of commutative $k$-algebras inside the category of $E_{\infty}$-algebras
induces a fully faithful morphism of Segal categories
$$i : LAff_{k}^{\sim,et} \longrightarrow LAff_{C^{-}(k)}^{\sim,et},$$
which admits a right adjoint
$$h^{0} : LAff_{C^{-}(k)}^{\sim,et} \longrightarrow LAff_{k}^{\sim,et}.$$
However, the morphism $i$ is not left exact (as
push-outs of commutative $k$-algebras do not coincide
in general with homotopy push-outs of $E_{\infty}$-algebras), and therefore
the adjoint pair $(i,h^{0})$ does not define a geometric morphism of Segal topoi.
The fact that $i$ does not commutes with fiber products is one of the
key feature of derived algebraic geometry: \emph{taking fiber products
of schemes might not be a scheme anymore}.

For any stack $F \in LAff_{C^{-}(k)}^{\sim,et}$, the stack $h^{0}(F)$ is called
the \textit{truncation of $F$}.
When $F$ is a strongly geometric stack, then $h^{0}(F)$ is an algebraic
stack in the sense of Artin, and the natural morphism $ih^{0}(F) \longrightarrow F$
is a \textit{closed embedding}. This is the general picture of a \textit{classical moduli
space} $h^{0}(F)$ sitting inside its \textit{derived version} $F$.
The closed embedding $ih^{0}(F) \hookrightarrow F$
behave very much as a formal thickening , and the geometry of the
two spaces are essentially the same whereas their structural sheaves
can be very different.

\item For a stack $F \in LAff_{C^{-}(k)}^{\sim,et}$ one can define its tangent stack $TF$
to be the stack of morphisms from $Spec\, k[\epsilon]$ to $F$. It comes
with a natural projection $TF \longrightarrow F$, whose fiber over
a global point $x$ of $F$ is defined to be the \textit{tangent space
of $F$ at $x$}. When the stack $F$ is strongly geometric,
its tangent space $T_{x}F$ at a point $x$ is a \textit{linear stack}. This means that
$T_{x}F$ correspond to a complex of $k$-modules, which is furthermore concentrated
in degree $[-1,\infty[$.

One should be careful that for a stack $F \in LAff_{k}^{\sim,et}$, and in particular
for a scheme, one does not have $iTF\simeq T(iF)$. Indeed, if $X$ is a scheme
then $iTX$ is the global space of the tangent sheaf on $X$ (i.e.
$iTX=Spec\, Symm(\Omega_{X}^{1})$), whereas $T(iX)$ is the global space
of Illusie's tangent complex, $T(iX)=\mathbb{R}Spec\, Symm(\mathbb{L}_{X})$.

\end{itemize}

We are now ready to define various derived moduli functors. As an examples
I will describe the stack $\mathbb{R}Ass$, classifying the associative algebra structures.
Other examples are described in \cite{to-ve4}.

Let $V$ be a projective $k$-module of finite type. For an $E_{\infty}$-algebra
$A$, one considers the category $Alg_{V}(A)$ whose objects are
associative $A$-algebras whose underlying $A$-module is, locally for the
\'etale topology on $A$, equivalent
to $A\otimes_{k}V$, and whose morphisms are equivalences of $A$-algebras.
The simplicial set $\mathbb{R}Ass_{V}(A)$ is defined to be the geometric realization
of the category $Alg_{V}(A)$. For a morphism of $E_{\infty}$-algebras
$A \longrightarrow A'$ one has a base change morphism
$\mathbb{R}Ass_{V}(A) \longrightarrow \mathbb{R}Ass_{V}(A')$
sending a $A$-algebra $B$ to the $A'$-algebra $B\otimes_{A}^{\mathbb{L}}A'$.
This defines a functor
$$\mathbb{R}Ass_{V} : Comm(C^{-}(k)) \longrightarrow SSet,$$
satisfying the conditions for being a stack,
and therefore and object in $LAff_{C^{-}(k)}^{\sim,et}$.

\begin{thm}{(To\"en-Vezzosi)}\label{t11}
The stack $\mathbb{R}Ass_{V}$ is strongly geometric. The tangent
space of $\mathbb{R}Ass_{V}$ at a global point $B \in
\mathbb{R}Ass_{V}(k)$,
corresponding to an associative
algebra structure on $V$, is naturally identified with
the shifted complex of derived derivation $\mathbb{R}Der(B,B)[1]$.
\end{thm}

The strongly geometric stack $\mathbb{R}Ass_{V}$ has a truncation $h^{0}(\mathbb{R}Ass_{V})$, which
is the usual Artin stack $Ass_{V}$ of associative algebra structures on the $k$-module $V$.
The stack $Ass_{V}$ sits as a closed sub-stack inside $\mathbb{R}Ass_{V}$. This is
a very general situation, for any strongly geometric stack $F$ the truncation $h^{0}(F)$ lives
in $F$ as a closed sub-stack. Furthermore, the usual tangent space
of $Ass_{V}$ at the point $B$ can be identified with
the truncation $\tau_{\leq 0}\mathbb{R}Der(B,B)[1]$, explaining, at
least infinitesimally, the additional information contained in
the derived version $\mathbb{R}Ass_{V}$.
The main advantage of
the derived version $\mathbb{R}Ass_{V}$ is that its tangent spaces can be explicitly
identified as the complexes of (derived) derivations, whereas
the cotangent complex of the truncation $Ass_{V}$ remains a mystery
(though is tangent space is understood).
This is, as far as I understand, one of the main feature of
derived moduli spaces; their tangent complexes have very natural
geometric descriptions which allows one to describe the local nature
of the moduli space, whereas the deformation theory of the truncated version
is very much unnatural and difficult to deal with. This observation
is apparently the main point of the \textit{derived deformation theory program} of
P. Deligne, V. Drinfeld and M. Kontsevich.

Without going into details, let me mention that one can define many other
derived moduli problems (e.g. local systems on a space, linear category structures,
vector bundles on a schemes \dots)
and prove that their are strongly geometric stacks in
many cases (see \cite[\S 5]{to-ve4}). The very general point of view of our approach
has also allowed us to define certain moduli stacks that have never been
considered before (as for example the $2$-geometric stack of linear categories
described in \cite[\S 5.3]{to-ve4}).

Clearly, derived algebraic geometry requires several basic
and fundamental results in order to be as easy to manipulate
as usual algebraic geometry. My feeling is that derived
geometric stacks behave in a quite similar fashion
than schemes and several (if not all) of the
results of the EGA's could be generalized to the
derived setting, as for example cohomology of projective spaces,
cohomology and base changes, the local structures of \'etale and smooth
morphisms, Grothendieck's existence theorem \dots.
Quite recently (see \cite{lucats2}), J. Lurie has announced a version
of Artin's representability theorem in the context of
derived algebraic geometry, which is definitely going to
be extremely useful for the construction of
more examples of geometric derived stacks.

\subsection{UDAG: Unbounded derived algebraic geometry}

In this part I want to present a natural extension of derived algebraic geometry
of the last paragraph, in which the base model
category of bounded above complexes $C^{-}(k)$ is replaced
with the model category of
unbounded complexes $C(k)$. This extension allows one to define and study
more moduli problems that do not naturally enter the setting of the
last paragraph, as for example the generalization of the stack
$\mathbb{R}Ass_{V}$ when $V$ is a complex of $k$-modules rather than just a $k$-module
(i.e. the classifying stack of dga structures on a complex $V$). Another
natural example is the stack of fiber functors of a $k$-tensor Segal category
mentioned in the motivations of this
chapter, and which will be investigated in \S 5. \\

Let now $M$ be the model category of unbounded complexes $C(k)$ with the
projective model structure, which is
a symmetric monoidal model category (see \cite[\S 2.3]{ho}).
The category $Aff_{C(k)}$ is now the opposite category of
unbounded $E_{\infty}$-algebras over $k$.
One defines the \'etale
topology on $LAff_{C(k)}$ exactly the same way as for the bounded case (note
that the definition given in the last paragraph
is also valid for unbounded $E_{\infty}$-algebras). Associated to this
topology is the Segal category of stacks $LAff_{C(k)}^{\sim,et}$.

The inclusion functor $i : C^{-}(k) \longrightarrow C(k)$ induces an adjunction
on the Segal categories of stacks
$$i_{!} : LAff_{C^{-}(k)}^{\sim,et} \longrightarrow LAff_{C(k)}^{\sim,et}
\qquad LAff_{C^{-}(k)}^{\sim,et} \longleftarrow LAff_{C(k)}^{\sim,et} : i^{*},$$
and the morphism $i_{!}$ is in fact fully faithful. The notion of stacks over
unbounded $E_{\infty}$-algebras is therefore a generalization of the notion
of stacks over bounded above $E_{\infty}$-algebras. In particular,
there exists natural fully faithful morphisms
$$\xymatrix{LAff_{k}^{\sim,et} \ar[r] &
LAff_{C^{-}(k)}^{\sim,et} \ar[r] &  LAff_{C(k)}^{\sim,et},}$$
each of them having right adjoints. This allows one to consider
schemes, algebraic stacks and strongly geometric stacks as defined
in the last paragraph as objects in $LAff_{C(k)}^{\sim,et}$.  \\

It turns out that the notion of strong smoothness used in the last
paragraph is not going to be
general enough for our new examples. Indeed, we will be lead to take
quotient by representable group stacks in $LAff_{C(k)}^{\sim,et}$, and
a representable stack which is strongly smooth over $Spec\, k$
is automatically a group scheme in the usual sense. We therefore
extend the notion of strongly smooth morphisms in the following way (this is what is
called \textit{standard smoothness} in \cite[\S 4.4]{to-ve4}). It is based on
the usual characterization of smooth morphisms of schemes as the morphisms
which \'etale locally looks like the projection of a vector bundle.

A morphism of $E_{\infty}$-algebra $f : A \longrightarrow B$ is
\textit{smooth} if there exists a commutative square (up to homotopy)
$$\xymatrix{A \ar[r]^-{f} \ar[d]_-{u} & B \ar[d]_-{v} \\
 A' \ar[r]_-{f'} & B',}$$
such that
\begin{itemize}
\item $v$ is an \'etale covering.

\item $f'$ is an \'etale morphism.

\item The $A$-algebra $A'$ is equivalent to
$A\otimes_{k}^{\mathbb{L}}L(E)$, where $E$ is a perfect complex
of $k$-modules and $L(E)$ is the free $E_{\infty}$-algebra
over $E$.
\end{itemize}

By taking $P$ to be the set of smooth morphisms as defined above one obtains
by definition Def. \ref{d9} the notion of a \textit{geometric stack} in $LAff_{C(k)}^{\sim,et}$.
One should be careful however that an object $F \in LAff_{C^{-}(k)}^{\sim,et}$
can be geometric as an object in $LAff_{C(k)}^{\sim,et}$ without being strongly geometric. To clarify, one has natural inclusions
$$\hspace*{-0.6cm} \left\lbrace Schemes \right\rbrace \hookrightarrow
\left\lbrace Alg. \; stacks \right\rbrace \hookrightarrow
\left\lbrace Strong. \; geom. \; stacks \right\rbrace
\hookrightarrow
\left\lbrace Geom. \; stacks \right\rbrace \hookrightarrow
LAff_{C(k)}^{\sim,et}$$

We can now generalize the construction of the stack $\mathbb{R}Ass$ in the following way.
Let $V$ be a perfect complex of $k$-modules.
For a $E_{\infty}$-algebra
$A$, one considers the category $Alg_{V}(A)$ whose objects are
associative $A$-algebras whose underlying $A$-module is, locally for the
\'etale topology on $A$, equivalent
to $A\otimes_{k}^{\mathbb{L}}V$, and whose morphisms are equivalences of $A$-algebras.
The simplicial set $\mathbb{R}Ass_{V}(A)$ is defined to be the geometric realization
of the category $Alg_{V}(A)$. For a morphism of $E_{\infty}$-algebras
$A \longrightarrow A'$ one has a base change morphism
$\mathbb{R}Ass_{V}(A) \longrightarrow \mathbb{R}Ass_{V}(A')$
sending a $A$-algebra $B$ to the $A'$-algebra $B\otimes_{A}^{\mathbb{L}}A'$.
This defines a functor
$$\mathbb{R}Ass_{V} : Comm(C(k)) \longrightarrow SSet,$$
satisfying the stack conditions,
and therefore and object in $LAff_{C(k)}^{\sim,et}$.

\begin{thm}{(To\"en-Vezzosi)}\label{t12}
The stack $\mathbb{R}Ass_{V}$ is geometric. The tangent
space of $\mathbb{R}Ass_{V}$ at a global point $B$, corresponding to a dga
structure on $V$, is naturally identified with
the shifted complex of derived derivation $\mathbb{R}Der(B,B)[1]$.
\end{thm}

\begin{rmk}\label{r5}
\emph{The geometric stack $\mathbb{R}Ass_{V}$ also has a truncation
$h^{0}\mathbb{R}Ass_{V}=:Ass_{V}$, which is a stack in
$LAff_{k}^{\sim,et}$. This stack is not
an algebraic stack in the sense of Artin anymore, and is actually
a non-truncated stack in general. However, one can show that
$Ass_{V}$ is equivalent to a quotient stack
$[F/G]$, where $F$ is an affine stack in the sense
of \cite{to2}, and $G$ is an affine group stack acting
on $F$. This shows that the stack $Ass_{V}$ still has
some kind of algebraic nature (see \cite[4.2.1]{to2}).}
\end{rmk}

It is sometimes useful to consider an associative $k$-algebra
as a $k$-linear category with a unique object. In the same way,
a dga can be considered as a dg-category (over $k$) with a unique object.
It would be too long to explain in details, but one can define
a stack $\widetilde{RCat_{V}}$ which classifies
structures of dg-categories having a unique object and
$V$ as complex of endomorphisms (see \cite[\S 5.3]{to-ve4}). This stack can not be geometric, but
one can show it is a \textit{$2$-geometric stack}.
Considering a dga as a dg-category
gives a morphism of stacks $\mathbb{R}Ass_{V} \longrightarrow \widetilde{RCat_{V}}$, which
is a smooth fibration. Indeed, the fiber of this morphism over
the image of dga $B$ can be identified with $K(Gl_{1}(B),1)$, the classifying stack
of the group stack of invertible elements in the dga $B$ (see \cite[\S 4.3]{to-ve4} for more details).
Therefore, the above projection morphism induces an exact triangle
on the level of tangent spaces, that can be written as
$$\xymatrix{ B[1] \ar[r] & \mathbb{R}Der_{k}(B,B)[1] \ar[r] &
C^{+}_{k}(B,B)[2], \ar[r]^-{+1} & }$$
where $C^{+}_{k}(B,B)$ is the complex of Hochschild cohomology of the dga $B$.
This is the fundamental triangle appearing in \cite[p. 59]{kon}. The fibration sequence
(in the Segal category $LAff_{C(k)}^{\sim,et}$)
$$K(Gl_{1}(B),1) \longrightarrow \mathbb{R}Ass_{V} \longrightarrow \widetilde{\mathbb{R}Cat_{V}}$$
gives a geometric explanation of the triangle \cite[p. 59]{kon}. The object
$\widetilde{\mathbb{R}Cat_{V}}$ seems to right think to consider
in order to explain deformation theory of associative dg-algebras
and dg-categories in the spirit of \cite{kon-soi,soi}.

Finally, the $2$-geometric stack of dg-categories $\widetilde{\mathbb{R}Cat}$
can have some interesting applications
in the context of mirror symmetry which some people think about
as some involution of the stack $\widetilde{\mathbb{R}Cat}$
suitably modified (one needs to consider dg-categories
which are \emph{non-commutative Calabi-Yau manifolds}
in the sense of \cite{soi}).

\subsection{BNAG: Brave new algebraic geometry}

The expression \textit{brave new algebra} was introduced
by F. Waldhausen. As far as I understand it reflects the fact
that ring spectra behave very much the same way as
ordinary rings, and that many basic constructions in the context
of rings (tensor product, $Tor$, $Ext$, algebra
of matrices, $K$-theory, Hochschild cohomology, \dots)
has reasonable analogs in the context of ring spectra (see \cite{vo} for
an introduction on the subject). We like to use the expression
\textit{brave new algebraic geometry} for the special case of
HAG in which the base model category $M$ is chosen to be the category
of spectra. One can say that \textit{brave new algebraic geometry is to
algebraic geometry what brave new algebra is to algebra}.

The setting of brave new algebraic geometry seems to be a very well suited setting for many
recent works in stable homotopy theory, as elliptic cohomology (see \cite{goe}),
moduli spaces of multiplicative
structures (see \cite{goe-hop,laz}) \dots. In this part I will present very basic constructions, as
we did not investigate yet more serious examples
(see however \cite{to-ve7} for some more details). \\

We let $M$ to be the model category $Sp^{\Sigma}$, of symmetric spectra as defined in \cite{ho-sh-sm}. It is
a symmetric monoidal model category for the smash product. We endow
the Segal category $LAff_{Sp^{\Sigma}}$
(the opposite Segal category of $E_{\infty}$-ring spectra\footnote{An very convenient
model for the homotopy theory of $E_{\infty}$-ring spectra is exposed in \cite{shi}.}) with an \'etale topology
defined in the following way. A morphism of $E_{\infty}$-ring spectra
$A \longrightarrow B$ is \textit{\'etale} if it satisfies the following
properties.
\begin{itemize}
\item The $A$-algebra is finitely presented (in the sense
of Segal categories).
\item The cotangent spectra $\mathbb{L}_{B/A}$ is acyclic.
\end{itemize}

A finite family of morphisms of $E_{\infty}$-ring spectra $\{A \longrightarrow B_{i}\}_{i \in I}$
is an \textit{\'etale covering} if each morphisms $A \longrightarrow B_{i}$ is \'etale, and if
the family of base change functors of the category of modules
$$\begin{array}{ccc}
Ho(Mod(A)) & \longrightarrow & Ho(Mod(B_{i})) \\
M & \mapsto & M\wedge_{A}^{\mathbb{L}}B_{i}
\end{array}$$
is a conservative family (i.e. a $A$-module $M$ is acyclic if and only if each
$B_{i}$-modules $M\wedge_{A}^{\mathbb{L}}B_{i}$ is acyclic).

This defines an \'etale topology on $LAff_{Sp^{\Sigma}}$,
and therefore one has an associated
Segal category of stacks $LAff_{Sp^{\Sigma}}^{\sim,et}$. The existence of this
Segal category already allows to make sense of certain constructions
that have been informally considered by J. Rognes (see \cite{ro}).
First of all, for a $E_{\infty}$-ring spectra $R$, one can consider
its small \'etale site $(\mathbb{R}Spec\, R)_{et}$, as well as its
associated Segal topos. By the general construction \ref{d7} one extract from
this Segal topos a pro-homotopy type, called the \textit{\'etale homotopy type
of the $E_{\infty}$-ring spectra $R$}. Following the same kind of ideas, one
can define \'etale $K$-theory of an $E_{\infty}$-ring spectra, $K_{et}(R)$,
which comes naturally with a localization morphism $K(R) \longrightarrow K_{et}(R)$
(see \cite{to-ve3}).
One could then try to follow the same guide line as for usual rings,
like defining a Thomason's style topological $K$-theory, state a Quillen-Lichtenbaum
conjecture \dots. As far as I understood from \cite{ro} and talks
by J. Rognes, one of the main goal of the whole machinery
would be to compute, at least partially, the algebraic $K$-theory of the
sphere spectra which appears naturally in geometric topology
(the sphere spectra $\mathbb{S}$ is a direct factor
of the algebraic $K$-theory of the point space $*$).

To finish this part, let me mention that one can also define
a notion of smooth morphisms of $E_{\infty}$-ring spectra,
in the exact same way that we have defined it for
$E_{\infty}$-algebras. Using this definition one obtains a useful notion of
geometric stacks in $LAff_{Sp^{\Sigma}}$ (though there are some
complications as the \'etale topology is not known to be sub-canonical).
The very first basic example of geometric stack is the spectra
version of $\mathbb{R}Ass_{V}$ presented in the previous paragraph, but where $V$ is now
a finite $\mathbb{S}$-module. I will not repeat the definition here as it
is very similar to the one in the linear context. We obtain this way a
stack $\mathbb{R}Ass_{V}$, classifying associative ring structures on
the spectra $V$, and one can prove that $\mathbb{R}Ass_{V}$ is a geometric stack
(see \cite{to-ve7}).
This geometric stack gives a way to approach geometrically similar
questions than the one stated in \cite{goe-hop}.

\section{Tannakian duality for Segal categories}

My wish to generalize Tannakian duality to the more
general setting of Segal categories comes back
to the spring 1999 while I was reading \emph{pursuing stacks}, and have been
my original interest in higher category theory.
My main objective was to understand in what sense
affine group schemes duals to Tannakian categories
are \emph{fundamental groups}, as is usually considered
in the literature. In other words, my original question
was \emph{of which kind of 'spaces' are they
the fundamental groups ?} This question was motivated by
the \emph{schematization problem} exposed in \cite{gr1}, and
the general feeling that the pro-algebraic completion of the
fundamental group of topological space $X$ is only
the $\pi_{1}$ of a more general object, called the
\emph{schematization}, and encoding higher homotopical information.
In this section I will present the general formalism
of Tannakian Segal categories, for which several
applications and examples will be discusses in the
next two sections. \\

The starting point of the Tannakian theory for Segal categories is
on one hand the well known analogies between Galois and Tannakian
duality, and on the other hand the generalization of
Galois theory to the Segal setting exposed in \S 3.2.
I
personably like to keep in
mind the following informal scheme
$$
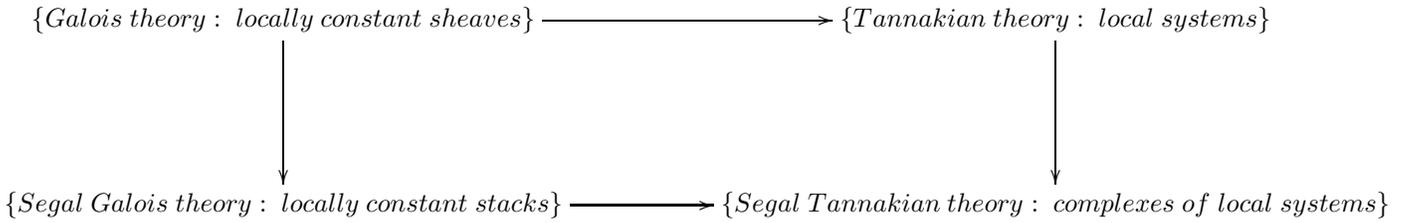
\hspace{-2cm}
\begin{footnotesize}
\xymatrix{
\left\lbrace  Galois \; theory : \;
locally \; constant \; sheaves  \right\rbrace \ar[rr] \ar[dd]& &
\left\lbrace Tannakian \; theory : \;
local \; systems \right\rbrace \ar[dd] \\
 & & \\
\left\lbrace Segal\; Galois \; theory : \;
locally \; constant \; stacks \right\rbrace \ar[rr] & &
\left\lbrace Segal \; Tannakian \; theory : \;
complexes\; of \; local\;  systems \right\rbrace}
\end{footnotesize}$$

In this last diagram, the vertical arrows represent
the passage to homotopical mathematics (i.e. from categories to
Segal categories), whereas the horizontal ones represent
the linearization process (i.e. passing from sets
to vector spaces and from homotopy types to
complexes of vector spaces). The diagram also says that
the coefficients of Galois theory are locally constant sheaves of sets,
whereas in the Tannakian theory they are local systems of
vector spaces. In the Segal setting the coefficients of Galois
theory are locally constant stacks, and in the Tannakian theory are
complexes of local systems (or rather complexes whose cohomology
are local systems). \\

In this section I will describe the general notions
of Segal version of the Tannakian formalism, as
tensor structures, fiber functors \dots. I will also state
the main theorem of Tannakian duality for Segal categories,
which unfortunately remains a conjecture. The situation
here is a bit frustrating as it seems no new ideas
are really required in order to prove it, and that the
standard arguments used in the proof of the usual
Tannakian duality should also make sense
for Segal categories. However, certain technical difficulties
concerning homotopy coherences appeared during the application of
Beck's theorem for Segal categories and have prevented me
to find a complete proof.

The material of this section is extracted from
\cite{to3}, with some slight modifications. For example,
the present version of the stack of fiber functors is
defined over all commutative ring spectra, whereas
in \cite{to3} I was only considering its restriction to
commutative algebras, which does not seem enough
in non-zero characteristic. \\

\subsection{Tensor Segal categories}

The fundamental notion is the one of symmetric monoidal
Segal category, already briefly used in \S 2.5, and more generally
the $2$-Segal category of symmetric monoidal categories.
It can be defined as follows.

Let $\underline{SeCat}$ be the $2$-Segal category of
Segal categories, as defined for example in \cite[\S ]{hi-si}.
We consider $\Gamma$, the category of pointed finite sets
and pointed morphisms between them, which is considered
as a $2$-category and therefore as a $2$-Segal category.
Let $\mathbb{R}\underline{Hom}(\Gamma,\underline{SeCat})$ be
the $2$-Segal category of (derived) morphisms of $2$-Segal categories
from $\Gamma$ to $\underline{SeCat}$, defined
using the model category structure on the category of
$2$-Segal precategories of \cite[\S 2, \S 11]{hi-si}.

For any integers $1 \leq n$, we consider
the finite set $[n]:=\{0,\dots,n\}$, pointed at $0$, as an object
of $\Gamma$. One has for any $1\leq i\leq n$ the
pointed morphisms
$$s_{i} : \{0,\dots,n\} \longrightarrow \{0,1\}$$
sending everything to $0$ and $i$ to $1$. This defines, for
any morphisms of $2$-Segal categories
$A : \Gamma \longrightarrow \underline{SeCat}$,
a well defined morphism of Segal categories
$$\prod_{i}s_{i} : A([n]) \longrightarrow A([1])^{n},$$
called the $n$-th Segal morphism (compare with \cite[\S 2]{hi-si}).

\begin{df}\label{dtan1}
The \emph{$2$-Segal category of symmetric monoidal Segal categories}
(\emph{$\otimes$-Segal categories} for short)
is defined to be the full sub-$2$-Segal category of
$\mathbb{R}\underline{Hom}(\Gamma,\underline{SeCat})$
consisting of morphisms $A : \Gamma \longrightarrow \underline{SeCat}$
such that
\begin{itemize}
\item $A([0])\simeq *$
\item for any integer $n$ the natural morphism
$A([n]) \longrightarrow A([1])^{n}$ is an equivalence of
Segal categories.
\end{itemize}
The $2$-Segal category of $\otimes$-Segal categories
will be denoted by $\otimes-\underline{SeCat}$. For two
$\otimes$-Segal categories $A$ and $B$, the Segal
category of symmetric monoidal morphisms (\emph{$\otimes$-morphisms}
for short) will be denoted by\footnote{
Here, $A_{(x,y)}$ denotes the Segal category
of morphisms between two objects $x$ and $y$ in a $2$-Segal category $A$.}
$$\mathbb{R}\underline{Hom}^{\otimes}(A,B):=
\otimes-\underline{SeCat}_{(A,B)}.$$
Finally, the underlying Segal category of a $\otimes$-Segal category $A$ is
$A([1])$, and will be still denoted by $A$.
\end{df}

The situation here is very much the same as the definition
of Segal categories, and the monoidal structure is encoded in the diagram
of Segal categories
$$\xymatrix{A([2]) \ar[r] \ar[d] & A([1]) \\
A([1])\times A([1]), & }$$
where the vertical arrow is an equivalence. This shows the existence
of a morphism $-\otimes - : A\times A \longrightarrow A$, well
defined in the homotopy category of Segal categories, which is the
first step of the structure of a $\otimes$-Segal category $A$.
Of course, this morphism is not enough to recover the whole
$A$, and higher coherencies are encoded in the whole
diagram over $\Gamma$. In any case, we will often use the
morphism $-\otimes - : A\times A \longrightarrow A$.

Another important remark is that when $A$ is
a $\otimes$-Segal category then its homotopy category $Ho(A)$
has a natural symmetric monoidal structure. Indeed, as
the functor $A \mapsto Ho(A)$ commutes with finite products, the
composite functor
$$Ho(A) : \xymatrix{\Gamma \ar[r] & \underline{SeCat} \ar[r]^-{Ho(-)} &
\underline{Cat} \hookrightarrow \underline{SeCat}}$$
still satisfies the condition to be a $\otimes$-Segal category. Of course,
the structure of $\otimes$-Segal category $A$ is more than just
a symmetric monoidal structure on the homotopy category
$Ho(A)$, as it encodes also higher homotopy coherences (in particular
for the commutativity). \\

Definition \ref{dtan1} is purely in terms of Segal categories, and in
practice it is very useful to have a more down to earth
description of the $2$-Segal category $\otimes-\underline{SeCat}$. This
is possible thanks to the $2$-Segal version of the
strictification theorem \ref{tstrict}, showing that
$\otimes-\underline{SeCat}$ can be described in the following way.
\begin{itemize}
\item Its objects are functors
$$A : \Gamma \longrightarrow SeCat$$
from $\Gamma$ to the (usual) category of Segal categories,
which satisfy the following two conditions.
\begin{itemize}
\item $A([0])=*$

\item For any $n\geq 1$, the natural morphism
$A([n]) \longrightarrow A([1])^{n}$ is an equivalence.

\end{itemize}

\item Recall the existence of the model category
of Segal precategories, $PrCat$ (see \cite{hi-si,to-ve1}).
$\otimes$-Segal categories can therefore be considered
as $\Gamma$-diagrams in $PrCat$, for which
there exists a well known model category structures
(equivalences and fibrations are defined levelwise,
see \cite{hir}).

Given two $\otimes$-Segal categories
$A$ and $B$, corresponding to two functors
$A,B : \Gamma \longrightarrow SeCat$, the
Segal category of $\otimes$-morphisms
$\mathbb{R}\underline{Hom}^{\otimes}(A,B)$
is equivalent to the Segal category of morphisms
$\underline{Hom}(QA,RB)$, where $QA$ is a cofibrant replacement
of the diagram $A$ and $RB$ is a fibrant replacement of the
diagram $B$.

\end{itemize}

\begin{rmk}\label{rtan1}
\begin{enumerate}
\item
\emph{The previous description of the $2$-Segal category
$\otimes$-Segal categories is actually a consequence of the
existence of a model category of $\otimes$-Segal precategories,
as described in \cite{to3}.}
\item \emph{Our $\otimes$-Segal categories are generalization of
}$ACU-\otimes$-categories \emph{in the sense of \cite{sa}. We could
also defined various notions of $\otimes$-Segal categories corresponding
to various combination of the constraints $A$ (associativity), $C$
(commutativity) and $U$ (unital). In the Segal setting the commutativity
constraint is the most interesting, because it has an infinity of
non-equivalent generalizations, the notion of $d$-fold $\otimes$-Segal
categories briefly encountered in \S 2.5. In some sense, our present notion
corresponds to the case $d=\infty$, and is the most commutativity condition
one could define.}
\end{enumerate}
\end{rmk}

There are many interesting examples of $\otimes$-Segal categories.

\begin{itemize}

\item As any category is a Segal category,
any $ACU-\otimes$ category ($\otimes$-category for short),
in the sense of \cite{sa} is a $\otimes$-Segal category.
Indeed, if $C$ a $\otimes$-category one can construct
a pseudo-functor $\Gamma \longrightarrow Cat$, sending
$[n]$ to $C^{n}$ and using the monoidal structure to define
transition functors (this construction is discussed in great details
in \cite{le2}). Strictifying this pseudo-functor gives a
functor $\Gamma \longrightarrow Cat$, and therefore
a functor $\Gamma \longrightarrow SeCat$ which is easily seen
to satisfy the condition to be a $\otimes$-Segal category.

Using this construction, we will always consider
$\otimes$-categories as $\otimes$-Segal categories. Furthermore,
this defines a natural morphism of $2$-Segal categories
$$\otimes-\underline{Cat} \longrightarrow \otimes-\underline{SeCat}$$
which is fully faithful. The essential image of this
morphism clearly consists of all $\otimes$-Segal categories whose
underlying Segal category is equivalent to a category.

\item Let $M$ be a symmetric monoidal model category, in the sense
of \cite[\S 4.3]{ho}. Then, the Segal category $LM$ has a natural
structure of a $\otimes$-Segal category. It can be described in the
following way. We consider $M^{c,1}$, the full sub-category of
$M$ of consisting of cofibrant objects and the unit
(which might be non-cofibrant). Then, $M^{c,1}$ is a category with a
notion of equivalences and a compatible monoidal structure.
This implies that the corresponding functor
$\Gamma \longrightarrow Cat$ is in fact a functor
from $\Gamma$ to the category of pairs $(C,W)$ consisting of
a category $C$ and a sub-category $W \subset C$. Composing with the
construction $(C,W) \mapsto L(C,W)$ one gets
a functor $\Gamma \longrightarrow SeCat$ which is
a $\otimes$-Segal category. Furthermore, the underlying Segal category
is $LM^{c,1}$, which is easily seen to be equivalent
to $LM$ via the natural embedding $M^{c,1} \hookrightarrow M$.

This construction gives a lot of interesting $\otimes$-Segal categories,
as for example the Segal categories of complexes and of symmetric spectra.

\item There is a natural morphism of $2$-Segal
categories $\Pi_{\infty} : Top \longrightarrow \underline{SeCat}$
sending a simplicial set to its Segal fundamental groupoid
(see \cite[\S 2]{hi-si}, where
our $\Pi_{\infty}$ is denoted by
$\Pi_{1,se}$), and which identifies $Top$ with the
full $2$-Segal sub-category of $\underline{SeCat}$ consisting of
Segal groupoids. Using this identification, one sees that
there is an equivalence between the $2$-Segal category of $\Gamma$-spaces
and the $2$-Segal category
of $\otimes$-Segal categories whose underlying Segal category is
a Segal groupoid. Therefore, this provides an equivalence between the
$2$-Segal category of connective spectra and the $2$-Segal category of
$\otimes$-Segal groupoids.

\end{itemize}

We now come to the cental definition of a tensor Segal category,
an analog of the notion of tensor category in the Segal setting.

\begin{df}\label{dtan2}
A \emph{tensor Segal category} is
a $\otimes$-Segal category $A$ satisfying the following two conditions.
\begin{itemize}
\item The underlying Segal category $A$ is stable (in the
sense of Def. \ref{d1}).

\item For any object $x \in A$, the morphism of Segal categories
$x\otimes - : A \longrightarrow A$ is exact.

\end{itemize}
For any tensor Segal categories $A$ and $B$, we denote by
$\mathbb{R}\underline{Hom}^{\otimes}_{ex}(A,B)$ the full sub-Segal category
of $\mathbb{R}\underline{Hom}^{\otimes}(A,B)$ consisting of exact morphisms.
Objects in $\mathbb{R}\underline{Hom}^{\otimes}_{ex}(A,B)$ will be called
\emph{tensor morphisms from $A$ to $B$}.

Tensor Segal categories and $\mathbb{R}\underline{Hom}^{\otimes}_{ex}$ form
a $2$-Segal category denoted by $\underline{TenSeCat}$.
\end{df}

Tensor Segal categories essentially comes from the following example.
Let $M$ be a symmetric monoidal model category. We also
assume that the model category $M$ is stable (in the sense of
\cite[\S 7]{ho}). Then, the $\otimes$-Segal category $LM$ is
a tensor Segal category. Indeed, clearly $LM$ is stable Segal category
(see Thm. \ref{t3} $(4)$). Furthermore the existence
of the adjunction equivalence $Map_{M}(x\otimes^{\mathbb{L}} y,z)\simeq
Map_{M}(x,\mathbb{R}\underline{Hom}(y,z))$, implies that
$x\otimes^{\mathbb{L}} -$ commutes with homotopy colimits
in $M$. This is equivalent to condition $(2)$ of Def. \ref{dtan2}. \\

\begin{rmk}
\emph{As a tensor Segal category $A$ is on one hand
a stable Segal category, and the other hand a $\otimes$-Segal
category, the homotopy category $Ho(A)$ is
endowed with a natural triangulated structure and a symmetric
monoidal structure. Of course, these two structures are compatible,
making $Ho(A)$ into a tensor triangulated
category.}
\end{rmk}

Let us now consider the category of symmetric spectra
$Sp^{\Sigma}$, endowed with its positive
model category structure described in \cite[Prop. 3.1]{shi}.
The model category $Sp^{\Sigma}$ has a natural
monoidal structure, $-\wedge -$, making into a
symmetric monoidal model category.
Let $A$ be a commutative (unital and associative) monoid
in $Sp^{\Sigma}$ (we will simply say that $A$ is a
commutative ring spectrum), and $A-Mod$ be its category of
modules. From \cite[Thm. 3.2]{shi} we know that
$A-Mod$ is a symmetric monoidal model category. Furthermore,
this model category is clearly stable. Therefore, the
Segal category $L(A-Mod)$ is then a tensor Segal category.
The tensor Segal categories $L(A-Mod)$ are very important, as they
are the Segal version of the tensor categories of modules over some
commutative rings, and actually they behave even better as
their monoidal structure is always exact. Furthermore,
it is important to notice that
if $k$ is a commutative ring, and $Hk$ is the corresponding
Eilenberg-MacLane commutative ring spectrum, the $L(Hk-Mod)$
is equivalent (as a tensor Segal category) to
$LC(k)$, where $C(k)$ is the symmetric monoidal model category
of complexes of $k$-modules.

Let $A \longrightarrow B$ be a morphism of commutative
ring spectra, one has a base change morphism
$$-\wedge^{\mathbb{L}}_{A}B : L(A-Mod) \longrightarrow L(B-Mod),$$
which is a morphism of tensor Segal categories.
If we denote by $Comm(Sp^{\Sigma})$ the model category of
commutative ring spectra (see \cite[Thm. 3.2]{shi}), this construction defines
a morphism of $2$-Segal categories
$$LComm(Sp^{\Sigma}) \longrightarrow \underline{TensSeCat}.$$
The first fundamental
conjecture is the following.

\begin{conj}\label{conjtan1}
The morphism of $2$-Segal categories
$$\begin{array}{ccc}
L(Comm(Sp^{\Sigma})) & \longrightarrow & \underline{TensSeCat} \\
A & \mapsto & L(A-Mod) \\
(A\rightarrow B) & \mapsto & -\wedge^{\mathbb{L}}_{A}B
\end{array}$$
is fully faithful.
\end{conj}

Of course, the natural approach to prove the conjecture
\ref{conjtan1} is by defining a morphism in the other way, sending
a tensor Segal category $A$ to $A_{(1,1)}$, where $1$ is the unit of the
tensor structure. However, to endow $A_{(1,1)}$ with a natural
structure of commutative ring spectrum requires to solve some
homotopy coherence problems which do not seem so obvious. In any case,
the conjecture seems to me clearly correct.

I would also like to mention that conjecture \ref{conjtan1} is surely false
if one considers the tensor triangulated categories $Ho(A-Mod)$ instead
of the Segal categories $L(A-Mod)$, as it might exists two
non-equivalent commutative ring spectra with equivalent
tensor triangulated categories of modules (similar examples
for non-commutative ring spectra are given in \cite[Example. 3.2.1]{schw-shi}). \\

Another important conjecture is a comparison between the
theories of commutative ring spectra and
of commutative monoids in the tensor Segal category $LSp^{\Sigma}$.
Let $FS$ be the category of finite sets and all morphisms
between them. The disjoint union makes $FS$ into a
symmetric monoidal category and therefore into a
$\otimes$-Segal category. For a $\otimes$-Segal category $A$, the
Segal category of commutative monoids in $A$ is defined to
be
$$Comm(A):=\mathbb{R}\underline{Hom}^{\otimes}(FS,A).$$
In \cite{to3} I have constructed a natural morphisms
of Segal categories
$$L(Comm(Sp^{\Sigma})) \longrightarrow Comm(LSp^{\Sigma}),$$
where $Comm(Sp^{\Sigma})$ is as before the model category of
commutative ring spectra.

\begin{conj}\label{contan2}
The natural morphism
$$L(Comm(Sp^{\Sigma})) \longrightarrow Comm(LSp^{\Sigma})$$
is an equivalence of Segal categories.
\end{conj}

Conjecture \ref{contan2} is the key result in order to have
a Segal category interpretation of the theory of ring spectra.
It has also a sense when $Sp^{\Sigma}$ is replaced by any
symmetric monoidal model category and $Comm(Sp^{\Sigma})$
by the model category of $E_{\infty}$-algebra in $M$, and
can be seen as a monoidal analog of the strictification theorem
Thm. \ref{tstrict}. I also explains why the notion of commutative
ring spectra is the \emph{right one}, as it corresponds
to \emph{commutative monoids in the monoidal $\infty$-category of spectra}. \\

We finish by the Segal version of the notion
of rigidity (see e.g. \cite{br,sa}). \\

\begin{df}\label{dtan3}
Let $A$ be a $\otimes$-Segal category $A$.
\begin{itemize}
\item We say that $A$ is  \emph{closed}  if
for any two objects $x$ and $y$ in $A$, the morphism of
Segal categories
$$\begin{array}{ccc}
A^{op} & \longrightarrow & Top \\
z & \mapsto & A_{(z\otimes x,y)}
\end{array}$$
is representable by an object $\mathbb{R}\underline{Hom}(x,y)\in A$
(in the sense explained in \S 2.1).

\item We say that $A$ is \emph{rigid} if it is closed and furthermore if
for any object $x$ in $A$, the natural morphism
in $Ho(A)$
$$\mathbb{R}\underline{Hom}(x,1)\otimes x \longrightarrow \mathbb{R}\underline{Hom}(x,x)$$
is an isomorphism in $Ho(A)$.
\end{itemize}
\end{df}

An important remark is that when $M$ is a symmetric monoidal model category then
the $\otimes$-Segal category $LM$ is always closed. Indeed, the existence of
derived $Hom$ objects in $M$ (see \cite[\S 4.3]{ho}) implies the existence of
the objects $\mathbb{R}\underline{Hom}(x,y)\in LM$.

In a rigid $\otimes$-Segal category $A$ one always has the
famous formula
$$\mathbb{R}\underline{Hom}(x\otimes x',y\otimes y')\simeq
\mathbb{R}\underline{Hom}(x,y)\otimes \mathbb{R}\underline{Hom}(x',y'),$$
and in particular
$$\mathbb{R}\underline{Hom}(x,y)\simeq \mathbb{R}\underline{Hom}(x,1)\otimes y.$$

\subsection{Stacks of fiber functors}

The right setting to state the Tannakian duality for Segal categories is
BNAG, exposed in \S 4.4. We refer to this paragraph for the general notions, and
we start by defining the Segal version of the well known adjunction between
tensor categories and stacks using the central notion fiber functors. \\

We consider the model category $Comm(Sp^{\Sigma})$ of commutative, associative and
unital monoids in the category of symmetric spectra $Sp^{\Sigma}$. By
\cite[Thm. 3.2]{shi}, we know that $Comm(Sp^{\Sigma})$ is endowed with a model
category structure for which fibrations and equivalences are defined on the level
of underlying spectra. Objects in $Comm(Sp^{\Sigma})$ will simply be called
\emph{commutative ring spectra}. We consider the Segal category
$LComm(Sp^{\Sigma})$ as well as its opposite category that we will denote by
$$LComm(Sp^{\Sigma})^{op}:=LAff_{Sp^{\Sigma}}.$$
We endow the Segal category $LAff_{Sp^{\Sigma}}$ with the
following topology. A finite family of morphisms of commutative
ring spectra $\{A \longrightarrow B_{i}\}_{i\in I}$ is called a \emph{strongly ffqc
covering} if it satisfies the following conditions.
\begin{itemize}
\item The family of morphisms of affine schemes
$\{Spec\, H^{0}(B_{i}) \longrightarrow Spec\, H^{0}(A)\}_{i\in I}$ is
faithfully flat.

\item For any $i$, the natural morphism of $H^{0}(B_{i})$-modules
$H^{*}(A)\otimes_{H^{0}(A)}H^{0}(B_{i}) \longrightarrow
H^{*}(B_{i})$ is an isomorphism.

\end{itemize}

One checks that this defines a Segal topology on the Segal category $LAff_{Sp^{\Sigma}}$ in the
sense of Def. \ref{d6}, and therefore on has an associated Segal category of stacks
$LAff_{Sp^{\Sigma}}^{\sim,sffqc}$ (\emph{sffqc} stands for \emph{strongly faithfully flat and
quasi-compact}). The Segal topology $sffqc$ is sub-canonical and therefore the Yoneda embedding
gives a fully faithful morphism of Segal categories
$$h : LAff_{Sp^{\Sigma}} \longrightarrow LAff_{Sp^{\Sigma}}^{\sim,sffqc}.$$
The image of a commutative ring spectra $A$ by $h$ is as usual denoted by
$\mathbb{R}Spec\, A$. \\

The Segal site $LAff_{Sp^{\Sigma}}$ has a natural stack in tensor Segal categories
$LParf$, which is a BNAG version of the usual stack of vector bundles. It is
defined in the following way.

For any commutative ring spectrum $A$, one has a closed $\otimes$-Segal
category $L(A-Mod)$ of $A$-modules. One can consider the full
sub-Segal category $L(A-Mod)^{rig}$ consisting of \emph{rigid
objects}. Recall here that an object $x$ in a closed $\otimes$-Segal category $T$
(in the sense of definition Def. \ref{dtan3} $(1)$) is \emph{rigid} if the
natural morphism $\underline{Hom}(x,1)\otimes x \longrightarrow \underline{Hom}(x,x)$
is an isomorphism in $Ho(T)$. Rigid objects in $L(A-Mod)$ are precisely the
$A$-modules which are \emph{strongly dualizable} in the sense
of \cite{ekmm}, and therefore consist of all retract of finite
cell $A$-modules. One also sees that rigid objects in $L(A-Mod)$ are
precisely the finitely presented objects (in the sense of Segal category).
The Segal category $L(A-Mod)^{rid}$ is called the Segal category
of perfect $A$-modules, and is denoted by $LParf(A)$. Clearly,
$LParf(A)$ is a rigid tensor Segal category.

For a morphism of commutative ring spectra $A \longrightarrow B$, one has
a base change morphism
$$\begin{array}{ccc}
LParf(A) & \longrightarrow & LParf(B) \\
M & \mapsto & M\wedge^{\mathbb{L}}_{A}B,
\end{array}$$
defining (after some standard strictification procedure)
a morphism of $2$-Segal
categories
$$\begin{array}{ccc}
LAff_{Sp^{\Sigma}}^{op} & \longrightarrow & \underline{TenSeCat}^{rig} \\
A & \mapsto & LParf(A) \\
(A\rightarrow B)  & \mapsto & -\wedge_{A}^{\mathbb{L}}B,
\end{array}$$
from the Segal site $LAff_{Sp^{\Sigma}}^{op}$ to the $2$-Segal category
of rigid tensor Segal categories. This morphism is in fact a \emph{stack in tensor Segal categories},
in the sense that the underlying pre-stack of Segal categories is a stack
(over the Segal site $LAff_{Sp^{\Sigma}}^{op}$).

\begin{df}\label{dtan4}
The \emph{stack of perfect modules} is the stack of tensor Segal categories
defined above. It is denoted by $LParf$.
\end{df}

The stack $LParf$ can be used in order to define
a morphism from the $2$-Segal category of stacks $LAff_{Sp^{\Sigma}}^{\sim,sffqc}$
to the $2$-Segal category $\underline{TenSeCat}^{rig}$ of rigid tensor Segal categories. For
a stack $F \in LAff_{Sp^{\Sigma}}^{\sim,sffqc}$, one can consider the Segal category
$LParf(F):=\mathbb{R}\underline{Hom}(F,LParf)$. This Segal category is
naturally a rigid tensor Segal category, and therefore one has a morphism of $2$-Segal categories
$$\begin{array}{ccc}
(LAff_{Sp^{\Sigma}}^{\sim,sffqc})^{op} & \longrightarrow & \underline{TenSeCat}^{rig} \\
F & \mapsto & LParf(F).
\end{array}$$
This morphism possesses a left adjoint
$$\begin{array}{ccc}
\underline{TenSeCat}^{rig} & \longrightarrow & (LAff_{Sp^{\Sigma}}^{\sim,sffqc})^{op} \\
T & \mapsto & FIB(T),
\end{array}$$
where $FIB(T) \in LAff_{Sp^{\Sigma}}^{\sim,sffqc}$
is the stack of fiber functors on $T$ defined as follows. For any
commutative ring spectra $A$ one defines
$$FIB(T)(A):=|\mathbb{R}\underline{Hom}_{ex}^{\otimes}(T,LParf(A))|,$$
the geometric realization of the Segal category of tensor morphisms from $T$ to
$LParf(A)$\footnote{Note that the rigidity condition on $T$ implies that
$\mathbb{R}\underline{Hom}_{ex}^{\otimes}(T,LParf(A))$ is actually a Segal
groupoid, which justifies considering its geometric realization.}. One has the following adjunction
equivalence for a stack $F \in LAff_{Sp^{\Sigma}}^{\sim,sffqc}$ and a rigid tensor Segal category $T$
$$\mathbb{R}\underline{Hom}(F,FIB(T))\simeq |\mathbb{R}\underline{Hom}_{ex}^{\otimes}(T,LParf(F))|,$$
saying that $T \mapsto FIB(T)$ is the left adjoint to $F \mapsto LParf(F)$. \\

The situation can also be generalized over a base commutative ring spectra
$A$ in the following way. Instead of the $2$-Segal category $\underline{TenSeCat}^{rig}$, of
rigid tensor Segal categories, one consider $A-\underline{TenSeCat}^{rig}$, the $2$-Segal
category $LParf(A)/\underline{TenSeCat}^{rig}$, of rigid tensor Segal categories
under $LParf(A)$. The $2$-Segal category $A-\underline{TenSeCat}^{rig}$ will be
called the \emph{$2$-Segal category of rigid $A$-tensor Segal categories}. Its objects
are simply rigid tensor Segal categories $T$ together with a tensor morphism
$LParf(A) \longrightarrow T$. For two objects
$LParf(A) \longrightarrow T$ and $LParf(A) \longrightarrow T'$ the
Segal category of morphisms from $T$ to $T'$ in $LParf(A)/\underline{TenSeCat}^{rig}$
sits naturally in a homotopy cartesian square
$$\xymatrix{
\mathbb{R}\underline{Hom}_{A-ex}^{\otimes}(T,T'):=A-\underline{TenSeCat}^{rig}_{(T,T')} \ar[rr] \ar[d] & &
\mathbb{R}\underline{Hom}_{ex}^{\otimes}(T,T') \ar[d] \\
\bullet \ar[rr] & & \mathbb{R}\underline{Hom}_{ex}^{\otimes}(LParf(A),T'). }$$
In other words $\mathbb{R}\underline{Hom}_{A-ex}^{\otimes}(T,T')$ is the (homotopy) fiber
of
$$\mathbb{R}\underline{Hom}_{ex}^{\otimes}(T,T') \longrightarrow
\mathbb{R}\underline{Hom}_{ex}^{\otimes}(LParf(A),T')$$
at the point corresponding to the structural morphism
$LParf(A) \longrightarrow T'$.

For a rigid $A$-tensor Segal category $LParf(A) \longrightarrow T$, one can define
its stack of fiber functors $FIB_{A}(T)$, which naturally lives in the Segal category
$LAff_{Sp^{\Sigma}}^{\sim,sffqc}/\mathbb{R}Spec\, A$ of stacks over $\mathbb{R}Spec\, A$\footnote{The Segal category
$LAff_{Sp^{\Sigma}}^{\sim,sffqc}/\mathbb{R}Spec\, A$ is equivalent to the
category of stacks over the Segal site of commutative $A$-algebras.
In other words, with the notations of \S 4
one has
$$LAff_{Sp^{\Sigma}}^{\sim,sffqc}/\mathbb{R}Spec\, A \simeq
LAff_{A-Mod}^{\sim,sffqc}$$
where $A-Mod$ is the model category of $A$-modules, and
$sffqc$ is the induced strongly flat topology defined as
for commutative ring spectra. In particular, objects
in $LAff_{Sp^{\Sigma}}^{\sim,sffqc}/\mathbb{R}Spec\, A$ can be
considered as functors
$$F : A/Comm(Sp^{\Sigma}) \longrightarrow SSet$$
from the category of commutative $A$-algebras to the category
of simplicial sets.}. This defines
a morphism of $2$-Segal categories
$$\begin{array}{ccc}
A-\underline{TenSeCat}^{rig} & \longrightarrow & (LAff_{Sp^{\Sigma}}^{\sim,sffqc}/\mathbb{R}Spec\, A)^{op} \\
T & \mapsto & FIB_{A}(T),
\end{array}$$
which also has a right adjoint
$$\begin{array}{ccc}
(LAff_{Sp^{\Sigma}}^{\sim,sffqc}/\mathbb{R}Spec\, A)^{op} & \longrightarrow & A-\underline{TenSeCat}^{rig} \\
F & \mapsto & LParf_{A}(F):=\left(LParf(A) \rightarrow LParf(F)\right).
\end{array}$$

\begin{df}\label{dtan5}
For a rigid $A$-tensor Segal category $T$, the \emph{stack of fiber functors (over $A$)} is
the stack $FIB_{A}(T) \in LAff_{Sp^{\Sigma}}^{\sim,sffqc}/\mathbb{R}Spec\, A$ defined above.
\end{df}

The adjunction $(T \mapsto FIB_{A}(T),F \mapsto LParf_{A}(F))$ is the fundamental object of study,
and the Tannakian duality for Segal categories describes conditions on stacks and
$A$-tensor Segal categories for which this adjunction induces a one-to-one correspondence. \\

\subsection{The Tannakian duality}

We let $A$ be a base commutative ring spectra. For example $A$ could be
of the form $Hk$, for a commutative ring $k$, which will be our main example
of applications. \\

In this part I will introduce the notions of \emph{$A$-Tannakian Segal categories} and
of \emph{Tannakian (Segal) gerbes over $A$}, as well as
a conjecture stating that the
constructions $F \mapsto LParf_{A}(F)$ and $T \mapsto FIB_{A}(T)$ induces an equivalence
between them. This conjecture is the Segal version of the Tannakian duality theorem
giving an equivalence between affine gerbes and Tannakian categories as stated in
\cite{de,sa}. However, some technical difficulties have prevented me to
actually prove this conjecture, thought I am convinced it is correct. \\

I start with the definition of affine gerbes in the Segal setting. For this I will need
the following notion of morphisms of positive Tor dimension between commutative ring spectra.
For a morphism $f : B \longrightarrow B'$ of commutative ring spectra, one says that
$f$ is of \emph{positive Tor dimension}\footnote{There are no mistakes here, has
for a complex of abelian groups $E$, viewed as a symmetric spectra $HE$, one has
$\pi_{i}(HE)\simeq H^{-i}(E)$.} if for any $B$-module $M$,
$$\left(\pi_{i}(M)=0 \; \forall \; i>0\right) \Rightarrow \left(
\pi_{i}(M\wedge^{\mathbb{L}}_{B}B')=0 \; \forall \; i>0\right).$$
Furthermore, one says that $f : B \longrightarrow B'$ is
\emph{faithfully of positive Tor dimension} if it is of positive Tor dimension
and if moreover the base change functor
$$-\wedge_{B}^{\mathbb{L}}B' : Ho(B-Mod) \longrightarrow
Ho(B'-Mod)$$
is conservative (i.e. $M\simeq *$ if and only if
$M\wedge_{B}^{\mathbb{L}}B'\simeq *$).

If $f : B \longrightarrow B'$ is a strongly flat morphism in the sense
explained in the last paragraph, then one has
$\pi_{i}(M\wedge^{\mathbb{L}}_{B}B')\simeq \pi_{i}(M)\otimes_{\pi_{0}(B)}\pi_{0}(B')$, and therefore
$f$ is of positive Tor dimension. Moreover, if $k \longrightarrow k'$
is a morphism of commutative rings then
$Hk \longrightarrow Hk'$ is of positive Tor dimension
if and only if $k'$ is flat over $k$.
However, the notion of positive Tor dimension morphisms
is much weaker than of strongly flat morphisms, as there exists
positive Tor dimension morphisms $Hk \longrightarrow B'$, where $k$ is commutative ring
but $B'$ has non-trivial negative homotopy groups.

As usual, the notion of positive Tor dimension morphisms extends to the case
of representable morphisms of stacks.

\begin{df}\label{dtan6}
A stack $F \in LAff_{Sp^{\Sigma}}^{\sim,sffqc}/\mathbb{R}Spec\, A$, is
an \emph{affine (Segal) gerbe} if it satisfies the following  two conditions.
\begin{itemize}
\item The stack $F$ is locally non-empty (i.e.
there exists a $sffqc$ covering $A \longrightarrow A'$ such that
$F(A')\neq \emptyset$).

\item For any commutative $A$-algebra $B$, and any two morphisms of stacks
$x,y : \mathbb{R}Spec\, B \longrightarrow F$, the stack of path from $x$ to $y$
$$\Omega_{x,y}F:=\mathbb{R}Spec\, B  \times_{x,F,y}\mathbb{R}Spec\, B$$
is representable and faithfully of positive Tor dimension over $\mathbb{R}Spec\, B$
(i.e. there exists a morphism which is faithfully of positive
Tor dimension $B \longrightarrow C$ such that
$\Omega_{x,y}F\simeq \mathbb{R}Spec\, C$).

\end{itemize}
An affine (Segal) gerbe $F$ is said to be \emph{neutral} if
there exists a morphism of stacks
$* \longrightarrow F$.
\end{df}

Clearly, a general affine gerbe over $A$ is
locally for the positive Tor dimension topology of the form $BG$ for $G$ an affine group stack
of positive Tor dimension over $A$. This last description is much closer to the usual
notion of affine gerbe encountered in algebraic geometry, which are stacks locally equivalent to
classifying stacks of flat affine group schemes (see the next
paragraph for more about the comparison with the usual
notion). \\

We are now ready to define \emph{Tannakian (Segal) gerbes}, which
are affine gerbes satisfying some cohomological
conditions, as well as \emph{Tannakian Segal categories}.

\begin{df}\label{dtan7}
\begin{itemize}
\item A morphism $f : F \longrightarrow F'$ in $LAff_{Sp^{\Sigma}}^{\sim,sffqc}$
is a $P$-\emph{equivalence} if the induced morphism of tensor
Segal categories
$$f^{*} : LParf(F') \longrightarrow LParf(F)$$
is an equivalence.
\item

A stack $F \in LAff_{Sp^{\Sigma}}^{\sim,sffqc}$, is
$P$-\emph{local}, if for any $P$-equivalence of stacks $G \longrightarrow G'$, the
induced morphism
$$(LAff_{Sp^{\Sigma}}^{\sim,sffqc})_{(G',F)} \longrightarrow
(LAff_{Sp^{\Sigma}}^{\sim,sffqc})_{(G,F)}$$
is an equivalence of simplicial sets.

\item A stack $F \in LAff_{Sp^{\Sigma}}^{\sim,sffqc}/\mathbb{R}Spec\, A$, is
a \emph{$A$-Tannakian (Segal) gerbe} if it is an affine gerbe over $A$, and
if the stack $F$ is $P$-local.

\item A $A$-tensor Segal category $T$ is a \emph{$A$-Tannakian Segal category}
if it is equivalent to some $LParf_{A}(F)$ for $F$ a $A$-Tannakian
gerbe. A $A$-tensor Segal category is \emph{neutral} if furthermore
the $A$-Tannakian gerbe $F$ above can be chosen to be neutral.

\end{itemize}
\end{df}

\begin{rmk}\label{rtstructure}
\emph{When $F$ is a $A$-Tannakian Segal gerbe, the tensor Segal category
$LParf_{A}(F)$ comes equipped with some kind of $t$-structure. Indeed,
one can chose a $sffqc$ covering $A \rightarrow A'$ and a point
$x : \mathbb{R}Spec\, A' \longrightarrow F$, and define an object
$E \in LParf_{A}(F)$ to be} positive \emph{if $x^{*}(E)$ is a positive
$A'$-module (i.e. $\pi_{i}(M)=0$ for all $i>0$).
The fact that this definition is independent of the choice
of $x$ uses that any two points of $F$ are locally equivalent up  to a
positive Tor dimension covering.}
\end{rmk}

The main conjecture is the following duality statement.

\begin{conj}\label{tanconj}
\begin{enumerate}

\item For any $A$-Tannakian Segal category $T$, the natural morphism
$$T \longrightarrow LParf_{A}(FIB_{A}(T))$$
is an equivalence of $A$-tensor Segal categories.

\item A rigid $A$-tensor Segal category $T$ is Tannakian if
it satisfies the following two conditions.
\begin{enumerate}
\item The structural morphism
$LParf(A) \longrightarrow T$ induces an equivalence
$$LParf(A)_{(1,1)} \longrightarrow T_{(1,1)},$$
where $1$ denotes the unit of the $\otimes$-structures.
\item There exists a $sffqc$ covering $A \rightarrow A'$,
and a $A$-tensor morphism
$$\omega : T \longrightarrow LParf(A')$$
satisfying the following two conditions.

\begin{enumerate}
\item The extension of $\omega$ to the Segal categories of $Ind$-objects
$$\overline{\omega} : Ind(T) \longrightarrow
Ind(LParf(A'))\simeq L(A'-Mod)$$
is conservative (i.e. $\overline{\omega}(M)\simeq *$ if and only if
$M\simeq *$).

\item Let $Ind(T)_{\geq 0}$ (resp. $Ind(T)_{<0}$) be the full sub-Segal category
of $Ind(T)$ consisting of objects $E$ such that
$\pi_{i}(\overline{\omega}(E))=0$ for all $i\geq 0$ (resp. for all $i<0$).
Then, and object $E \in Ind(T)$ lies in $Ind(T)_{\geq 0}$ if and only if
for any $F \in Ind(T)_{<0}$ one has
$T_{(F,E)}\simeq *$.
\end{enumerate}
\end{enumerate}
\end{enumerate}
\end{conj}

Clearly, conjecture \ref{tanconj} implies that the two morphisms
of $2$-Segal categories $FIB_{A}$ and $LParf_{A}$ induce
an equivalence between the $2$-Segal category of $A$-Tannakian
Segal categories and the $2$-Segal category of $A$-Tannakian
gerbes. \\

\begin{rmk}\label{rtan}
\begin{itemize}
\item
\emph{Condition $(2-b-ii)$ of conjecture \ref{tanconj} defines
some kind of $t$-structure on the stable Segal category
$Ind(T)$, and it could very well be that this additional structure
is an important part of data for Tannakian Segal categories that
I have been neglecting up to now. For example it could very well be
that the equivalence $T \longrightarrow LParf_{A}(FIB_{A}(T))$
of \ref{tanconj} $(1)$ is only correct if one take into account the natural
$t$-structure on $Ind(T)$, and if one replaces
$FIB_{A}(T)$ by its sub-stack of} $t$-positive \emph{fiber functors (i.e.
the one that preserve positive objects). I will however stay with the conjecture
\ref{tanconj} as is, as it seems that for $A=Hk$ and $k$ a field any fiber
functor is automatically $t$-positive.}
\item
\emph{The conjecture \ref{tanconj} is quite general, and
I am not sure it is so important to have it for
a general base commutative ring spectra $A$. Up to now, I only
consider it serious when $A=Hk$ for some commutative ring
$k$.}
\end{itemize}
\end{rmk}

To finish this part I a going to describe the general steps for a proof
of the main point of conjecture \ref{tanconj} in the neutral case. \\

Let me start by some linear algebra notions in the context
of commutative ring spectra.

\begin{itemize}
\item
A $H_{\infty}$-Hopf $A$-algebra is a co-simplicial diagram
of commutative $A$-algebras
$$B_{*} : \Delta \longrightarrow A/Comm(Sp^{\Sigma})$$
such that $B_{0}=A$, and for each $n$ the natural
morphism (the dual version of the $n$-th Segal morphism)
$$\underbrace{B_{1}\wedge_{A}^{\mathbb{L}}B_{1} \dots \wedge_{A}^{\mathbb{L}}B_{1}}_{n\; times}
\longrightarrow B_{n}$$
is an equivalence. Clearly, $H_{\infty}$-Hopf $A$-algebras
correspond via the $\mathbb{R}Spec$ functor to
affine group stacks in $LAff_{Sp^{\Sigma}}^{\sim,sffqc}/\mathbb{R}Spec\, A$.

\item
Any $H_{\infty}$-Hopf $A$-algebra $B_{*}$ has a rigid $A$-tensor Segal category of
perfect (or rigid) comodules $L(B_{*}-Comod^{rig})$, defined
as the limit (of Segal categories)
$$L(B_{*}-Comod^{rig}) \longrightarrow
Holim_{n\in Delta}LParf(B_{n}).$$
Clearly, one has a natural equivalence
$$LParf(B\mathbb{R}Spec\, B_{*})\simeq L(B_{*}-Comod^{rig}),$$
where $B\mathbb{R}Spec\, B_{*}:=Hocolim_{n\in \Delta^{op}}\mathbb{R}Spec\, B_{n}$
is the stack whose loop stack
is the affine group scheme $\mathbb{R}Spec\, B_{1}$.
In the same way, one has a full (non-rigid) $A$-tensor Segal category of
comodules
$$L(B_{*}-Comod) \longrightarrow
Holim_{n\in Delta}L(B_{n}-Mod).$$
\end{itemize}

Now, let $T$ be a rigid $A$-tensor Segal category satisfying the
conditions of \ref{tanconj} $(b)$, but with $A'=A$
(this is the neutral case). We first consider the induced
morphism on the Segal categories of $Ind$-objects
$$\overline{\omega} : Ind(T) \longrightarrow
L(A-Mod),$$
which is a morphism of non-rigid $A$-tensor Segal categories.
By some general principle, the morphism $\overline{\omega}$
have a right adjoint $p$ (this follow from the fact that
it commutes with colimits and that
$Ind(T)$ has a small sets of small generators). We next consider
the $A$-module $B:=\overline{\omega}(p(A))$, which is expected to
have a natural structure
of an $H_{\infty}$-Hopf $A$-algebra. The comultiplication is
given by the adjunction morphism
$$B=\overline{\omega}p(A)
\longrightarrow B\wedge_{A}^{\mathbb{L}}B \simeq \overline{\omega}p\overline{\omega}p(A)$$
induced by $Id \rightarrow p\overline{\omega}$, and the multiplication is given
by using the natural morphism $p(A)\otimes p(A) \rightarrow p(A)$
adjoint to $\overline{\omega}p(A)\otimes \overline{\omega}p(A) \rightarrow
A\wedge_{A}^{\mathbb{L}}A\simeq A$ (Though this structure seems clear
from an heuristic point of view, controlling all the homotopy coherences
is a problem). Furthermore, any object $M \in Ind(T)$ gives
rise to a $B_{*}$-comodule $\overline{\omega}(M)$, via
the co-action
$$\overline{\omega}(M) \longrightarrow
B\wedge_{A}^{\mathbb{L}}\overline{\omega}(M)\simeq \overline{\omega}p\overline{\omega}(M),$$
induced by the adjunction $Id \rightarrow p\overline{\omega}$.
One should show this way that the adjunction $(\overline{\omega},p)$ induces
a new adjunction
$$\overline{\omega} : Ind(T) \longrightarrow L(B_{*}-Comod) \qquad
Ind(T) \longleftarrow L(B_{*}-Comod) : p.$$
A Segal version of the Barr-Beck theorem should now be applied, thanks
to the conservative property of $\overline{\omega}$, which shows that
this last adjunction is an equivalence. This equivalence identifies the
sub-Segal categories of rigid objects.

In conclusion, we have shown that for any rigid $A$-tensor Segal category
$T$ which satisfies the conditions $(b)$ of conjecture \ref{tanconj}
there exists an equivalence between $Ind(T)$
and $L(B_{*}-Comod)$ which transforms the fiber functor $\overline{\omega}$
into the forgetful functor $L(B_{*}-Comod) \longrightarrow L(A-Mod)$.
We can therefore replace $T$ by $L(B_{*}-Comod^{rig})$, which implies that
$T$ is of the form $LParf(F)$ for $F$ the classifying stack
of the group stack $\mathbb{R}Spec\, B$. Finally,
the second assertion on $\overline{\omega}$ will imply that the $A$-algebra
$B=\overline{\omega}p(A)$ is faithfully of positive Tor dimension, showing that
$F$ is an affine gerbe. This would be the main point
for a proof of conjecture \ref{tanconj}, and the remaining statements
should follow quite formally from it. \\

\subsection{Comparison with usual Tannakian duality}

We now suppose that $k$ is a commutative ring, and that
$A=Hk$ is its Eilenberg-MacLane commutative ring spectra\footnote{In
this situation $A$-tensor Segal categories will simply be called
$k$-tensor Segal categories.}.
On one hand one has the Segal topos $LAff_{Sp^{\Sigma}}^{\sim,sffqc}/\mathbb{R}Spec\, Hk$, whose objects can be
seen as functors
$$F : Hk/Comm(Sp^{\Sigma}) \longrightarrow SSet$$
from the category of commutative $Hk$-algebras to the category
of simplicial sets (and satisfying the usual stack
conditions). One the other hand
one has the usual Grothendieck site $(Aff/k,ffqc)$, of affine $k$-schemes
with the faithfully flat topology, and its
associated Segal topos $LAff_{k}^{\sim,ffqc}$. Sending
a commutative $k$-algebra $k'$ as
a commutative $Hk$-algebra $Hk'$, induces a well defined
morphism of Segal categories (it is not a geometric morphism
of Segal topoi)
$$i : LAff_{k}^{\sim,ffqc} \longrightarrow LAff_{Sp^{\Sigma}}^{\sim,sffqc}/\mathbb{R}Spec\, Hk,$$
which on affine objects sends $Spec\, k'$ to $\mathbb{R}Spec\, Hk'$.
This morphism is actually fully faithful, and possesses
a right adjoint
$$h^{0} : LAff_{Sp^{\Sigma}}^{\sim,sffqc}/\mathbb{R}Spec\, Hk
\longrightarrow LAff_{k}^{\sim,ffqc}$$
defined for a stack $F \in LAff_{Sp^{\Sigma}}^{\sim,sffqc}/\mathbb{R}Spec\, Hk$
by the restriction, $h^{0}(F)(k'):=F(Hk')$.

The very first observation is that if $F$ is an affine
gerbe in the sense of \cite{de}, then $i(F) \in LAff_{Sp^{\Sigma}}^{\sim,sffqc}/\mathbb{R}Spec\, Hk$ is
an affine Segal gerbe in the sense of Def. \ref{dtan6}, which is neutral
if and only if $F$ is. This way, affine gerbe in the sense
of \cite{de} form a full sub-Segal category of
the Segal category of affine Segal gerbes. Our notion of
affine Segal gerbe is therefore a reasonable generalization
of the usual notion. \\

Let us now suppose that $k$ is a field, and let
$T$ be a Tannakian category over $k$ in the sense
of \cite{de}. We let $\mathcal{G}$ be the gerbe of fiber functors
on $T$, which is an object in $LAff_{k}^{\sim,ffqc}$.

We consider $C^{b}(T)$, the category of bounded complexes
in $T$, and we let $LParf(T):=LC^{b}(T)$ be the Segal category obtained by localizing $C^{b}(T)$ along the quasi-isomorphisms (and which we will call
the Segal category of perfect complexes in $T$). As $T$ is
Tannakian over $k$, there is a natural $k$-tensor functor
$k-Vect \longrightarrow T$, sending a $k$-vector space $V$
to the external product $V\otimes 1$, where $1 \in T$ is the unit
(as usual, for $E,F \in T$, $V\otimes E$ is determined by the adjunction formula
$Hom_{T}(V\otimes E,F)\simeq Hom_{k-Vect}(V,Hom_{T}(E,F))$).
Passing to the Segal categories of perfect complexes one finds
a morphism of $\otimes$-Segal categories
$$LParf(k-Vect)\simeq LParf(Hk) \longrightarrow
LParf(T),$$
making $LParf(T)$ into a rigid $k$-tensor Segal category.

Resuming the notations, one has a $k$-Tannakian category $T$
together with its gerbe of fiber functors $\mathcal{G}$, and
a rigid $k$-tensor Segal category $LParf(T)$ together
with its stack of fiber functors $FIB_{k}(LParf(T))$. The claim is
then the following.

\begin{prop}\label{ptan1}
The $k$-tensor Segal category $LParf(T)$ is Tannakian, and furthermore
one has a natural equivalence of stacks in
$LAff_{Sp^{\Sigma}}^{\sim,sffqc}/\mathbb{R}Spec\, Hk$
$$i(\mathcal{G})\simeq FIB_{k}(LParf(T)).$$
\end{prop}

The main corollary of proposition Prop. \ref{ptan1}
is the following, showing that our notion of
$k$-Tannakian Segal category is a reasonable generalization of
the usual notion.

\begin{cor}\label{ctan1}
The construction $T\mapsto LParf(T)$ defines a fully faithful morphism
from the $2$-Segal category of $k$-Tannakian categories to the
$2$-Segal category of $k$-Tannakian Segal categories.
\end{cor}

One can also go further and characterize the image
of the morphism $T \mapsto LParf(T)$. Indeed, if $T'$ is
any $k$-Tannakian Segal category, then $T'$ has a natural
$t$-structure defined as in \ref{tanconj} $(b-i)$ whose heart
will be denoted by $\mathcal{H}(T')$. The heart $\mathcal{H}(T')$
is always a $k$-Tannakian category, and there exists a natural
morphism of $k$-tensor Segal categories
$$LParf(\mathcal{H}(T')) \longrightarrow T'.$$
Then, $T'$ is of the form $LParf(T)$ for some $k$-Tannakian category
$T$ if and only if this last morphism is an equivalence (and then
furthermore one has $T\simeq \mathcal{H}(T')$). \\

In the general situation, but still when $k$ is a field,
for any $k$-Tannakian Segal category $T$, one has
its heart $\mathcal{H}(T)$, which is a $k$-Tannakian category,
together with the natural morphism of $k$-tensor Segal categories
$$j : LParf(\mathcal{H}(T)) \longrightarrow T.$$
It is not hard to show that after passing to $Ind$-objects,
$j$ has a right adjoint
$$q : Ind(T) \longrightarrow Ind(LParf(\mathcal{H}(T)))\simeq
LC(Ind(\mathcal{H}(T))),$$
where $LC(Ind(\mathcal{H}(T)))$ is the Segal category of complexes
in the abelian category of $Ind$-objects $Ind(\mathcal{H}(T))$.
Using the adjunction $(j,q)$ one should then construct
a commutative monoid $A:=q(1)$ in $LC(Ind(\mathcal{H}(T)))$, and
an equivalence between $Ind(T)$ and the Segal category $A-Mod$
of $A$-modules in $LC(Ind(\mathcal{H}(T)))$. Condition
$(2-ii)$ of conjecture \ref{tanconj} would then implies that
$A$ is cohomologically concentrated in non-negative degrees.
This construction
would give a structure theorem for $k$-Tannakian Segal categories
when $k$ is a field, that we state as a conjecture.

\begin{conj}\label{tanconj2}
Let $k$ be a field and $T$ be a $k$-Tannakian Segal category.
Then there exists a $k$-Tannakian category $\mathcal{H}(T)$,
and a non-negatively graded $E_{\infty}$-algebra $A$
in the $k$-tensor category
$Ind(\mathcal{H}(T))$, such that $T$ is equivalent to the
Segal category of $E_{\infty}$-modules over $A$ in $Ind(\mathcal{H}(T))$.
\end{conj}

One of the main consequence of conjecture \ref{tanconj2} would be that
any Tannakian Segal category over a field has a nice model
as a model category, and that they can all be constructed as Segal categories
of modules over an $E_{\infty}$-algebra in a usual Tannakian category.
Essentially, this means that Tannakian Segal categories over fields
are essentially the same thing as
$E_{\infty}$-algebras equipped with an action of an
affine group scheme $G$. Geometrically, the affine group
scheme $G$ is the fundamental group of
the associated Tannakian Segal gerbe, whereas
the $E_{\infty}$-algebra $A$ is a model for the homotopy type of
its universal covering (i.e. is its $E_{\infty}$-algebra of cohomology).
This last picture actually appears in practice, where one constructs
directly certain equivariant $E_{\infty}$-algebras
without using the notion of Tannakian Segal
categories (see for example \cite{ha,ka-pa-to,ol}). \\

\section{Homotopy types of algebraic varieties}

The cohomology groups  of algebraic varieties are endowed with
additional structures reflecting the algebraic nature of the space
(e.g. Hodge structures, Galois actions, crystal structures \dots).
More generally, one expects that algebraic varieties have
not only cohomology groups but full homotopy types also having additional structures extending
the one on cohomology. In this section I will present some
of my works around the very basic question:
\textit{How to define interesting homotopy types of algebraic varieties, and
what are their additional structures ?} The references for this
section are \cite{ka-pa-to,to2,to3}.

I have been considering seriously only the case of complex algebraic varieties.
Indeed, in the complex situation, algebraic varieties has an
underlying topological space of complex point (with the analytic topology)
and therefore a given homotopy type. Of course, as Hodge structures
does not exist directly on the integral cohomology groups but only after
tensoring with $\mathbb{C}$, one could not expect the topological homotopy type
of a complex variety to have a reasonable Hodge structure. One first problem
was therefore to extract from the topological homotopy type the part
that has to do with algebraic geometry and that one can endow
with an additional Hodge structure. One possible answer to this
is the solution to the schematization problem of Grothendieck presented below.
I will then describe how the schematization of a complex
smooth and projective variety can be endowed with
a certain action of $\mathbb{C}^{\times}$ playing the role of the
\emph{Hodge decomposition}. Finally, a Tannakian point of view
on this construction will be given.

\subsection{The schematization problem and one solution}

In \textit{pursuing stacks} Grothendieck considers what he calls
the \textit{schematization problem}. The questions he asked are not
very well defined and not very precise (e.g. he was considering
higher stacks without even having defined higher groupoids !),
which makes the understanding of the problem hazardous.
I personally understood the \textit{schematization
problem} in the following way (I refer to the appendix
of \cite{to2} for more details). \\

\begin{center} \textit{The schematization problem} \end{center}

\begin{enumerate}

\item For a ring $k$, there should exist a notion
of \textit{schematic homotopy types over $k$}. These are expected
to be \textit{$\infty$-stacks} (i.e. presheaves of $\infty$-groupoids, or equivalently
presheaves of homotopy types)
on the site of $k$-schemes satisfying certain algebraicity conditions.
It is expected for example that the Eilenberg-MacLane stack $K(\mathbb{G}_{a},n)$ is a schematic
homotopy type. It is also expected that the category of schematic homotopy types
is stable by certain standard constructions as fiber products, extensions \dots

\item For any topological space $X$, and any ring $k$,
there should exist a schematic homotopy type $(X\otimes k)^{sch}$, called
the \textit{schematization of $X$ over $k$}. The stack $(X\otimes k)^{sch}$
is required to be the universal schematic homotopy type receiving
a morphism from $X$ (considered as a constant presheaf of homotopy types and
therefore as a constant $\infty$-stack).

\item If $k=\mathbb{Q}$, simply connected schematic homotopy types over
$\mathbb{Q}$ are expected to be models for simply connected rational homotopy types.
Furthermore, for any simply connected space $X$,
the schematization $(X\otimes \mathbb{Q})^{sch}$
is expected to be a model for the rational homotopy type of the space
$X$ as defined by D. Quillen and D. Sullivan.

\end{enumerate}

In \cite{to2} I propose two solutions to the schematization problem. I will
present here the second one, which is for the moment only valid
when the base ring $k$ is a field. \\

In order to state the main definition, let us recall that
for any ring $k$ one has the site of affine $k$-schemes $Aff/k$
endowed the ffqc (faithfully flat and quasi-compact) topology.
The Segal category of stacks on the site $Aff/k$
will be denoted by $LAff_{k}^{\sim,ffqc}$, and its object will be called
\textit{stacks over $k$}. Note that the homotopy category
of $LAff_{k}^{\sim,ffqc}$ is simply the homotopy category of
simplicial presheaves on $Aff/k$. Let us also recall that
a stack over $k$, $F$, has a sheaf of connected component $\pi_{0}(F)$,
and sheaves of homotopy groups $\pi_{i}(F,*)$ well defined
for any global base point $* \longrightarrow F$.

\begin{df}\label{d8}
Let $k$ be a field.
A \emph{neutral schematic homotopy type over $k$} is a stack
$F \in LAff_{k}^{ffqc}$
which satisfies the following three conditions.
\begin{enumerate}
\item There exist a morphism of stacks $* \longrightarrow F$
inducing an isomorphism of sheaves $*\simeq \pi_{0}(F)$.

\item The homotopy sheaf $\pi_{1}(F,*)$ is representable by
an affine group scheme (for any choice of global base point).

\item For any $i> 1$, the homotopy sheaf $\pi_{i}(F,*)$
is representable by an affine and unipotent group scheme.
\end{enumerate}

A \emph{schematic homotopy type over $k$} is a stack over $k$ which after
a field base change $k'/k$ becomes a neutral schematic
homotopy type over $k'$.
\end{df}

A solution to the schematization problem of Grothendieck is given by the
following theorem.

\begin{thm}{(\cite[Thm. 3.3.4]{to2})}\label{t9}
Let $k$ be any field.
Let $X$ be any connected simplicial set, considered as
a constant simplicial presheaf and therefore as an object in
$Ho(k)$. There exist a schematic homotopy type $(X\otimes k)^{sch}$ over $k$
(automatically neutral),
and a morphism $X \longrightarrow (X\otimes k)^{sch}$
which is universal among morphisms towards schematic homotopy types
over $k$.
\end{thm}

For a connected simplicial set $X$ and $x\in X$ one can use
the universal property of the schematization in order to
prove the following fundamental properties (see \cite{to2}).

\begin{enumerate}
\item The sheaf $\pi_{1}((X\otimes k)^{sch},x)$ is naturally isomorphic
to $\pi_{1}(X,x)^{alg}$, the pro-algebraic hull of the discrete group $\pi_{1}(X,x)$ (relative to the
field $k$).

\item For any finite dimensional local system of $k$-vector spaces $L$ on $X$,
the morphism $X \longrightarrow (X\otimes k)^{sch}$ induces isomorphism in cohomology
with local coefficients
$$H^{*}((X\otimes k)^{sch},L) \simeq H^{*}(X,L)$$
(this makes sense because of property $(1)$ above;
$L$ can also be considered as a local system on
the stack $(X\otimes k)^{sch}$).

\item Let assume that $X$ is a finite simply connected
simplicial set. Then, one has natural isomorphisms
$$\pi_{i}((X\otimes k)^{sch},x)\simeq \pi_{i}(X,x)\otimes \mathbb{G}_{a} \quad if \quad car\, k=0,$$
$$\pi_{i}((X\otimes k)^{sch},x)\simeq \pi_{i}(X,x)\otimes \hat{\mathbb{Z}}_{p} \quad if \quad car\, k=p>0.$$

\item The functor $X \mapsto (X\otimes k)^{sch}$, when restricted to rational (resp. $p$-complete)
simply connected simplicial sets of finite type, is fully faithful.

\end{enumerate}

\begin{rmk}\label{ld8}
\emph{Schematic homotopy types as defined in Def. \ref{d8} are
very close to the notion of \emph{Tannakian Segal gerbes} introduced
in \S 5.3, and actually the two notions are
more or less equivalent. In the same way, the object $(X\otimes k)^{sch}$
of theorem \ref{t9} is conjecturally the dual Tannakian gerbe of some
Tannakian Segal category of perfect complexes of $X$.
These relations will be precised later
in \S 6.3.}
\end{rmk}

\subsection{Schematization and Hodge theory}

A complex algebraic variety $X$ has an underlying topological space
of complex point with the analytic topology $X^{top}$. It is well
known that the topology of $X^{top}$ is in general not easy to
describe in purely algebraic terms, as for example if $X$ is defined
over a number field $K$ the topology of $X^{top}$ depends
non trivially on the embedding of $K$ into $\mathbb{C}$.
However, the schematization $(X^{top}\otimes \mathbb{C})^{sch}$
can be described purely in algebraic terms and without referring
to the analytic topology (this is in fact an incarnation
of the Riemann-Hilbert correspondence). In fact,
as shown in \cite{ka-pa-to}, the stack $(X^{top}\otimes \mathbb{C})^{sch}$ has
an explicit algebraic model which uses algebraic de Rham complexes
of forms with coefficients in various flat bundles. This explicit description
in terms of differential forms allows one to use
Simpson's non-abelian Hodge correspondence in order to
endow the stack $(X^{top}\otimes \mathbb{C})^{sch}$ with a \textit{Hodge decomposition},
which reflects how the algebraic nature of the manifold $X$ interacts with his homotopy type.

\begin{thm}{(Katzarkov-Pantev-To\"en, \cite{ka-pa-to})}\label{t10}
Let $X$ be a smooth projective complex algebraic variety. There exists a natural
action of the discrete group $\mathbb{C}^{\times}$, called the
Hodge decomposition, on the
stack $(X^{top}\otimes \mathbb{C})^{sch}$ satisfying the following conditions.
\begin{enumerate}
\item The induced action on $H^{*}((X^{top}\otimes \mathbb{C})^{sch}),\mathcal{O})\simeq
H^{n}(X^{top},\mathbb{C})$ is compatible with the Hodge decomposition (i.e.
$\mathbb{C}^{\times}$ acts by weight $q$ on the factor $H^{p}(X,\Omega_{X}^{q})$).

\item The induced action of $\mathbb{C}^{\times}$ on $\pi_{1}((X^{top}\otimes \mathbb{C})^{sch}),x)\simeq
\pi_{1}(X,x)^{alg}$ is the one constructed by C. Simpson in \cite{si2}.

\item Assume that $X$ is simply connected. The induced action of $\mathbb{C}^{\times}$
on $\pi_{i}((X^{top}\otimes \mathbb{C})^{sch}))\simeq \pi_{i}(X)\otimes \mathbb{G}_{a}$
is compatible with the Hodge decomposition on the rational homotopy type
defined in \cite{mo}.

\end{enumerate}
\end{thm}

The above theorem gives a way to unify the Hodge decomposition on
the rational homotopy type of \cite{mo} with the Hodge decomposition
on the fundamental group of \cite{si2}. In a way, everything was already
contained in the non-abelian Hodge correspondence \cite{si2}, and
the new feature of theorem \ref{t10} is to give a homotopy theory interpretation of this
correspondence, based on the notion of schematic homotopy types and
the schematization functor. The notion of schematic homotopy types was apparently
the missing part in order to relate the various works on non-abelian
Hodge theory (e.g. \cite{si3,si4,ha,ka-pa-si}) to actual homotopy theory. \\

Theorem \ref{t10} possesses two important consequences. First of
all, it is not difficult to produce examples of finite CW complexes
$X$ such that the stack $(X\otimes \mathbb{C})^{sch}$ can not be endowed
with a $\mathbb{C}^{\times}$-action satisfying the conditions of
theorem \ref{t10}. In particular, this gives new examples of homotopy types
which are not realizable by projective manifolds, and the interesting new
feature is that obstructions to realizability lie in higher
homotopical invariants (precisely the action of the fundamental group
on the higher homotopy groups, see \cite{ka-pa-to} for details). Another important consequence
is the degeneracy of the \textit{Curtis spectral sequence}, starting from
the homology of $(X^{top}\otimes \mathbb{C})^{sch})$ with coefficients in the
universal reductive local system, and converging to its
homotopy groups (this is a generalization of
the degeneracy of the Bar spectral sequence). To be a bit more precise,
for any neutral and pointed schematic homotopy type $*\rightarrow F$, one can consider
its universal reductive covering $F^{0}$, corresponding to the
unipotent radical of $\pi_{1}(F,*)$. The stack $F^{0}$ is clearly a pointed and
connected affine stack in the sense of \cite{to2} and can therefore
by represented as $BG_{*}$, where $G_{*}$ is a simplicial affine group scheme
such that each $G_{n}$ is a free unipotent group scheme.
By considering the central lower
series of $G_{*}$ one constructs a tower of fibrations
$$\xymatrix{& \dots B(G_{*}/\mathcal{C}^{n}G_{*}) \ar[r] &
B(G_{*}/\mathcal{C}^{n-1}G_{*}) \ar[r] & \dots &
B(G_{*}/[G_{*},G_{*}])=B(G_{*})_{ab},}$$ whose limit is the stack
$F^{0}$. Corresponding to this tower of fibrations is a spectral
sequence in homotopy as explained in \cite{bou-kan}, which by definition
is the Curtis spectral sequence of $F$. The
$E_{1}^{p,q}$-term of this spectral sequence is
$\pi_{p}(B(\mathcal{C}^{q}G_{*}/\mathcal{C}^{q+1}G_{*})$, and its
abutment is $\pi_{*}(F^{0})$. Furthermore, one can see that
$\pi_{p}(B(\mathcal{C}^{q}G_{*}/\mathcal{C}^{q+1}G_{*}))$ only depends on the
graded group $\pi_{*}((G_{*})_{ab})$ which is the
homology of $F^{0}$, or in other words the
homology of $F$ with coefficients in the universal
reductive local system $\mathcal{L}^{red}$. This implies that the
$E_{1}$-term of the Curtis spectral sequence for $F$
is given by the free
Lie algebra over $H_{*>0}(F,\mathcal{L}^{red})$. The fact that
when $F=(X^{top}\otimes \mathbb{C})^{sch}$, for a smooth projective complex manifold
$X$, the Curtis spectral sequence degenerates at $E_{2}$ simply follows from
some \textit{weight property} of the action of $\mathbb{C}^{\times}$
on it. The conclusion is
the following result.

\begin{cor}{(Katzarkov-Pantev-To\"en, \cite{ka-pa-to})}\label{ct10}
\begin{enumerate}
\item For any pointed schematic homotopy type $F$, there exists
a Curtis spectral sequence $E_{r}^{p,q}$, whose
$E_{1}$-term is the
free Lie algebra over $H_{*>0}(F,\mathcal{L}^{red})$ and whose abutment is
$\pi_{*}(F)$.

\item Let $X$ by a (pointed) smooth projective complex manifold.
The Curtis spectral sequence of $(X^{top}\otimes \mathbb{C})^{sch}$
degenerates at $E_{2}$.
\end{enumerate}
\end{cor}

\begin{rmk}\label{r4}
\emph{Directly related
to this is the formality theorem stated in \cite{ka-pa-to} which generalizes the well
known fact that the rational homotopy type of
a smooth projective manifold is formal.}
\end{rmk}

The conclusion is that the existence of the Hodge decomposition of theorem \ref{t10}
has very strong consequences on the
schematization $(X^{top}\otimes \mathbb{C})^{sch}$ and this leads to
a striking fact. In the general situation, the schematization
$(X\otimes \mathbb{C})^{sch}$ of a space $X$ seems very hard (if not
impossible) to compute, as its homotopy invariants can be very far
from the original homotopy invariants of $X$ (see e.g. \cite[\S 3.4]{to2}).
However, the degeneracy of the Curtis spectral sequence and the formality theorem
imply that the schematization of a smooth projective manifold is much more
simple than for an arbitrary space, and that one can expect to
compute certain homotopy invariants which seem out of reach in the general situation.
Another result going in the same direction is the fact that
fundamental groups of Artin's neighborhood are \textit{algebraically good groups} , and
therefore that Zariski locally a smooth projective manifold has a very simple
schematization. Together with the Van-Kampen theorem this
again implies that the schematization of a smooth projective manifold looks
more simple than for a general space (see \cite{ka-pa-to2} for more
on this). \\

\subsection{The Tannakian point of view}

For a space $X$, the fundamental group of its schematization
$\pi_{1}((X\otimes k)^{sch},x)$ is naturally isomorphic
with the pro-algebraic completion of the discrete
group $\pi_{1}(X,x)$ over the base field $k$. Hence, the
functor $X \mapsto (X\otimes k)^{sch}$ is an
extension of the pro-algebraic completion functor.

It is very well known that the pro-algebraic
group $\pi_{1}(X,x)^{alg}$ has a Tannakian interpretation, as
the Tannaka dual of the category of local systems of finite dimensional
$k$-vector spaces on $X$ (see \cite[\S 10.24]{de2}). In the same way
the whole stack $(X\otimes k)^{sch}$ does have a Tannakian
interpretation in the Segal sense of \S 5,
at least conjecturally. Actually, it turned out that
this Tannakian point of view was originally the way
the schematization functor was introduced for the first time (in \cite{to3}).
This point of view has been also used in order to define
the Hodge decomposition, way before the explicit construction presented
in \cite{ka-pa-to} has been considered. \\

Let $X$ be a connected finite CW complex, $k$ be any ring (commutative with unit), and
$C(X,k)$ be the category of complexes of presheaves of $k$-modules on $X$.
In the category $C(X,k)$ there is a notion of \textit{local quasi-isomorphism}
(morphisms inducing quasi-isomorphisms on each stalks), and $C(X,k)$ can be
made into a model category for which equivalences are
local quasi-isomorphisms. Furthermore, the tensor product
of complexes makes $C(X,k)$ into a symmetric monoidal model
category. By  localization
one gets a $\otimes$-Segal category $LC(X,k)$. By definition,
the Segal category $LParf(X,k)$ is the full sub-Segal category of
$LC(X,k)$ consisting of rigid objects (they are exactly
the perfect complexes, i.e. are locally
on $X$ quasi-isomorphic to a constant complex of presheaves associated with a
bounded complex of projective $k$-modules  of finite type).

Pulling back from the point gives a morphism of
$\otimes$-Segal categories
$$LParf(k):=LParf(*,k) \longrightarrow LParf(X,k)$$
making $LParf(X,k)$ into a rigid $k$-tensor Segal
category. Furthermore, taking a base point $x \in X$, and considering
the pull-back along $x : * \longrightarrow X$ gives rise to
a $k$-tensor morphism
$$\omega_{x} : LParf(X,k) \longrightarrow LParf(k).$$
Clearly, $(LParf(X,k),\omega_{x})$ satisfies the
conditions of conjecture Conj. \ref{tanconj}, and therefore
should be a Tannakian Segal category. \\

When $k$ is a field, one expects the following
conjecture. Recall that one has full embedding
$$i : LAff_{k}^{\sim,ffqc} \longrightarrow
LAff_{Sp^{\Sigma}}^{\sim,sffqc}/\mathbb{R}Spec\, Hk.$$

\begin{conj}\label{conj1}
Let $k$ be a field and $(X\otimes k)^{sch}
\in LAff_{k}^{\sim,ffqc}$ be the stack
defined as in Thm. \ref{t9}.
\begin{enumerate}
\item There is a natural morphism
$$i(X\otimes k)^{sch} \longrightarrow FIB_{k}(LParf(X,k))$$
which is a $P$-equivalence (see Def. \ref{dtan7}). In particular, the natural
morphism of $k$-tensor Segal categories
$$LParf(X,k) \longrightarrow LParf(i(X\otimes k)^{sch})$$
is an equivalence.
\item
If $k$ is of characteristic $0$,
then the natural morphism
$$i(X\otimes k)^{sch} \longrightarrow FIB_{k}(LParf(X,k))$$
is an equivalence.
\end{enumerate}
\end{conj}

Conjecture \ref{conj1} essentially states that the stack
$(X\otimes k)^{sch}$ is the Tannakian dual
of the Tannakian Segal category $LParf(X,k)$
\footnote{In positive characteristic the object
$(X\otimes k)^{sch}$ is not exactly
the dual of $LParf(X,k)$ but still determines it. In fact,
as shown in \cite{ka-pa-to2} $(X\otimes k)^{sch}$ is of the form
$B\mathbb{R}Spec\, B$, for $B$ a co-simplicial Hopf algebra. Viewing
$B$ as a $H_{\infty}$-Hopf $Hk$-algebra, conjecture
\ref{conj1} predicts that
$B\mathbb{R}Spec\, B\simeq FIB_{k}(LParf(X,k))$.}. This is the
(conjectural) Tannakian interpretation of the schematization functor.
It is important to note that the object $FIB_{k}(LParf(X,k))$ is itself
defined for any commutative ring $k$, which gives (up to
conjecture \ref{conj1}) an extension of the schematization functor
over any base ring. \\

A first application of conjecture \ref{conj1} and the whole
Tannakian formalism, is an alternative construction of the Hodge
decomposition of \cite{ka-pa-to}, more in the style of
the construction of the Hodge structure on the fundamental group given
in \cite{si2}. Indeed, let
$X$ be a connected smooth and projective complex variety. One has
two $\mathbb{C}$-tensor Segal categories
$LParf(X^{top},\mathbb{C})$, and $LParf(X_{Dol},\mathcal{O})$. The
first one is the Tannakian Segal category of perfect complexes of
$\mathbb{C}$-vector spaces on the topological space of complex
point $X^{top}$. The latter is defined to be the Segal category of
complexes of quasi-coherent $\mathcal{O}$-modules on the stack
$X_{Dol}$ (i.e. complexes of quasi-coherent sheaves endowed with
an integrable Higgs field, see \cite{si2,si4,si5}) which are cohomologically
bounded and whose cohomology sheaves (which are Higgs coherent
sheaves) are semi-stable Higgs bundles of degree $0$ on $X$. The
equivalence of tensor dg-categories given in \cite[Lem. 2.2]{si2} can be
enhanced as
an equivalence of $\mathbb{C}$-tensor Segal categories
$$LParf(X^{top},\mathbb{C}) \simeq LParf(X_{Dol},\mathcal{O}).$$
Now, the discrete group $\mathbb{C}^{\times}$ acts on the stack
$X_{Dol}$, and therefore on the Tannakian Segal category
$LParf(X_{Dol},\mathcal{O})$. Using the equivalence above one gets
an action of $\mathbb{C}^{\times}$ on the Tannakian Segal category
$LParf(X^{top},\mathbb{C})$, and therefore on its Tannakian
dual $(X^{top}\otimes\mathbb{C})^{top}$. This construction was the
original idea of the construction of the Hodge decomposition
on the schematization, and have been abandoned by using
the more explicit approach taken in \cite{ka-pa-to}. I think however
that it is interesting to keep the Tannakian point of view, it helps
understanding things better, even if the Tannakian formalism for Segal
categories is still conjectural. \\

To finish this comparison let me mention the following conjecture,
establishing a relation between schematic homotopy types and
Tannakian Segal gerbes.

\begin{conj}\label{conj1'}
Let $k$ be a field and $F
\in LAff_{k}^{\sim,ffqc}$ be a schematic homotopy type in the
sense of Def \ref{d8}.
\begin{enumerate}

\item The $k$-tensor Segal category
$LParf(i(F))$ is a $k$-Tannakian Segal category, and therefore
$FIB_{k}(LParf(i(F)))$ is a $k$-Tannakian gerbe.

\item There is a natural morphism
$$i(F) \longrightarrow FIB_{k}(LParf(i(F)))$$
which is a $P$-equivalence (see Def. \ref{dtan7}). In particular, the natural
morphism
$i(F) \longrightarrow FIB_{k}(LParf(i(F)))$ is universal
among morphisms from $i(F)$ towards
$k$-Tannakian Segal gerbes.

\item
If $k$ is of characteristic $0$,
then the natural morphism
$$i(F) \longrightarrow FIB_{k}(LParf(i(F)))$$
is an equivalence. In particular $i(F)$ is
a $k$-Tannakian Segal gerbe dual to the $k$-Tannakian
Segal category $LParf(i(F))$. Furthermore, the morphism
$F \mapsto i(F)$ induces and equivalence between the Segal category of
schematic homotopy types and the Segal category of $k$-Tannakian
Segal gerbes.
\end{enumerate}
\end{conj}

Conjecture \ref{conj1'} essentially says that schematic homotopy types
are duals to Tannakian Segal categories, at least in the characteristic zero
case. In positive characteristic however, schematic homotopy types are only
approximations of Tannakian Segal gerbes, and the morphism
$F \mapsto i(F)$ does not provide an equivalence between the two notions
anymore. This last morphism is actually very similar to the
normalization functor $N : Ho(Alg_{k}^{\Delta}) \longrightarrow
Ho(E_{\infty}-Alg_{k})$, going from the
homotopy category of co-simplicial commutative $k$-algebras to the
homotopy category of $E_{\infty}$-$k$-algebras, which is
well known not to be an equivalence (not even fully faithful)
when $k$ is of positive characteristic (see for example
\cite{kr}). As co-simplicial algebras are strict forms of
$E_{\infty}$-algebras, schematic homotopy types are
strict form of Tannakian Segal gerbes. \\

\subsection{Other homotopy types in algebraic geometry}

Inspired by the case of complex projective varieties and
their Hodge decomposition, one can also define and
study various schematic homotopy types
related to other cohomological theories, as
the $l$-adic and crystalline theories. In practice these
schematic homotopy types can be constructed using
some formalism of equivariant co-simplicial algebras
(see \cite[\S 3.5]{to2} for an overview), but they
are conjecturally duals to natural Tannakian Segal categories in the
same way that $(X\otimes k)^{sch}$ is conjecturally the dual
to the Tannakian Segal category $LParf(X,k)$. As the schematization
$(X^{top}\otimes \mathbb{C})^{sch}$ is endowed with a Hodge decomposition,
these new homotopy types have additional structures as for example
action of the Galois group of the ground field, or $F$-isocrystal structures,
which are in general given by an action of some group. This action is
expected to capture interesting geometric and arithmetic information,
as for example a rational point of a variety will give
rise to a natural homotopy fixed point on the corresponding
homotopy type. Actually, the map sending a rational point to a homotopy
fixed point can reasonably considered as a non-abelian version of
the Abel Jacobi maps (see \cite[\S 3.5.3]{to2} for more details),
and its seems to me a very interesting invariant to study in the future.
But this is another story, and I will finish this m\'emoire
here \dots .

\newpage

\begin{appendix}

\section{Some thoughts about $n$-category theory (a letter to P. May)}

In this appendix, I reproduce a letter to Peter May, written during the fall 2001,
for the purpose of some kind of NSF proposal around the general theme \textit{higher category
theory and applications}. We wrote it together with L. Katzarkov, T. Pantev and
C. Simpson in order to explain our interests in higher categories and its applications to
algebraic geometry. I thought that it is very much related to
the main subject of this m\'emoire and that it could be of some use to
present it in this appendix. \\

\bigskip

Dear Peter,

\

\smallskip

In this short note we would like to share some of our thoughts on
higher categories. We will be mainly concerned with applications, and
so we will adopt a utilitarian and pragmatic point of view.  We will
try to explain why having a unified theory of $n$-categories is of
fundamental importance for the applications. Though not directly
related to the \textit{unification program}, we hope that these
considerations will be of some help for your proposal.

It seems to us that the fundamental question in the subject is not
about the correct definition of an $n$-category (always in the weak
sense, and always with $n \in [0,\infty]$).  In practice, at least for
our interest in the subject, we have found it more important to ask
the following question: What is the $(n+1)$-category of
$n$-categories ? Of course, this sounds like a joke, because being able
to define the $(n+1)$-category of $n$-categories requires a notion of
$(n+1)$-categories, and therefore of $n$-categories (this is maybe one
of the possible explanation of the word \textit{Pursuing} in
Grothendieck's \textit{Pursuing stacks}). However, we hope the
following few pages will clarify our point of view.

We begin by describing some of the examples of higher categories which
are of interest to us. After that we explain why a model category
approach to higher categories seems to be very well suited to treat
these examples, and also for the study of higher categories
themselves. Finally, we will discuss the relevance of a unified theory
for the purpose of applications.

In the following, let us suppose we have a good theory of
$n$-categories (again, $n$ could be $\infty$).  By this, we mean for
example that there exist a $1$-category of $n$-categories with strict
functors, denoted by $n-Cat$, as well as a well defined
$(n+1)$-category of $n$-categories with lax functors, denoted by
$n-\underline{Cat}$.  For two $n$-categories $A$ and $B$, the
$n$-category of lax functors between $A$ and $B$ will be denoted by
$\underline{Hom}(A,B)$.  By a \textit{good theory} we also mean that
all expected properties are satisfied.  We apologize for not being
very precise about which theory of higher categories we use, but in
fact every construction we will consider can be made rigorous by
employing the $n$-Segal categories of Tamsamani-Simpson (see \cite{ta,hi-si}). We also
apologize in advance for not paying the proper attention to the
set-theorical complications that arise (most of the time at least two
universes $\mathbb{U}\in \mathbb{V}$ have to be chosen). \\

\begin{center} \textit{Some examples of $n$-categories} \end{center}

Here are some examples of higher categories we are interested in. The
first three are $\infty$-categories for which the $i$-morphisms are
invertible (up to $(i+1)$-morphisms) as soon as $i>1$ (in other words,
the $Hom$'s are $\infty$-groupoids). This type of higher
categories is of fundamental importance to us.  \

\smallskip

\textbf{Simplicial categories:}
Simplicial categories are categories enriched over the category of
simplicial sets (i.e. their $Hom$ sets are endowed with the
additional structure of simplicial sets, such that compositions are
morphisms of simplicial sets).  These type of categories is very
important as they are models for a certain kind of
$\infty$-categories, and appear naturally in topology. More precisely,
if $T$ is a simplicial category, one can produce an $\infty$-category
$\Pi_{\infty}(T)$ by keeping the same set of objects, and replacing
the simplicial sets of morphisms in $T$ by their $\infty$-fundamental
groupoids\footnote{We assume such a construction exists for
our theory of higher categories}. In other words, for two objects $x$ and $y$,
$$
Hom_{\Pi_{\infty}(T)}(x,y):=\Pi_{\infty}(Hom_{T}(x,y)).
$$
By this construction, every simplicial category will be considered as
an object in $\infty-Cat$, or in $\infty-\underline{Cat}$.
Some basic examples of this kind are:

\begin{itemize}
\item For a $1$-category $C$, and $S\subset Mor(C)$ a sub-set of
morphisms, one can consider the simplicial localization of Dwyer and
Kan $L(C,S)$.  The construction $(C,S) \mapsto L(C,S)$ is the left
derived functor of the usual localization $(C,S) \mapsto S^{-1}C$.  If
one considers $L(C,S)$ as an $\infty$-category, then there exists a
lax functor $C \longrightarrow L(C,S)$, which is universal for lax
functors $f : C \longrightarrow A$, with $A$ an $\infty$-category, and
such that $f$ sends $S$ to equivalence in $A$.

If $M$ is a simplicial model category, and $Equiv$ is the
sub-set of equivalence in $M$, then $LM := L(M,Equiv)$ is
equivalent to the simplicial category of fibrant and cofibrant objects
in $M$. When $M$ is a non-simplicial model category, one can still
compute the simplicial sets of morphisms in $LM$ by using the mapping
spaces of $M$.
\item One defines $Top := LSSet$, the simplicial
localization of the model category of simplicial sets. As an
$\infty$-category, the objects of $Top$ are the fibrant
simplicial sets, the $1$-morphisms are the morphisms of simplicial
sets, the $2$-morphisms are the homotopies, the $3$-morphisms are the
homotopies between homotopies \ldots and so on. The $\infty$-category
$Top$ seems to us as fundamental as the $1$-category
$Set$.
\item For a ringed topos $(T,\mathcal{O}_{T})$, one can consider
$C(T)$, the category of complexes of sheaves of
$\mathcal{O}_{T}$-modules in $T$. Consider further the simplicial
localization $LC(T)$ of $C(T)$ with respect to the quasi-isomorphisms.
As before, the objects of $LC(T)$ are the fibrant complexes of sheaves
of $\mathcal{O}_{T}$-modules (for a reasonable model structure), the
$1$-morphisms are the morphisms of complexes, the $2$-morphisms are
the homotopies, the $3$-morphisms homotopies between homotopies \ldots
and so on.
\end{itemize}

\

\smallskip

\textbf{Stacks of complexes:}
For a $1$-category $C$, one can consider the $(n+1)$-category
$\underline{Hom}(C^{op},n-\underline{Cat})$. This is the
$(n+1)$-category of $n$-prestacks on $C$.  When $C$ is a Grothendieck
site, one can also consider the full sub-$(n+1)$-category $n-St(C)$ of
$\underline{Hom}(C^{op},n-\underline{Cat})$ consisting of
\textit{$n$-stacks}. These are the $n$-prestacks which satisfy the descent
conditions for hyper-coverings (due to the lack of space
we will not make this precise here).

For a scheme $X$, let $\mathcal{O}_{X}-Mod$ be the ringed topos of
sheaves of $\mathcal{O}$-modules on the big site on $X$ with the
faithfully flat and quasi-compact topology. Taking the simplicial
localization of each category of complexes $C(\mathcal{O}_{X}-Mod)$,
one obtains an $\infty$-prestack on the category of schemes
$$
\begin{array}{cccc} LC(-,\mathcal{O}) : & Sch^{op} & \longrightarrow
& n-\underline{Cat} \\ & X & \mapsto & LC(\mathcal{O}_{X}-Mod).
\end{array}
$$
The fact that this $\infty$-prestack is actually an $\infty$-stack was
proved by A.Hirschowitrz and C.Simpson (see \cite{hi-si}). This $\infty$-stack contains
some very important full $\infty$-sub-stacks, as for example
$\underline{Perf}$, the $\infty$-stack of perfect complexes of
$\mathcal{O}$-modules. The $\infty$-stack of perfect complexes plays a
crucial role in C.Simpson's non-abelian Hodge theory. For
example, if $X$ is a scheme over $\mathbb{C}$, one can associate to $X$
its de Rham shape $X_{DR}$ and Dolbeault shape $X_{Dol}$. The shapes
$X_{DR}$ and $X_{Dol}$ are
sheaves (and therefore $\infty$-stacks) on the big flat site of schemes.
The following two  $\infty$-categories
$$\underline{Hom}_{{\infty}-St(Sch)}(X_{DR},
\underline{Perf}), \qquad
\underline{Hom}_{{\infty}-St(Sch)}(X_{Dol},
\underline{Perf}),
$$
are of
particular interest (here $\underline{Hom}_{{\infty}-St(Sch)}$
denotes
the $\infty$-categories of morphisms in the $\infty$-category of
$\infty$-stack).  Actually, the main objects of study of the
non-abelian Hodge theory are the $\infty$-stacks
$$
\mathcal{HOM}_{{\infty}-St(Sch)}(X_{DR},\underline{Perf}),
\qquad
\mathcal{HOM}_{{\infty}-St(Sch)}(X_{Dol},\underline{Perf}),
$$
(here
$\mathcal{HOM}_{{\infty}-St(Sch)}$ denotes the $\infty$-stack of
morphisms, i.e.  the internal $Hom$ in $\infty-St(Sch)$).  The
$\infty$-stacks of this type possess rich additional structures, like
symmetric monoidal structures, linear structures, duality, rigidity
\ldots In order to understand these structures properly one has to
make sense of a  very advanced theory of
$\infty$-categories. Note however, that all the $\infty$-categories
appearing in these constructions come from simplicial categories.

Using these kind of constructions, C. Simpson has developed
non-abelian Hodge theory, and already several results have been proven
by him and his collaborators. Let us mention for example the
construction of the non-abelian Hodge filtration, Gauss-Manin connection and the proof
of its
regularity (see \cite{si4,si5}), the higher Kodaira-Spencer deformation classes
(see \cite{si7}), the non-abelian $(p,p)$-classes theorem (see \cite{ka-pa}),
a non-abelian analogue of the density of
the monodromy (see \cite{ka-pa-si1}), the notion of non-ablian mixed Hodge structure
(see \cite{ka-pa-si}),
and some new  restrictions on homotopy types of projective manifolds (see \cite{ka-pa-to}).
All these results where guessed
and proved using higher stack and/or higher category theory.

\

\smallskip

\textbf{Homotopy Galois theory:}
Let us go back to the $\infty$-category $Top$, obtained as the
simplicial localization of $SSet$ with respect to all
equivalences. For a nice enough space $X$
(e.g. a $CW$-complex), one can consider the $\infty$-category
$Top(X)$ of
locally constant stacks of $\infty$-groupoids on $X$. This is an
$\infty$-version of the category of locally
constant sheaves of sets on $X$ and has many analogous properties.  In
particular, it was shown by B. To\"en and G. Vezzosi (see \cite{to1,to-ve1}) that
the $\infty$-category $Top(X)$ can be used to reconstruct the whole
homotopy type of $X$\footnote{This statement was first mentioned by
A.Grothendieck in one of his letters to L.Breen.} much in the same way
as the fundamental groupoid
$\Pi_{1}(X)$ of $X$ can be reconstructed from  the
category of locally constant sheaves of sets on $X$.  Moreover, there is an explicit formula:
$$
\Pi_{\infty}(X)\simeq
\underline{Hom}^{geom}(Top(X),Top),
$$
where the
right hand side is the full sub-$\infty$-category of
$\underline{Hom}(Top(X),Top)$ of \textit{geometric} (i.e.
which are exact and possess a right adjoint)
lax functors. Combined with the observation that there exists a
natural equivalence $Top(X)\simeq
\underline{Hom}(\Pi_{\infty}(X),Top)$, the above formula can be
rewritten  as
$$
\Pi_{\infty}(X)\simeq
\underline{Hom}^{geom}(\underline{Hom}(\Pi_{\infty}(X),Top),
Top).
$$
In
this last formula, it is interesting to note that the
$\infty$-category $Top$ plays the role of a \textit{dualizing
object}. This observation is the starting point of a research program
on very general higher Tannaka dualities, for which the dualizing
object is replaced by the $\infty$-stack of perfect complexes (see \cite{to3}). Once
again, one should note the all the $\infty$-categories involved in
these considerations are associated to some simplicial categories, and are endowed with
various additional structures (linear, monoidal, tensorial \dots).

The reconstruction result quoted above is also the starting point of a
generalization of the categorical Galois theory (as developed by
A.Grothendieck) to the setting of $\infty$-categories. Such a theory
may find applications in homotopy theory (e.g. a new point of view on
pro-finite, pro-nilpotent \ldots localizations), as well as in
algebraic geometry (e.g. in etale homotopy type of schemes).

\

\smallskip

\textbf{Monoidal $n$-categories:}
We have already mentioned that there are examples of higher
categories possessing interesting extra structures. The symmetric
monoidal structures form an important type of such structures and
deserve special attention. One can
try to make sense of
the theory of symmetric monoidal $n$-categories in the following way.

Let $\Gamma$ be the category of pointed finite sets. Define the
$(n+1)$-category of symmetric monoidal $n$-categories $n-Sym$ to be
the full sub-$(n+1)$-category of
$\underline{Hom}(\Gamma,n-\underline{Cat})$ consisting of
lax functors
$F : \Gamma \longrightarrow n-\underline{Cat}$ which satisfy the
\textit{Segal conditions} (i.e. the same conditions satisfied by
$\Gamma$-spaces).

An important class of symmetric monoidal $\infty$-categories comprises
the symmetric monoidal model categories. In particular, if $M$  is such a model
category, then the associated $\infty$-category $LM$ possesses a natural
symmetric monoidal structure. Using this, one immediately checks that our
$\infty$-stacks $LC(-,\mathcal{O})$, $\underline{Perf}$ \ldots
of complexes are in fact symmetric monoidal $\infty$-stacks. For example, this
monoidal structure may be of some importance for the new
higher stack interpretation of Grothendieck duality obtained  recently
by A.Hirschowitz.

The higher monoidal categories may be applicable to many other
situations. For example they provide
a new point of view on $E_{\infty}$-ring structures. As
this is a kind of strictification result, we will discuss this in the
next section (see \textit{Strictification and monoidal structures}).

\

\smallskip

\textbf{Extended $TQFT$:} For $n\in [0,\infty]$ consider  the $n$-category $n-Cob$ of
cobordisms in dimension less than $n$. We do not know of any rigourous
construction of this category, but intuitively its objects are
oriented $0$-dimensional
compact varieties, the $1$-morphisms are the oriented $1$-dimensional compact
varieties with boundary, the $2$-morphisms are oriented compact surfaces
with corners  and so on. The $n$-category $n-Cob$ has a symmetric
monoidal structure induced by disjoint union of varieties.

Following J.Baez and J.Dolan, one can formulate the general purpose of extended
topological quantum field theories as the study of (higher) categories
of representations of $n-Cob$. Classically, one introduces the
$n$-category $n-Vect$, of $n$-vector spaces. Heuristically it is defined
by induction on $n$. First, one defines $1-Vect$ to be the
$1$-category of finite
dimensional vector spaces (over some base field). This category
possesses direct sums and a symmetric monoidal structure (the tensor
product) which make it into a \textit{rig}\footnote{
A rig $n$-category is an $n$-category equipped with two different symmetric
monoidal structures ($\oplus$ and $\otimes$),
such that the second one is distributive over the first one. The words
'rig' refers to 'ring' without 'n' (i.e. a ring without
a substraction)}  $1$-category.  Inductively,
if the rig $(n-1)$-category $(n-1)-Vect$ is defined,
$n-Vect$ will be
the $n$-category of $(n-1)$-categories which are free modules of
finite type over the rig $(n-1)$-category $(n-1)-Vect$. In
contrast with $n-Cob$, there is a rigorous construction (due to
B.Toen)  of
$n-Vect$
which utilizes  Tamsamani's $n$-categories.  In particular, $n-Vect$ is constructed as a rig $n$-category,
and hence as a symmetric monoidal $n$-category (if one forget the first monoidal structure).

By definition, the $n$-category of $n$-dimensional topological quantum
field theories is then given by
$$
n-TQFT :=
\underline{Hom}^{\otimes}(n-Cob,n-Vect).
$$
Unlike
previous examples, the higher categories appearing in $TQFT$ do not
come from simplicial categories. \\

\begin{center} \textit{Model categories and $n$-categories of lax functors} \end{center}

The above examples of higher categories are based in an essential way
on  the
existence of $n$-categories of lax functors, i.e the
existence of the $(n+1)$-category $n-\underline{Cat}$.
In this section we would like to
present a general approach due to C.Simpson, which allows one to construct
$n$-categories of lax functors (see \cite{hi-si,pe}). This  approach
uses  model categories, which seems to be significant in at least two
ways.   Firstly, it is directly related to your proposal  to use model
categories to compare different theories
of $n$-categories. Secondly, it gives a way to state and prove various
strictification results, which are of fundamental importance for many
applications. Actually, these strictification results may also be
directly relevant to the  comparison of different theories of $n$-categories.

\

\smallskip

\textbf{Model categories of $n$-categories:} The general idea is to define a model structure on the category
$n-Cat$ for which the weak equivalences are precisely the equivalences of
$n$-categories, and which is \textit{internal} (i.e. $n-Cat$ is a
symmetric monoidal model category for the monoidal structure given by
direct products). In other words, we require the existence of internal
$Hom$'s together with a compatibility condition with the model
structure. Unfortunately, it seems that such a model structure can not
exist directly on $n-Cat$\footnote{At least, the known or expected model
structures producing theories of $n$-categories do not seem to be
internal.}. For example, the category $SimpCat$, of simplicial
categories (which are models for $\infty$-categories for which
the $\infty$-categories or morphisms are $\infty$-groupoids) possesses a
model structure defined by B.Dwyer and D.Kan, which has the correct weak
equivalence but which is not internal (because the product of two
cofibrant objects is not cofibrant anymore). It seems  that any
approach to $n$-categories which is of operadic nature might
have this problem.

In order to circumvent this difficulty, C.Simpson has introduced a
notion of $n$-precategory. With his definition
the Tamsamani $n$-categories (or more generally the $n$-Segal
categories) can be viewed as $n$-precategories satisfying some special
conditions. The relation between $n$-precategories and $n$-categories
is very similar to the relation between prespectra and spectra,
presheaves and sheaves, prestacks and stacks, pre-$\Gamma$-spaces and
$\Gamma$-spaces, etc.. Remarkably, C.Simpson showed that the category
$n-PrCat$,
of $n$-precategories possesses an internal model structure.

From now on let  $n-PrCat$ denote the
model category of $n$-Segal precategories (here $n$ must be finite),
and let $\underline{Hom}_{n-PrCat}$ denote its internal
$Hom$.   We are not
going to recall its definition, but let us mention that objects of
$n-PrCat$ can be thought off as systems of generators and relations for
$n$-Segal categories, which in turn are models for
$\infty$-categories whose $\infty$-category of $n$-morphisms are
$\infty$-groupoids (i.e. all $i$-morphisms are invertible up to
$(i+1)$-morphisms for $i>n$).

An important property of $n-PrCat$ is that its fibrant objects are
all $n$-Segal categories (the converse is not true), and that every object
is cofibrant.  Let $n-SeCat$ be the full sub-category of
$n-PrCat$
consisting of $n$-Segal categories.  Since the model category
$n-PrCat$
is internal one can define the $(n+1)$-Segal category of $n$-Segal
categories $n-\underline{SeCat}$, by taking fibrant objects in
$n-PrCat$, and by forming the $Hom$'s between them with values in
$n-PrCat$.  Therefore, $n-\underline{SeCat}$ is a category
with values
in $n-SeCat$, which can be seen in an obvious way as a $(n+1)$-Segal
category. In other words, to compute the correct $n$-Segal categories
of functors between two $n$-Segal (pre)categories $A$ and $B$, one has
to consider $\underline{Hom}_{n-PrCat}(A,B')$, where $B'$
is a fibrant
 model for $B$. In conformity with tradition one writes
$\mathbb{R}\underline{Hom}_{n-PrCat}(A,B)$ for
$\underline{Hom}_{n-PrCat}(A,B')$.

This construction is actually more general. For any model category
$M$, which is enriched over $n-PrCat$ (as a model category), one can
define a $(n+1)$-Segal category $Int(M)$, whose objects are fibrant
and cofibrant objects in $M$ together with its $Hom$'s with values in
$n-SeCat$ coming from the enrichment.

\

\smallskip

\textbf{Strictification:}
Let $I$ be a $(n+1)$-precategory and let $M$ be a cofibrantly generated
model category enriched over $n-PrCat$. As the category
$n-PrCat$  is
acting on $M$, one can define the notion of a representation of $I$ in
$M$, or $I$-modules in $M$. These $I$-modules form a category $M^{I}$.
Since $M$ is
cofibrantly generated, it is likely that $M^{I}$ is again a
cofibrantly generated  model category for the levelwise model structure
(i.e. fibrations and equivalence are defined levelwise). Furthermore,
the model category $M^{I}$ is naturally enriched over $n-PrCat$,
and one
can consider the associated $(n+1)$-Segal category $Int(M^{I})$.

On the other hand, one can first consider $Int(M)$, and then the
$(n+1)$-Segal category of functors
$\mathbb{R}\underline{Hom}_{n-PrCat}(I,Int(M))$.  The expected
strictification theorem would be the existence of a natural
equivalence of $(n+1)$-Segal categories
$$
\mathbb{R}\underline{Hom}_{n-PrCat}(I,Int(M))
\simeq Int(M^{I}).
$$

In applications, this theorem is very fundamental, and for example is
of essential use for the few results mentioned in the first part of
this note. Some essential special cases of this theorem are proven by
A.Hirschowitz and C.Simpson (see \cite{hi-si}), but the general case is still a
conjecture. It is important to stress that this particular case of
the strictification theorem is used in the proofs of the few results
mentioned above.

One of the really important consequences of the strictification theorem is the
Yoneda lemma. Indeed, for any $n$-Segal category $A$  C.Simpson has
defined a Yoneda embedding
$$
h : A \longrightarrow
\mathbb{R}\underline{Hom}_{n-PrCat}(A^{op},(n-1)-SeCat).
$$
By definition the
essential image of this embedding consists of \textit{representable
functors}. A way to define this essential image
would be to show that every $n$-Segal
category is equivalent to a strict category with values in
$(n-1)-PrCat$. For such a category, there is an obvious morphism of
$n$-precategories $A\times A^{op} \longrightarrow (n-1)-PrCat$, which
induces the functor $h$ by adjunction. However, the fact that $h$ is
fully faithful is still a conjecture for $n>1$ (it is proved
for $1$-Segal categories, see \cite{si1})).

As $(n-1)-SeCat$ is by definition equivalent to
$Int((n-1)-PrCat)$,
one may actually find a fully faithful functor
$$
h : A \longrightarrow Int((n-1)-PrCat^{A^{op}}),
$$
which shows that
every $n$-Segal category $A$ should embeds into a $n$-Segal category
of the form $Int(M)$, for $M$ an $(n-1)-PrCat$-enriched model
category.

A direct consequence of this fact is the following method for
\textit{computing} the $n$-Segal categories of lax functors. Let $A$ and
$B$ be two $n$-Segal categories, and let us consider the model
category $(n-1)-PrCat^{A\times B^{op}}$. Then, the $n$-Segal category
$\mathbb{R}\underline{Hom}_{n-PrCat}(A,B)$ can now be
identified with
the full sub-$n$-Segal category of $Int((n-1)-PrCat^{A\times
B^{op}})$, consisting of functors $F : A\times B^{op} \longrightarrow
(n-1)-PrCat$ such that each evaluation at an object $a$ in $A$, $F_{a}
: B^{op} \longrightarrow (n-1)-PrCat$ is equivalent to a
\textit{representable functor}. This gives a way to systematically
reduce the computation
to the case of $n$-Segal categories of the form $Int(M)$, which
is again a very powerful tool for applications.

Unfortunately, it seems unlikely that the general approach of
C.Simpson could be generalized or imitated for other theories of
$n$-categories. Indeed, we have already mentioned that the model
category of simplicial categories is not internal, and as far as we
know the expected model structures for other theories may have the
same problem. But of course, we do not have a proof that these model
structures can not exist, and this is just a general feeling.

\

\smallskip

\textbf{Strictification and monoidal structures:}
To finish this section, let us mention another conjectural
strictification result related to $E_{\infty}$-ring structures.

For this, let $M$ be a symmetric monoidal closed model category, which
is assumed to satisfy some nice but reasonable properties (e.g.  is
cofibrantly generated and left proper). One can make sense of
$E_{\infty}$-algebras in $M$. These form a category denoted by
$E_{\infty}-Alg(M)$. It was shown by M.Spitzweck that the category
$E_{\infty}-Alg(M)$ is almost a model category (for its natural model
structure), and actually is a model category in many cases. Whatever
it is, there is a natural notion of a weak equivalence in
$E_{\infty}-Alg(M)$, and so one can consider the $\infty$-category
$LE_{\infty}-Alg(M)$, obtained by applying the simplicial localization
of Dwyer and Kan to the set of weak equivalences.

On the other hand, one can consider first $LM$, which is a symmetric
monoidal $\infty$-category, and the $\infty$-category
$$
Comm(LM):=\underline{Hom}^{\otimes}(FS,LM)
$$
 of commutative
(unital and associative, as usual) monoids in it. Here $FS$ is the
symmetric monoidal category of finite sets (where the monoidal
structure is given by disjoint union).

It was shown by B. To\"en (see \cite{to3}) that
there exists a natural functor
$$
Comm(LM) \longrightarrow LE_{\infty}-Alg(M),
$$
which is conjectured
to be an equivalence of $\infty$-categories. This conjecture would
identified $E_{\infty}$-algebra as commutative monoids in certain
$\infty$-categories (this point of view was also considered by
T.Leinster). \\

\begin{center} \textit{The need of a unification} \end{center}

What we wanted to stress out in the previous discussions is that for
purposes of applications, a good theory of $n$-categories requires in a
very essential way the existence of $n$-categories of lax functors, as
well as some strictification results. This
is the reason why we have been very much interested in $n$-Segal
categories instead of any other theories, for which we know, or at
least we expect, the required properties to hold.  Of course,
there is a price to pay for using $n$-Segal categories, and this is
what we would like to discuss in this last part.

When we look at T.Leinster's list of definition of $n$-categories,
the definitions seem to belong to two different classes. There are
definitions of $n$-categories for which the composition of morphisms is
well defined but associativity only holds up to coherent homotopies,
and definitions for which composition is defined only up to an
equivalence (i.e. composition is not defined in the conventional
sense).  Simplicial categories,
and all operadic definitions seem to belong to the first class, and
Tamsamani-Simpson definitions to the second one. For the purpose of
application and concrete manipulation each of these classes has its
own advantages.

\begin{itemize}
\item
The big advantage of the first type of definitions is that composition
being well defined allows a more easy manipulation of $n$-categories
themselves, and is also closer to the usual intuition of
categories. Also, they can be easily generalized to general contexts,
for example to deal with linear higher categories, or more exotic
enrichment. For example, in your talk at Morelia you are already
dealing with a general base category.  Such notions of linear, or
enriched higher categories are much more unclear for the second type
of definitions, and it may even be that there are no reasonable
analogs.

On the other side, already the example of simplicial categories seem to show that
defining higher categories of lax functors is quite difficult in this context.
It also shows that even if one can hope to define these categories of lax
functors, they could be difficult to use in practice.

\item
We have already mentioned the advantage of the second type of
definitions, which is the existence of a good theory of internal
$Hom$, and already explained this is of fundamental importance, even
to define the objects we would like to study. The main problem one encounter using these
kind definitions is the lack of computability. For example, when one wants to compute
the $n$-category of functors $\underline{Hom}(A,B)$, between two $n$-Segal categories, one
need first to consider a fibrant replacement $A'$ of $A$. Usually, this fibrant replacement
is highly non-explicit, and therefore very \textit{un-computable}. This fact
makes the standard categorical techniques (e.g. the Yoneda embedding) difficult to use
in this context.

\end{itemize}

We have already mentioned that various higher categories (or higher stacks) we have encountered in
applications are endowed with additional structure. For example, fur the purpose of
higher Tannakian duality and non-abelian Hodge theory the stack of perfect complexes
$\underline{Perf}$ has to be considered as a stack of \textit{tensor} $\infty$-categories. Therefore, the
first class of definitions of higher categories could be useful for us to give a sense to these
additional structures (as the linearity for example). On the other hand, internal $Hom$'s and
strictifications results can not be avoided to manipulate, define and even prove things. This makes
the second class of definitions of higher categories difficult to avoid for the purpose of our applications.
The same kind of remark can be made concerning the monoidal structure. We have already mentioned
one approach to monoidal structures in this letter, which is very much suited when one use
the second class of definitions. There exist also another approach, which consist of considering
$E_{\infty}$-algebras in some model categories of higher categories, and which is more in the style
of the second class of definitions. Again, the two approaches have their own advantages and dis-advantages, and
probably both are needed to really prove non-trivial theorems.

This situation with the stack $\underline{Perf}$ presented here is only one particular example, and one can find
many other situations where the two classes of definitions of higher categories seem to be
needed. Therefore, for us, the unification of these two classes of definitions
would be a great progress as far applicability is concerned. It would allow
one to choose the most suited of the two models for each particular
situation. Already for the theories of simplicial categories and
$1$-Segal categories, which are know to be equivalent, this principle
is highly used, and without such an equivalence many results presented in this letter would
have been unreachable to us. There is no doubt that a unified
theory will make accessible some expected results in algebraic
geometry concerned with higher categorical structures.\\

\bigskip

These are the reasons why we are very much interested in participating in
a unification program of all the different theories of higher
categories, and this is why we believe that such a unification
may also be useful to many other end-users of higher
categories.

\

All the very best,

\

\bigskip

L.Katzarkov, T.Pantev, C.Simpson, B.Toen.

\bigskip

\section{Comparing various homotopy theories}

As stated in section 1, there exist many different theory in order to do homotopy theory.
Motivated by P. May's project to unify the various point of views
on higher category theory, I have been interested in comparing four of these theories which
I was interested in. They are the theories of $S$-categories, Segal categories,
complete Segal spaces and quasi-categories. They are supposed to be all equivalent to each others, and
in this appendix I present comparison adjunctions which are conjecturally
Quillen equivalences. \\

\begin{center} \textit{Notations} \end{center}

The category of Segal precategories, $PrCat$, is the category of
functors
$$A : \Delta^{op} \longrightarrow SSet$$
such that $A_{0}$ is a discrete simplicial set. It is endowed with
the model structure described in \cite{hi-si,pe}.

The model category of quasi-categories, $QCat$, is the category of
simplicial sets, endowed with the model structure defined by
Joyal.

The model category of complete Segal spaces, $CSS$, is the
category of functors
$$X : \Delta^{op} \longrightarrow SSet$$
endowed with the model structure defined in \cite{re}.

Finally, $S-Cat$ is the category of $S$-categories (i.e.
simplicially enriched categories). It is endowed with the model
structure described in \cite[XII-48]{dw-hir-ka}\footnote{
The proof given there is not correct, as the generating trivial
cofibrations are not even equivalences. However, it seems that the
model structure still exists.}.

\begin{center} \textit{Segal categories and quasi-categories} \end{center}

One defines a functor
$$\phi : PrCat \longrightarrow QCat$$
by the following way. For $A \in PrCat$, the simplicial
set $\phi(A)$ is
$$\begin{array}{cccc}
\phi(A) : & \Delta^{op} & \longrightarrow & Set \\
 & [n] & \mapsto & (A_{n})_{0}.
\end{array}$$

Let $h_{n}$ be the object in $PrCat$ which is represented by $[n] \in \Delta$.
This defines a cosimplicial object $[n] \mapsto h_{n}$, in $PrCat$.
One has
$$\phi(A) \simeq Hom(h_{-},A).$$

This shows that the functor $\phi$ possesses a left adjoint
$$\psi : QCat \longrightarrow PrCat.$$

More precisely, for $X \in QCat$, $\psi(X)$ is the co-equalizer of the diagram in $PrCat$
$$\coprod_{[m] \rightarrow [p] \in \Delta}X_{p}\times h_{m}
\rightrightarrows \coprod_{[n] \in \Delta}X_{n}\times h_{n}.$$

This defines an adjunction
$$\psi : QCat \longrightarrow PrCat \qquad QCat \longleftarrow PrCat : \phi.$$

\begin{conj}\label{conj2}
The previous adjunction is a Quillen equivalence.
\end{conj}

\begin{center} \textit{Segal categories and complete Segal spaces} \end{center}

For any $[n] \in \Delta$ one denotes by $\overline{I}(n)$ the category
with $n+1$ objects and a unique isomorphism between them. Functors
$\overline{I}(n) \longrightarrow C$ are then is one-to-one correspondence
with strings of $n$ composable isomorphisms in $C$.

Considering categories as object in $PrCat$ (via their nerves), one gets
a co-simplicial object
$$\begin{array}{ccc}
\Delta^{op} & \longrightarrow & PrCat \\
 n & \mapsto & \overline{I}(n).
\end{array}$$

This functor extends in a standard way to a functor
$$\Pi : SSet \longrightarrow PrCat.$$

We define a simplicial structure on $PrCat$ via the functor $\Pi$. For $A$ and
$B$ objects in $PrCat$, the simplicial set of morphisms $Hom_{s}(A,B)$
is given by the following formula
$$Hom_{s}(A,B)_{n}:=Hom(A\times \overline{I}(n),B).$$
Tensors and co-tensors are define by
$$X\otimes A := \Pi(X)\times A \qquad A^{X}:=Hom_{s}(\Pi(X),A),$$
for any $A \in PrCat$ and $X \in SSet$.
The category $PrCat$ becomes  a simplicial model category.

We define a functor
$$\phi : PrCat \longrightarrow CSS$$
sending $A$ to the bi-simplicial set
$$\begin{array}{cccc}
\phi(A) : & \Delta^{op} & \longrightarrow & SSet \\
& [n] & \mapsto & Hom_{s}(h_{n},A).
\end{array}$$
Here, $h_{n}$ is again the object of $PrCat$ represented by $[n] \in \Delta$. \\

The functor $\phi$ possesses a left adjoint $\psi$ defined as follows. For $X \in CSS$,
the object $\psi(X)$ is the co-equalizer of the diagram
$$\coprod_{[m] \rightarrow [p] \in \Delta}X_{p}\otimes h_{m} \rightrightarrows
\coprod_{[n] \in \Delta}X_{n}\otimes h_{n}.$$

This defines an adjunction
$$\psi : CSS \longrightarrow PrCat \qquad CSS \longleftarrow  PrCat : \phi.$$

\begin{conj}\label{conj3}
The previous adjunction is a Quillen equivalence.
\end{conj}

\begin{center} \textit{Complete Segal spaces and $S$-categories} \end{center}

The category $S-Cat$ being a model category, one can
find a cofibrant resolution functor $\Gamma^{*}$ in the
sense of \cite[\S 5.2]{ho}.

Let $I(n)$ be the category with $n+1$ objects and a unique morphism between them.
Functors $I(n) \longrightarrow C$ are in one-to-one correspondence
with strings of $n$ composable morphisms in $C$. \\

We define a functor
$$\phi : S-Cat \longrightarrow CSS$$
as follows. For $T \in S-Cat$, the bi-simplicial set $\phi(T)$ is defined by
$$\begin{array}{cccc}
\phi(T) : & \Delta^{op} & \longrightarrow & SSet \\
 & [n] & \mapsto & Hom(\Gamma^{*}(I(n)),T).
\end{array}$$

This functor possesses a left adjoint $\psi : CSS \longrightarrow S-Cat$.
For $X \in CSS$, the $S$-category $\psi(X)$ is the co-equalizer of the
diagram
$$\coprod_{[m] \rightarrow [p] \in \Delta}X_{p}\otimes_{\Gamma} I(m) \rightrightarrows
\coprod_{[n] \in \Delta}X_{n}\otimes_{\Gamma} I(n).$$
Here, $\otimes_{\Gamma}$ is the tensor product over $\Gamma$. It is defined
for a simplicial set $Y$ and an $S$-category $T$ by the co-equalizer of the
following diagram
$$\coprod_{[m] \rightarrow [p] \in \Delta}\coprod_{Y_{p}} \Gamma^{m}(T) \rightrightarrows
\coprod_{[n] \in \Delta}\coprod_{Y_{n}} \Gamma^{n}(T).$$

This defines an adjunction
$$\psi : CSS \longrightarrow S-Cat \qquad CSS \longleftarrow  S-Cat : \phi.$$

\begin{conj}\label{conj4}
The previous adjunction is a Quillen equivalence.
\end{conj}

\end{appendix}

\end{document}